\documentclass[11pt,a4paper]{article}
\usepackage{amsmath,amssymb,amsthm,mathtools}
\usepackage{geometry}
\usepackage{enumitem}
\usepackage{xcolor}
\usepackage{hyperref}
\hypersetup{hidelinks,pdftitle={Minimax Growth and Morse Indices for the Critical Choquard Equation on the Sphere}}
\numberwithin{equation}{section}
\newtheorem{theorem}{Theorem}[section]
\newtheorem{proposition}[theorem]{Proposition}
\newtheorem{lemma}[theorem]{Lemma}
\newtheorem{corollary}[theorem]{Corollary}
\newtheorem{remark}[theorem]{Remark}
\theoremstyle{definition}

\title{Minimax Growth and Morse Indices\\
	for the Critical Choquard Equation on the Sphere}
\author{Wenjing Chen\footnote{E-mail addresses are \texttt{wjchen@swu.edu.cn} (Wenjing Chen)
		and \texttt{zhangle1430@163.com} (Le Zhang).}\hspace{1em}and Le Zhang\\[6pt]
	School of Mathematics and Statistics, Southwest University,\\
	Chongqing 400715, P. R. China}
\date{}

\begin{document}
	\maketitle
	
	\begin{abstract}
		We study the conformally critical Choquard equation on $\mathbb S^n$,
		\[
		-\Delta_{\mathbb S^n} v(\xi) + \frac{n(n-2)}{4} v(\xi)
		= \left(\int_{\mathbb S^n} \frac{|v(\eta)|^p}{|\xi-\eta|^\mu} dV_\eta\right) |v(\xi)|^{p-2} v(\xi),
		\]
		where $n\ge 3$, $0<\mu<n$, and $p=(2n-\mu)/(n-2)$ denotes the conformally critical Choquard exponent.
		For the standard action of \(G=O(k)\times O(\ell)\) with \(k,\ell\ge 2\) and \(k+\ell=n+1\), we use the cohomological index and orbit interval partition to build attained minimax levels with sharp growth \(c_j\asymp j^{2p/(p-1)}\). All \(G\)-invariant solutions at \(c_j\) change sign for \(j\ge 2\). If \(p\ge 2\), positivity of the Riesz form on nodal truncations gives \(N_G(v)\le m_G(v)\), and \(N_G(v)\le j\) whenever \(E(v)=c_j\). Here \(N_G\) counts nodal domains, and \(m_G,\nu_G\) are the invariant Morse index and nullity, respectively. In particular, \(c_2\) is the least energy of \(G\)-invariant sign-changing solutions, with such solutions having two nodal domains. Under the additional assumption that every critical point of the reduced functional \(\Psi\) at level \(c_j\) is isolated in the full critical set of \(\Psi\), the critical group theory gives \(m_G\leq j\leq m_G+\nu_G\) for a selected critical solution. For \(p\ge 2\), only finitely many angular blocks contribute to the full index and nullity, and every nonconstant \(G\)-invariant solution has transverse negative directions. When \(k=\ell\) and \(\mu>0\) is sufficiently small, we compute the radial index and nullity in the \((1,0),(0,1),(1,1)\) blocks of solutions from minimizing the reduced energy over \(G\)-invariant functions odd under factor interchange.
	\end{abstract}
	
	
	
	\section{Introduction}
	\noindent
	
    In this work, we consider the critical Choquard equation
	\begin{equation}
		Av=J_\mu(|v|^p)|v|^{p-2}v,\qquad v\in H^1(\mathbb S^n),
		\label{eq:intro-sphere}
	\end{equation}
	where
	\begin{equation}
		A=-\Delta_{\mathbb S^n}+\frac{n(n-2)}4,\qquad
		J_\mu h(\xi)=\int_{\mathbb S^n}\frac{h(\eta)}{|\xi-\eta|^\mu}\,dV_\eta,
		\label{eq:conformal-laplacian}
	\end{equation}
	and
	\begin{equation}
		n\geq3,\qquad 0<\mu<n,\qquad p=\frac{2n-\mu}{n-2}.
		\label{eq:intro-parameters}
	\end{equation}
	Here \(|\xi-\eta|\) is the chordal distance. We fix the standard product
	action
	\begin{equation}
		G=O(k)\times O(\ell),\qquad k,\ell\geq2,\qquad k+\ell=n+1
		\label{eq:intro-group}
	\end{equation}
	on \(\mathbb S^n\subset\mathbb R^k\times\mathbb R^\ell\). The space of
	\(G\)-invariant functions in \(H^1(\mathbb S^n)\) is denoted by \(H\).
	
	Under stereographic projection, \eqref{eq:intro-sphere} becomes
	\begin{equation}
		-\Delta u(x)
		=\left(\int_{\mathbb R^n}\frac{|u(y)|^p}{|x-y|^\mu}\,dy\right)
		|u(x)|^{p-2}u(x),\qquad u\in D^{1,2}(\mathbb R^n).
		\label{eq:intro-euclidean}
	\end{equation}
	Here \(D^{1,2}(\mathbb R^n)\) is the completion of
	\(C_c^\infty(\mathbb R^n)\) in the Dirichlet norm. The critical
	scaling is
	\[
	u(x)\longmapsto\lambda^{(n-2)/2}u(\lambda x),\qquad\lambda>0.
	\]
	Both terms in the energy are invariant under this scaling, so
	\eqref{eq:intro-euclidean} has the structure of a critical variational
	problem. For the variational background of the Choquard equation, see
	\cite{Lieb1977,MorozVanSchaftingen2017}. Positive solutions of
	\eqref{eq:intro-euclidean} were classified in part of the parameter
	range in \cite{DuYang2019,GuoHuPengShuai2019}, while the full range
	in \eqref{eq:intro-parameters} follows from
	\cite[Corollary~1.5]{DaiLiuQin2021}.
	
	A natural local counterpart is the critical equation
	\[
	-\Delta u=|u|^{2^*-2}u
	\quad\text{in }\mathbb R^n,
	\]
	where $2^*={2n}/{(n-2)}$. Its conformal interpretation is described in \cite{LeeParker1987}.
The positive solutions were classified in
\cite{CaffarelliGidasSpruck1989}. The linearized analysis of the
standard bubble \cite[Appendix, Lemma~A.1]{BianchiEgnell1991}
shows that it has one negative direction and an
\((n+1)\)-dimensional kernel generated by conformal transformations.
Ding \cite{Ding1986} used product symmetries to construct an
unbounded sequence of sign-changing solutions of the local critical
equation. 
	Almutairi and Gluck \cite[Theorem~1.1]{AlmutairiGluck2026} established
	the corresponding existence result for \eqref{eq:intro-euclidean} for
	every \(0<\mu<n\). They used both the standard action and an action
	with a sign character. Under the standard action, their argument for
	sign change used the positive classification in dimensions \(3,4\)
	or in dimensions \(n\geq5\) with \(\mu\leq4\).
Related variational results concern least-action nodal solutions of
Choquard equations with a mass term
\cite{GhimentiVanSchaftingen2016}, including the quadratic case
\cite{GhimentiMorozVanSchaftingen2017}. Morse-index estimates for
minimax solutions based on critical groups were developed in
\cite{BahriLions1988,Fang1995}.

Building on these developments, the present paper studies the cohomological minimax levels in the fixed \(G\)-invariant space \(H\). We prove that these levels are attained, establish their sharp growth behaviour, and obtain infinitely many pairs of sign-changing solutions. We further relate nodal domains to the invariant Morse index and minimax levels, and analyze the kernel and transverse negative directions of the full Hessian.
	
Equip \(H^1(\mathbb S^n)\) with the energy norm
\begin{equation}
	\|v\|^2=\int_{\mathbb S^n}
	\left(|\nabla_{\mathbb S^n}v|^2+\frac{n(n-2)}{4}v^2\right)dV.
	\label{eq:intro-operators}
\end{equation}
Its scalar product is denoted by \(\langle\cdot,\cdot\rangle_A\), and by
\(\langle\cdot,\cdot\rangle_H\) when restricted to \(H\). Let \(\|\cdot\|_s\) denote
the \(L^s(\mathbb S^n)\) norm. {The exponent \(2^*\) is as defined above.} Define
\begin{equation}
	\mathcal D(v)=\int_{\mathbb S^n}J_\mu(|v|^p)|v|^p\,dV,\qquad
	E(v)=\frac12\|v\|^2-\frac1{2p}\mathcal D(v).
	\label{eq:intro-energy}
\end{equation}
For \(w\in\mathcal S=\{w\in H:\|w\|=1\}\), let
\begin{equation}
	\Psi(w)=\max_{t>0}E(tw)=\frac{p-1}{2p}\mathcal D(w)^{-1/(p-1)}.
	\label{eq:intro-reduced}
\end{equation}
	The $\mathbb{Z}_2$-cohomological index $i$ due to Fadell and Rabinowitz \cite{FadellRabinowitz1978} is employed. For a finite-dimensional subspace $Y\subset H$, let $\mathbb S(Y)=\{y\in Y:\|y\|=1\}$. The following properties of this index will be used.
	\begin{enumerate}
		\item[(i)] If there is a continuous odd map between closed symmetric
		sets \(B\to C\), then \(i(B)\leq i(C)\).
		\item[(ii)] If \(B,C\) are closed symmetric sets, then
		\(i(B\cup C)\leq i(B)+i(C)\).
		\item[(iii)] One has \(i(\mathbb S(Y))=\dim Y\).
		\item[(iv)] Every compact symmetric set has a closed symmetric
		neighborhood with the same index.
	\end{enumerate}
These properties are given in \cite[Theorem~5.1]{FadellRabinowitz1978}.
For \(j\geq 1\), define
\begin{equation}
	\Gamma_j=\left\{{B\subset\mathcal S:B\text{ is compact},\ B=-B,\ i(B)\geq j}\right\}.
	\label{eq:admissible-family}
\end{equation}
The associated minimax levels are then
\begin{equation}
	c_j=\inf_{B\in\Gamma_j}\sup_{w\in B}\Psi(w).
	\label{eq:intro-levels}
\end{equation}
	Throughout the paper, \(C>0\) stands for a constant whose value may vary from one estimate to another. Subscripts are used to specify its dependencies. Let
	\begin{align}
	\kappa_0=\int_{\mathbb S^n}|\xi-\eta|^{-\mu}\,dV_\eta,
	\qquad
	a_0=\left(\frac{n(n-2)}{4\kappa_0}\right)^{1/(2p-2)}.
	\label{eq:constant-state}
	\end{align}
	By rotational invariance, this integral is independent of \(\xi\),
	and \(0<\mu<n\) implies \(0<\kappa_0<\infty\). The constants \(a_0\) and \(-a_0\) furnish constant-valued solutions to \eqref{eq:intro-sphere}.

	\begin{theorem}
		\label{thm:levels-intro}
		Assume \eqref{eq:intro-parameters} and \eqref{eq:intro-group}.
		There are constants \(C_1,C_2>0\), independent of \(j\), such that
		\begin{equation}
			C_1j^{2p/(p-1)}\leq c_j\leq C_2j^{2p/(p-1)},\qquad j\geq2.
			\label{eq:intro-j-bounds}
		\end{equation}
		Each \(c_j\) is attained by a solution in \(H\). The first level is
		attained only by the two constant solutions \(\pm a_0\) of
		\eqref{eq:constant-state}, and \(c_2>c_1\). For every \(j\geq2\), each
		invariant solution of energy \(c_j\) changes sign.
	\end{theorem}
	
	For \(p\geq2\), let \(m_G(v)\) and \(\nu_G(v)\) be the index and nullity of
	\(E''(v)\) on \(H\), and let \(m(v)\) and \(\nu(v)\) be the corresponding
	full-space quantities. 
	\begin{remark}
		The growth and attainment assertions hold for the whole range
		\(0<\mu<n\). The second-variation arguments below require
		\(p\geq2\). Indeed, for \(1<p<2\), the second variation at a nodal
		solution involves the singular coefficient \(|v|^{p-2}\), so the
		Hessian is not covered by the present framework. At \(p=2\), this
		coefficient is identically one, including on the nodal set.
	\end{remark}
	
	\begin{theorem}
		\label{thm:nodal-intro}
		Assume \(p\geq2\), and let \(v\in H\setminus\{0\}\) solve
		\eqref{eq:intro-sphere}. Then \(v\in C^{2,\gamma}(\mathbb S^n)\) for some
		\(\gamma\in(0,1)\). Its profile \(f(t)=v(\cos t\,\omega,\sin t\,\theta)\)
		has finitely many zeros, all simple and in \((0,\pi/2)\).
		If these zeros are \(t_1<\cdots<t_q\), then
		\begin{equation}
			\{v=0\}=\bigcup_{i=1}^{q}
			\mathbb S^{k-1}(\cos t_i)\times\mathbb S^{\ell-1}(\sin t_i).
			\label{eq:intro-nodal-set}
		\end{equation}
		Let \(N_G(v)\) be the number of connected components of \(\{v\ne0\}\).
		Then
		\begin{equation}
			N_G(v)=q+1\leq m_G(v),\qquad c_{N_G(v)}\leq E(v).
			\label{eq:intro-nodal-morse}
		\end{equation}
		For \(1\leq j<N_G(v)\) one has \(c_j<E(v)\). In particular,
		\(E(v)=c_j\) implies \(N_G(v)\leq j\).
		The least energy of sign-changing solutions in \(H\) is \(c_2\).
		Every solution in \(H\) at \(c_2\) has exactly two nodal domains.
	\end{theorem}
	
	The next result compares the minimax label with the index of a selected
	solution. Here isolation is taken in the entire critical set of \(\Psi\).
	
	\begin{theorem}
		\label{thm:morse-intro}
		Assume \(p\geq2\). {Suppose every critical point of
			\(\Psi\) at level \(c_j\) is isolated in the entire critical set of \(\Psi\).} Then there is a solution \(v_j\in H\) at 
		\(c_j\) such that
		\begin{equation}
			m_G(v_j)\leq j\leq m_G(v_j)+\nu_G(v_j).
			\label{eq:intro-minimax-morse}
		\end{equation}
		In particular, \(\nu_G(v_j)=0\) implies \(m_G(v_j)=j\).
	\end{theorem}
	
	\begin{remark}
		The isolation hypothesis in Theorem~\ref{thm:morse-intro} is imposed
		at the specified level \(c_j\). It has not been proved along the
		whole minimax sequence. Under this hypothesis, Palais--Smale
		compactness implies that the critical set at \(c_j\) is finite
		and that \(c_j\) is isolated among critical values. These are
		consequences, rather than additional assumptions.
	\end{remark}
	
For the full-space Hessian, we use the angular degrees \((a,b)\)
of the two factors, and denote the corresponding radial index and
nullity by \(m_{a,b}(v)\) and \(\nu_{a,b}(v)\). These quantities are
defined in Section~9. Under stereographic projection, the constant solutions \(a_0\) and
\(-a_0\) correspond respectively to the standard positive bubble
and its negative. Their linearized operators coincide. The
nondegeneracy of the positive bubble was studied in
\cite{DuYang2019,GiacomoniWeiYang2020} and established for the full
parameter range in \cite{LiLiuTangXu2026}.
	
	\begin{theorem}
		\label{thm:conformal-kernel-intro}
		Assume \(p\geq2\), and let \(v\in H^1(\mathbb S^n)\setminus\{0\}\) solve
		\eqref{eq:intro-sphere}. Set \(u=Pv\), where \(P\) is the stereographic
		isometry in \eqref{eq:conformal-map}. The Euclidean linearized kernel
		contains \(n+1\) linearly independent elements
		\[
		Z_i=\partial_{x_i}u\quad (1\leq i\leq n),\qquad
		Z_0=\frac{n-2}{2}u+x\cdot\nabla u .
		\]
		Thus \(\nu(v)\geq n+1\). If \(v\in H\) is nonconstant, the tangent space
		to its conformal orbit is a subspace of the kernel of dimension
		\(n+1+k\ell\). For such a nonconstant \(G\)-invariant solution,
		\begin{equation}
			m_{1,0}(v)+m_{0,1}(v)\geq1,\qquad
			m(v)\geq m_G(v)+\min\{k,\ell\}.
			\label{eq:intro-transverse-index}
		\end{equation}
		For the constant solutions, the conformal tangent space exhausts the kernel.
	\end{theorem}
	
When $k=\ell$, let $\mathcal T(\xi',\xi'')=(\xi'',\xi')$ and
$H^-=\{v\in H:v(\mathcal T\xi)=-v(\xi)\}$. As shown in the next theorem, a minimizer of $\Psi$ on
$\mathcal S\cap H^-$ gives rise to a solution via the rescaling
\eqref{eq:rescaling-critical}, with precise spectral and nodal properties described below.

	\begin{theorem}
		\label{thm:swap-blocks}
		Suppose \(k=\ell\). There exists \(0<\mu_*=\mu_*(k)<\min\{4,n\}\) such that for
		\(0<\mu<\mu_*\), a minimizer on \(\mathcal S\cap H^-\) yields a
		sign-changing solution \(v_\mu\) with two nodal domains. For every
		solution obtained from such a minimizer,
		\begin{equation}
			\begin{aligned}
				m_{1,0}(v_\mu)=m_{0,1}(v_\mu)&=1, & m_{1,1}(v_\mu)&=0,\\
				\nu_{1,0}(v_\mu)=\nu_{0,1}(v_\mu)&=1, & \nu_{1,1}(v_\mu)&=1.
			\end{aligned}
			\label{eq:intro-swap-inertia}
		\end{equation}
		In each of these three blocks the kernel is generated by the
		corresponding symmetry profile.
	\end{theorem}
	
\begin{remark}
	The minimization in Theorem~\ref{thm:swap-blocks} is performed
	in the proper subspace \(H^-\subset H\), where \(n=2k-1\).
	Therefore, the energy of the resulting solutions is not identified
	with the minimax level \(c_2\) in \(H\). The nullities in
	\eqref{eq:intro-swap-inertia} refer to the radial operators for fixed
	angular modes; after angular multiplicities are included, the
	corresponding full-space kernel dimensions are \(k\), \(\ell\), and
	\(k\ell\). The theorem determines the index and nullity only in the
	three sectors \((1,0)\), \((0,1)\), and \((1,1)\), rather than the
	total full-space index and nullity.
\end{remark}

Sections~\ref{sec:conformal}--\ref{sec:attainment} develop the
conformal and invariant variational framework and prove the growth
and attainment of the minimax levels.
Sections~\ref{sec:regularity}--\ref{sec:morse-label} establish the
nodal geometry, the comparison with the minimax levels, and the
critical-group relation between the minimax label, index, and
nullity. Section~\ref{sec:full-hessian} develops the angular
decomposition of the full Hessian and proves the corresponding
kernel and index results.
	\section{Conformal equivalence and invariant compactness}
	\label{sec:conformal}
	\noindent
	
	Let \(S=(0,\ldots,0,-1)\) and let
	\(\pi:\mathbb S^n\setminus\{S\}\to\mathbb R^n\) be stereographic
	projection,
	\[
	\pi(\xi)=\frac{(\xi_1,\ldots,\xi_n)}{1+\xi_{n+1}},\qquad
	\pi^{-1}(x)=\left(\frac{2x}{1+|x|^2},
	\frac{1-|x|^2}{1+|x|^2}\right).
	\]
	Set \(\Omega(x)=2/(1+|x|^2)\). Define
	\begin{equation}
		Pv(x)=\Omega(x)^{(n-2)/2}v(\pi^{-1}(x)).
		\label{eq:conformal-map}
	\end{equation}
	
	\begin{proposition}
		\label{prop:conformal-identities}
		The map \(P\) extends uniquely to an isometry from
		\(H^1(\mathbb S^n)\), equipped with the norm
		\eqref{eq:intro-operators}, onto \(D^{1,2}(\mathbb R^n)\).
		For every \(v\in H^1(\mathbb S^n)\),
		\begin{equation}
			\int_{\mathbb R^n}|\nabla Pv|^2\,dx=\|v\|^2,
			\qquad
			\iint_{\mathbb R^n\times\mathbb R^n}
			\frac{|Pv(x)|^p|Pv(y)|^p}{|x-y|^\mu}\,dx\,dy
			=\mathcal D(v).
			\label{eq:conformal-identities}
		\end{equation}
		Consequently, \(u=Pv\) solves \eqref{eq:intro-euclidean} if and only
		if \(v\) solves \eqref{eq:intro-sphere}.
	\end{proposition}
	
\begin{proof}
	All assertions are precisely the conformal equivalence established
	in \cite[Section~3]{AlmutairiGluck2026}.
\end{proof}
	
The two singular orbits are
\(\mathbb S^{k-1}\times\{0\}\) and
\(\{0\}\times\mathbb S^{\ell-1}\), of dimensions \(k-1\) and
\(\ell-1\), respectively. Every principal orbit is diffeomorphic to
\(\mathbb S^{k-1}\times\mathbb S^{\ell-1}\) and has dimension
\(k+\ell-2=n-1\). Put
\[
d=\min\{k-1,\ell-1\},\qquad
N=n-d=\max\{k,\ell\}.
\]
Then \(d\geq1\) and \(2\leq N<n\).
	
	\begin{lemma}
		\label{lem:compact-embedding}
		If \(N>2\), then
		\[
		H\hookrightarrow L^s(\mathbb S^n)
		\quad\text{continuously for }2\leq s\leq\frac{2N}{N-2},
		\]
		and the embedding is compact for \(2\leq s<2N/(N-2)\).
		If \(N=2\), it is continuous and compact for every finite
		\(s\geq2\). In particular, \(H\hookrightarrow L^{2^*}(\mathbb S^n)\)
		is compact, and the embedding remains compact for some
		\(s>2^*\).
	\end{lemma}
	
	\begin{proof}
		Apply the invariant Sobolev embedding in
		\cite[Corollary~1]{HebeyVaugon1997} with minimum orbit dimension
		\(d\). For \(N>2\), its critical exponent is \(2N/(N-2)\).
		Since \(N<n\),
		\[
		2^*=\frac{2n}{n-2}<\frac{2N}{N-2}.
		\]
		We may choose \(s\) strictly between these two exponents.
		
		For \(N=2\), fix any finite \(s\geq2\). Choose \(q\in(1,2)\)
		so close to \(2\) that \(s<2q/(2-q)\). The sphere has finite
		volume, and hence
		\[
		\|v\|_{W^{1,q}(\mathbb S^n)}\leq C_q\|v\|_{H^1(\mathbb S^n)}.
		\]
		Since \(s<2q/(2-q)\), the invariant \(W^{1,q}\) embedding into
		\(L^s\) is compact. Composing it with the preceding inclusion gives
		the compact embedding of \(H\) into \(L^s\), for every finite
		\(s\), including \(s>2^*\).
	\end{proof}
	
	Write
	\[
	q_\mu=\frac{2n}{2n-\mu},\qquad pq_\mu=2^*.
	\]
	The HLS inequality gives a bounded map
	\(J_\mu:L^{q_\mu}(\mathbb S^n)\to L^{2n/\mu}(\mathbb S^n)\),
	and
	\begin{equation}
		\mathcal D(v)\leq C\|v\|_{2^*}^{2p}.
		\label{eq:interaction-hls}
	\end{equation}
	The map \(v\mapsto |v|^p\) is of class
	\(C^1:L^{2^*}\to L^{q_\mu}\). Its derivative at \(v\) is
	\(\phi\mapsto p|v|^{p-2}v\phi\).
	Indeed, H\"older's inequality uses
	\((p-1)/2^*+1/2^*=1/q_\mu\), and the continuity of this derivative
	follows from the power estimates below. Thus
	\(\mathcal D,E\in C^1(H^1(\mathbb S^n),\mathbb R)\), with
	\begin{equation}
		\mathcal D'(v)[\phi]
		=2p\int_{\mathbb S^n}J_\mu(|v|^p)|v|^{p-2}v\phi\,dV.
		\label{eq:D-first}
	\end{equation}
	The kernel is positive, so \(\mathcal D(v)>0\) for every nonzero
	\(v\). We write \(H^*\) for the continuous dual of \(H\).
	
	\begin{lemma}
		\label{lem:derivative-compact}
		If \(v_m\rightharpoonup v\) in \(H\), then
		\begin{equation}
			\mathcal D'(v_m)\longrightarrow\mathcal D'(v)
			\quad\text{in }H^*.
			\label{eq:D-compact}
		\end{equation}
	\end{lemma}
	
\begin{proof}	
{
	Lemma \ref{lem:compact-embedding} yields
	\(
	v_m\to v~\text{in }L^{2^*}(\mathbb S^n).
	\)
	Since $pq_\mu=2^*$, the estimate
	\[
	\big\||a|^p-|b|^p\big\|_{q_\mu}
	\leq C\big(\|a\|_{2^*}^{p-1}
	+\|b\|_{2^*}^{p-1}\big)
	\|a-b\|_{2^*}
	\]
	implies
	\(
	|v_m|^p\to|v|^p
	~\text{in }L^{q_\mu}(\mathbb S^n).
	\)
	By the HLS inequality, we have
	\begin{equation}
		\|J_\mu(|v_m|^p-|v|^p)\|_{{2n}/{\mu}}
		\leq C\left\||v_m|^p-|v|^p\right\|_{q_\mu}=o(1).
		\label{eq:potential-convergence}
	\end{equation}
	Define $g(s)=|s|^{p-2}s$.
	For $p\geq2$, we have
	\[
	|g(a)-g(b)|\leq C(|a|+|b|)^{p-2}|a-b|,
	\]
	while for $1<p<2$, we have
	\[
	|g(a)-g(b)|\leq C|a-b|^{p-1}.
	\]
	Applying Hölder's inequality for $p\geq2$ and direct integration for $1<p<2$, we obtain
	\begin{equation}
		g(v_m)\longrightarrow g(v)
		\quad\text{in }L^{{2^*}/{(p-1)}}(\mathbb S^n).
		\label{eq:power-convergence}
	\end{equation}
	Moreover, the sequences
	\(
	\big\{J_\mu(|v_m|^p)\big\}_{m\geq1}\)
	and
	\(\big\{g(v_m)\big\}_{m\geq1}
	\)
	are bounded in
	\(
	L^{{2n}/{\mu}}(\mathbb S^n)\)
	and
	\(L^{{2^*}/{(p-1)}}(\mathbb S^n).
	\)
	Since
	\[
	\frac{\mu}{2n}
	+\frac{p-1}{2^*}
	+\frac1{2^*}
	=1,
	\]
	the product in \eqref{eq:D-first} can be estimated directly by Hölder's inequality.
	For \(\phi\in H\), Hölder's inequality gives
	\[
	\begin{aligned}
		\big|(\mathcal D'(v_m)-\mathcal D'(v))[\phi]\big|
		&\leq
		2p\,
		\big\|J_\mu(|v_m|^p-|v|^p)\big\|_{{2n}/{\mu}}
		\|g(v_m)\|_{{2^*}/{(p-1)}}
		\|\phi\|_{2^*} \\
		&\quad+
		2p\,
		\|J_\mu(|v|^p)\|_{{2n}/{\mu}}
		\|g(v_m)-g(v)\|_{{2^*}/{(p-1)}}
		\|\phi\|_{2^*}.
	\end{aligned}
	\]
	Hence, by \eqref{eq:potential-convergence},
	\eqref{eq:power-convergence}, and
	\(\|\phi\|_{2^*}\leq C\|\phi\|\), we obtain
	\[
	\left|
	(\mathcal D'(v_m)-\mathcal D'(v))[\phi]
	\right|
	\leq o(1)\|\phi\|.
	\]
	Taking the supremum over \(\|\phi\|\leq1\) proves
	\eqref{eq:D-compact}.
}
\end{proof}
	
	\begin{lemma}
		\label{lem:psi-ps}
		The functional \(\Psi:\mathcal S\to(0,\infty)\) is even, bounded
		away from zero, and of class \(C^1\). It satisfies the
		Palais--Smale condition.
	\end{lemma}
	
	\begin{proof}
		Since \(\mathcal D\) is positive on \(\mathcal S\), even, and
		of class \(C^1\), formula \eqref{eq:intro-reduced} gives the evenness
		and \(C^1\) regularity of \(\Psi\).
		For \(w\in\mathcal S\), \eqref{eq:interaction-hls} and Sobolev
		embedding give \(\mathcal D(w)\leq C\). Thus
		\[
		\inf_{\mathcal S}\Psi\geq\frac{p-1}{2p}C^{-1/(p-1)}>0.
		\]
		Let \(w_m\in\mathcal S\) be a Palais--Smale sequence.
		After taking a subsequence,
		\[
		\Psi(w_m)\to c>0,\qquad
		\|d(\Psi|_{\mathcal S})(w_m)\|\to0,\qquad
		w_m\rightharpoonup w\quad\text{in }H.
		\]
		By \eqref{eq:intro-reduced},
		\[
		\mathcal D(w_m)\longrightarrow
		D_c:=\left(\frac{p-1}{2pc}\right)^{p-1}>0.
		\]
		For
		\(T_w\mathcal S=\{\phi\in H:\langle w,\phi\rangle_H=0\}\),
		differentiation gives
		\[
		d(\Psi|_{\mathcal S})(w)[\phi]
		=-\frac1{2p}\mathcal D(w)^{-p/(p-1)}\mathcal D'(w)[\phi].
		\]
		It follows that
		\[
		\sup_{\substack{\phi\in T_{w_m}\mathcal S\\\|\phi\|\leq1}}
		|\mathcal D'(w_m)[\phi]|\longrightarrow0.
		\]
		Thus, the tangential part of \(\mathcal D'(w_m)\) tends to zero.
Euler's identity gives
	\[
	\mathcal D'(w_m)[w_m]=2p\mathcal D(w_m).
	\]
	The decomposition of \(H\) into the tangential subspace
	\(T_{w_m}\mathcal S\) and the radial subspace
	\(\operatorname{span}\{w_m\}\) therefore yields
	\begin{equation}
		\mathcal D'(w_m)
		-
		2p\mathcal D(w_m)\langle w_m,\cdot\rangle_H
		\longrightarrow0
		\quad\text{in }H^*.
		\label{eq:ps-lagrange}
	\end{equation}
		Let \(\mathcal R:H^*\to H\) be the Riesz isomorphism.
		Applying it to \eqref{eq:ps-lagrange} yields
		\[
		w_m=\frac{1}{2p\mathcal D(w_m)}\mathcal R\mathcal D'(w_m)+o(1)
		\quad\text{in }H.
		\]
		Lemma~\ref{lem:derivative-compact} now implies
		\[
		w_m\longrightarrow\frac{1}{2pD_c}\mathcal R\mathcal D'(w)
		\quad\text{strongly in }H.
		\]
		This limit equals the weak limit \(w\). Consequently
		\(\|w\|=\displaystyle \lim_{m\to\infty}\|w_m\|=1\), and \(w\in\mathcal S\).
		The Palais--Smale condition follows.
	\end{proof}
	
	\begin{lemma}
		\label{lem:radial-reduction}
		A point \(w\in\mathcal S\) is a critical point of \(\Psi\) if and only if
		\begin{equation}
			v=t(w)w,\qquad t(w)=\mathcal D(w)^{-1/(2p-2)},
			\label{eq:rescaling-critical}
		\end{equation}
		is a nonzero critical point of \(E|_H\). Every nonzero critical point
		of \(E|_H\) has this form, and \(E(v)=\Psi(w)\).
	\end{lemma}
	
	\begin{proof}
		For fixed \(w\in\mathcal S\), homogeneity gives
		\[
		E(tw)=\frac{t^2}{2}-\frac{t^{2p}}{2p}\mathcal D(w),\qquad
		\partial_tE(tw)=t-t^{2p-1}\mathcal D(w).
		\]
	The function \(t\mapsto E(tw)\) has a unique positive critical
	point at \(t=t(w)\), where it attains its global maximum.
		Suppose first that \(d\Psi(w)=0\).
		The derivative formula in the preceding proof shows that
		\(\mathcal D'(w)\) vanishes on \(T_w\mathcal S\).
		Its radial component is \(2p\mathcal D(w)\). Therefore
		\begin{equation}
			\mathcal D'(w)[\phi]
			=2p\mathcal D(w)\langle w,\phi\rangle_H,
			\qquad \phi\in H.
			\label{eq:D-lagrange-exact-rewritten}
		\end{equation}
		Put \(t=t(w)\). Then \(t^{2p-2}\mathcal D(w)=1\), and
		\[
		\begin{aligned}
			E'(tw)[\phi]
			&=t\langle w,\phi\rangle_H-\frac{t^{2p-1}}{2p}\mathcal D'(w)[\phi]\\
			&=t\big(1-t^{2p-2}\mathcal D(w)\big)\langle w,\phi\rangle_H=0.
		\end{aligned}
		\]
		Thus \(tw\) is a nonzero critical point of \(E|_H\).
		Testing the equation with \(tw\) gives
		\begin{equation}
			\|tw\|^2=\mathcal D(tw),\qquad
			E(tw)=\frac{p-1}{2p}\|tw\|^2.
			\label{eq:critical-testing}
		\end{equation}
		In particular,
		\[
		E(tw)=\frac{p-1}{2p}t^2
		=\frac{p-1}{2p}\mathcal D(w)^{-1/(p-1)}=\Psi(w).
		\]
		
		Conversely, let \(v\in H\setminus\{0\}\) satisfy \((E|_H)'(v)=0\).
		Set \(t=\|v\|\) and \(w=v/t\in\mathcal S\).
		Testing the equation with \(v\) gives \(\|v\|^2=\mathcal D(v)\). Hence
		\[
		t^2=t^{2p}\mathcal D(w),\qquad
		t^{2p-2}\mathcal D(w)=1,
		\]
		which proves \(t=t(w)\). For \(\phi\in T_w\mathcal S\),
		\[
		0=E'(tw)[\phi]
		=t\langle w,\phi\rangle_H-\frac{t^{2p-1}}{2p}\mathcal D'(w)[\phi]
		=-\frac{t^{2p-1}}{2p}\mathcal D'(w)[\phi].
		\]
		Thus \(\mathcal D'(w)\) vanishes on the tangent space, and the
		formula for \(d\Psi\) gives \(d\Psi(w)=0\).
		The energy equality follows from \eqref{eq:critical-testing}.
	\end{proof}
	
	\section{The orbit interval and the Jacobi spectrum}
	\label{sec:orbit}
	\noindent
	
	In this section, \(A\) denotes the self-adjoint realization of
	\eqref{eq:conformal-laplacian} in \(L^2(\mathbb S^n)\), with domain
	\(H^2(\mathbb S^n)\). Outside the two singular orbits, write
	\begin{equation}
		\xi=(\cos t\,\omega,\sin t\,\theta),\qquad
		0<t<\frac\pi2,\quad \omega\in\mathbb S^{k-1},\quad
		\theta\in\mathbb S^{\ell-1}.
		\label{eq:orbit-coordinates}
	\end{equation}
	The metric and volume element are
	\[
	dt^2+\cos^2t\,g_{\mathbb S^{k-1}}+\sin^2t\,g_{\mathbb S^{\ell-1}},
	\qquad
	dV=\rho(t)\,dt\,dV_\omega\,dV_\theta,
	\]
	where
	\begin{equation}
		\rho(t)=(\cos t)^{k-1}(\sin t)^{\ell-1}.
		\label{eq:orbit-weight}
	\end{equation}
	A \(G\)-invariant function has a profile \(v(\xi)=f(t)\).
	For smooth invariant functions,
	\begin{equation}
		\begin{aligned}
			\|v\|_2^2
			&=|\mathbb S^{k-1}||\mathbb S^{\ell-1}|
			\int_0^{\pi/2}|f(t)|^2\rho(t)\,dt,\\
			\|v\|^2
			&=|\mathbb S^{k-1}||\mathbb S^{\ell-1}|
			\int_0^{\pi/2}\left(|f'(t)|^2+\frac{n(n-2)}4|f(t)|^2\right)
			\rho(t)\,dt.
		\end{aligned}
		\label{eq:invariant-profile-norm}
	\end{equation}
	For \(v\in H\), take smooth approximants in \(H^1\) and average
	them over \(G\). Averaging is a contraction in the energy norm.
	The averaged sequence is smooth, invariant, and converges to \(v\),
	so the identities extend to \(H\). On every compact interval
	\(I\subset(0,\pi/2)\), the weight is bounded above and below by
	positive constants. The corresponding profiles converge in
	\(H^1(I)\). In particular, the profile of \(v\) has a continuous
	representative in the open orbit interval.
	
	\begin{lemma}
		\label{lem:reduced-laplacian}
		For a \(C^2\), \(G\)-invariant function \(v(\xi)=f(t)\),
		\begin{equation}
			\Delta_{\mathbb S^n}v
			=\rho(t)^{-1}(\rho(t)f'(t))'
			=f''(t)+\big((\ell-1)\cot t-(k-1)\tan t\big)f'(t).
			\label{eq:reduced-laplacian}
		\end{equation}
	\end{lemma}
	
	\begin{proof}
		All angular derivatives vanish. The divergence formula and the
		volume element above give
		\[
		\Delta_{\mathbb S^n}v
		=\frac1{\rho(t)}\frac d{dt}(\rho(t)f'(t))
		=f''(t)+\frac{\rho'(t)}{\rho(t)}f'(t).
		\]
		Differentiating \eqref{eq:orbit-weight}, we obtain
		\[
		\rho'(t)=-(k-1)(\cos t)^{k-2}(\sin t)^\ell
		+(\ell-1)(\cos t)^k(\sin t)^{\ell-2}.
		\]
		Division by \(\rho(t)>0\) gives
		\[
		\frac{\rho'(t)}{\rho(t)}
		=-(k-1)\frac{\sin t}{\cos t}
		+(\ell-1)\frac{\cos t}{\sin t}
		=-(k-1)\tan t+(\ell-1)\cot t.
		\]
		Substitution proves \eqref{eq:reduced-laplacian}.
	\end{proof}
	
	Set \(\alpha=(\ell-2)/2\) and \(\beta=(k-2)/2\).
	Let \(P_m^{(\alpha,\beta)}\) be the Jacobi polynomial of degree
	\(m\), with weight \((1-x)^\alpha(1+x)^\beta\)
	\cite{Szego1975}.
	
\begin{proposition}
	\label{prop:jacobi-spectrum}
	For \(m\geq0\), define
	\begin{equation}
		e_m(t)=\gamma_mP_m^{(\alpha,\beta)}(\cos 2t),
		\label{eq:jacobi-eigenfunction}
	\end{equation}
	where \(\gamma_m>0\) is chosen so that \(\|e_m\|_2=1\).
	Write \(L_G^2(\mathbb S^n)\) for the closed subspace of
	\(G\)-invariant functions in \(L^2(\mathbb S^n)\).
	The sequence \((e_m)_{m\geq0}\) is an orthonormal basis of
	\(L_G^2(\mathbb S^n)\) and an orthogonal basis of \(H\).
	The functions \(e_m\) satisfy the eigenvalue relations
	\begin{equation}
		-\Delta_{\mathbb S^n}e_m=2m(2m+n-1)e_m,
		\qquad
		Ae_m=\Lambda_me_m,
		\label{eq:invariant-eigenvalues}
	\end{equation}
	where
	\begin{equation}
		\Lambda_m
		=2m(2m+n-1)+\frac{n(n-2)}4.
		\label{eq:invariant-eigenvalue}
	\end{equation}
	In particular,
	\begin{equation}
		4m^2\leq\Lambda_m\leq C_n(m+1)^2.
		\label{eq:lambda-comparison}
	\end{equation}
\end{proposition}
	
	\begin{proof}
		Set \(x=\cos 2t\). Then
		\[
		\sin^2t=\frac{1-x}{2},
		\qquad
		\cos^2t=\frac{1+x}{2},
		\qquad
		|dx|=4\sin t\cos t\,dt.
		\]
		It follows that
		\[
		\begin{aligned}
			\rho(t)\,dt
			&=\frac14(\cos t)^{k-2}(\sin t)^{\ell-2}|dx|\\
			&=2^{-(k+\ell)/2}
			(1-x)^{(\ell-2)/2}(1+x)^{(k-2)/2}|dx|\\
			&=2^{-(k+\ell)/2}
			(1-x)^\alpha(1+x)^\beta|dx|.
		\end{aligned}
		\]
		Thus, for every integrable \(F\),
		\[
		\int_0^{\pi/2}F(\cos 2t)\rho(t)\,dt
		=
		2^{-(k+\ell)/2}
		\int_{-1}^1F(x)(1-x)^\alpha(1+x)^\beta\,dx.
		\]
		The correspondence
		\[
		v(\xi)=f(t)
		\quad\longleftrightarrow\quad
		F(x)=f\left(\frac12\arccos x\right)
		\]
		defines a linear bijection between \(L_G^2(\mathbb S^n)\) and
		\(
		L^2\bigl((-1,1),(1-x)^\alpha(1+x)^\beta\,dx\bigr),
		\)
		and the corresponding norms satisfy
		\[
		\|v\|_2^2
		=
		C_{k,\ell}
		\int_{-1}^1
		|F(x)|^2(1-x)^\alpha(1+x)^\beta\,dx,
		\]
		where
		\[
		C_{k,\ell}
		=
		2^{-(k+\ell)/2}
		|\mathbb S^{k-1}||\mathbb S^{\ell-1}|.
		\]
		Since \(\alpha,\beta>-1\), Jacobi orthogonality gives
		\[
		\begin{aligned}
			\langle e_m,e_j\rangle_{L^2}
			&=
			C_{k,\ell}\gamma_m\gamma_j
			\int_{-1}^1
			P_m^{(\alpha,\beta)}(x)
			P_j^{(\alpha,\beta)}(x)
			(1-x)^\alpha(1+x)^\beta\,dx=\delta_{mj}.
		\end{aligned}
		\]
		If \(v\in L_G^2(\mathbb S^n)\) is orthogonal to every \(e_m\),
		its transformed profile \(F\) is orthogonal to every Jacobi
		polynomial. Their completeness implies \(F=0\), and hence
		\(v=0\). Thus \((e_m)_{m\geq0}\) is an orthonormal basis of
		\(L_G^2(\mathbb S^n)\).
		
	Set \(y(x)=P_m^{(\alpha,\beta)}(x)\). The Jacobi differential
	equation \cite[Eq.~(4.2.4)]{Szego1975} gives
	\begin{equation}
		(1-x^2)y''
		+\bigl(\beta-\alpha-(\alpha+\beta+2)x\bigr)y'
		+m(m+\alpha+\beta+1)y=0.
		\label{eq:jacobi-equation}
	\end{equation}
	 The chain rule gives
		\[
		\frac d{dt}y(\cos 2t)
		=-2\sin 2t\,y'(x),
		\qquad
		\frac{d^2}{dt^2}y(\cos 2t)
		=4(1-x^2)y''(x)-4xy'(x).
		\]
		Substitution into \eqref{eq:reduced-laplacian} gives
		\[
		\Delta_{\mathbb S^n}[y(\cos 2t)]
		=
		4\left[
		(1-x^2)y''
		+\bigl(\beta-\alpha-(\alpha+\beta+2)x\bigr)y'
		\right].
		\]
    	Using \eqref{eq:jacobi-equation}, the preceding identity becomes
		\[
		-\Delta_{\mathbb S^n}[y(\cos 2t)]
		=
		4m(m+\alpha+\beta+1)y(\cos 2t).
		\]
		Since
		\(
		\alpha+\beta+1={(n-1)}/{2}
		\)
		and \(e_m(t)=\gamma_m y(\cos 2t)\), it follows that
		\[
		-\Delta_{\mathbb S^n}e_m
		=
		4m\left(m+\frac{n-1}{2}\right)e_m
		=
		2m(2m+n-1)e_m.
		\]
	Moreover, writing
	\[
	\xi=(\xi',\xi'')
	=(\cos t\,\omega,\sin t\,\theta)
	\in\mathbb R^k\times\mathbb R^\ell,
	\]
	we have
	\[
	\cos 2t
	=\cos^2t-\sin^2t
	=|\xi'|^2-|\xi''|^2.
	\]
	Thus \(e_m\) is the restriction to \(\mathbb S^n\) of a polynomial
	in the ambient coordinates. It is therefore smooth across both
	singular orbits, and the eigenvalue identity holds on the whole
	sphere. Adding the constant term in \(A\) gives
	\(
	Ae_m=\Lambda_me_m,
	\)
	with \(\Lambda_m\) as in
	\eqref{eq:invariant-eigenvalue}. This proves
	\eqref{eq:invariant-eigenvalues}.
		
		Rotations commute with \(A\), and averaging over \(G\) is the
		orthogonal projection onto \(L_G^2(\mathbb S^n)\). Hence
		\(L_G^2(\mathbb S^n)\) reduces \(A\).
		Let \(A_G\) be the corresponding restriction. Its form domain is
		\(H\), and the spectral theorem gives
		\[
		H=
		\left\{
		v=\sum_{m\geq0}\widehat v_me_m:
		\sum_{m\geq0}\Lambda_m|\widehat v_m|^2<\infty
		\right\},
		\qquad
		\widehat v_m=\langle v,e_m\rangle_{L^2}.
		\]
		In particular,
		\begin{equation}
			\|v\|^2
			=
			\sum_{m\geq0}\Lambda_m|\widehat v_m|^2,
			\qquad
			\|v\|_2^2
			=
			\sum_{m\geq0}|\widehat v_m|^2.
			\label{eq:invariant-parseval}
		\end{equation}
		The identity
		\[
		\langle e_m,e_j\rangle_H
		=
		\langle A^{1/2}e_m,A^{1/2}e_j\rangle_{L^2}
		=
		\Lambda_m\delta_{mj}
		\]
		proves energy orthogonality. For \(v\in H\), its partial sums
		converge in \(H\), since
		\[
		\left\|
		v-\sum_{m=0}^{M}\widehat v_me_m
		\right\|^2
		=
		\sum_{m>M}\Lambda_m|\widehat v_m|^2
		\longrightarrow0.
		\]
		Finally,
		\[
		\Lambda_m
		=
		4m^2+2(n-1)m+\frac{n(n-2)}4,
		\]
		which gives \eqref{eq:lambda-comparison}.
	\end{proof}

	\section{Sharp growth of the minimax levels}
	\label{sec:level-estimates}
	\noindent
	
Throughout this section, we use the minimax family \(\Gamma_j\)
defined in \eqref{eq:admissible-family}. In addition to the HLS estimate
\eqref{eq:interaction-hls}, the bound \(|\xi-\eta|\leq2\) gives
\begin{align}
	\mathcal D(w)
	&=\iint_{\mathbb S^n\times\mathbb S^n}
	\frac{|w(\xi)|^p|w(\eta)|^p}{|\xi-\eta|^\mu}
	\,dV_\xi\,dV_\eta\nonumber\\
	&\geq
	2^{-\mu}
	\left(\int_{\mathbb S^n}|w|^p\,dV\right)^2
	=
	2^{-\mu}\|w\|_p^{2p}.
	\label{eq:kernel-lower}
\end{align}
{For every subinterval \(I\subset[0,\pi/2]\), define}
\[
L_\rho(I)
=
\int_I\rho(t)^{-1/(n-1)}\,dt,
\qquad
e_f(I)
=
\int_I|f'(t)|^2\rho(t)\,dt.
\]
The weight \(\rho\) in \eqref{eq:orbit-weight} satisfies
\[
{
\rho(t)\asymp t^{\ell-1}
\quad(0<t\leq \pi/4),
\qquad
\rho(t)\asymp(\pi/2-t)^{k-1}
\quad(\pi/4\leq t<\pi/2).}
\]
Since \(k,\ell<n\), we have
\[
{\int_0^{\pi/2}\rho(t)^{-1/(n-1)}\,dt<\infty.}
\]
Only traces at interior points will be prescribed below. No trace
condition is imposed at either singular orbit.

\begin{lemma}
	\label{lem:orbit-partition-bound}
	Let \(f\) be the profile of a function in \(H\). If
	{\(I=[a,b]\subset[0,\pi/4]\)}, \(a<b\), and \(f(b)=0\), then
	\begin{equation}
		\|f\|_{L^{2^*}(I,\rho\,dt)}
		\leq C L_\rho(I)^{(n-1)/n} e_f(I)^{1/2}.
		\label{eq:orbit-cell-bound}
	\end{equation}
	Here \(C>0\) depends only on \(n\), \(k\), and \(\ell\). The same estimate holds on {\(I=[a,b]\subset[\pi/4,\pi/2]\)} with
	\(f(a)=0\).
\end{lemma}

\begin{proof}
	{On the left half, the profile belongs to
		\(H^1\) on each compact subinterval of {\((0,\pi/2)\)}.} For \(a<t<b\),
	\[
	f(t)=-\int_t^b f'(r)\,dr,
	\qquad
	|f(t)|\leq e_f(I)^{1/2}
	\left(\int_t^b\rho(r)^{-1}\,dr\right)^{1/2}.
	\]
	It remains to estimate
	\[
	\left(\int_a^b \rho(t)
	\left(\int_t^b\rho(r)^{-1}\,dr\right)^{2^*/2}dt\right)^{1/2^*}.
	\]
	Suppose first that \(a\geq b/2\). Then
	\(\rho(t)\asymp b^{{\ell}-1}\) on \([a,b]\), with constants independent
	of the interval. Consequently,
	\[
	\begin{aligned}
		\int_a^b\rho(t)
		\left(\int_t^b\rho(r)^{-1}\,dr\right)^{2^*/2}dt
		&\leq C b^{({\ell}-1)(1-2^*/2)}\int_a^b(b-t)^{2^*/2}\,dt\\
		&\leq C b^{({\ell}-1)(1-2^*/2)}(b-a)^{1+2^*/2}.
	\end{aligned}
	\]
	The identities
	\[
	\frac12+\frac1{2^*}=\frac{n-1}{n},
	\qquad \frac1{2^*}-\frac12=-\frac1n
	\]
	give
	\[
	\|f\|_{L^{2^*}(I,\rho\,dt)}
	\leq C e_f(I)^{1/2}
	(b-a)^{(n-1)/n}b^{-({\ell}-1)/n}.
	\]
	Also,
	\[
	L_\rho(I)\asymp (b-a)b^{-({\ell}-1)/(n-1)}.
	\]
	These two estimates imply \eqref{eq:orbit-cell-bound}.
	
	Suppose next that \(a<b/2\). Since \([b/2,b]\subset I\),
	\[
	L_\rho(I)\geq C^{-1}
	\int_{b/2}^b t^{-({\ell}-1)/(n-1)}\,dt
	\geq C^{-1}b^{(n-{\ell})/(n-1)}.
	\]
	If \({\ell}>2\), then
	\[
	\int_t^b\rho(r)^{-1}\,dr
	\leq C\int_t^b r^{1-{\ell}}\,dr
	\leq Ct^{2-{\ell}}.
	\]
	Moreover,
	\[
	{\ell}+\frac{(2-{\ell})2^*}{2}
	=\frac{2(n-{\ell})}{n-2}>0.
	\]
	We therefore obtain
	\[
	\begin{aligned}
		\int_a^b\rho(t)
		\left(\int_t^b\rho(r)^{-1}\,dr\right)^{2^*/2}dt
		&\leq C\int_0^b
		t^{{\ell}-1+(2-{\ell})2^*/2}\,dt\\
		&\leq Cb^{{\ell}+(2-{\ell})2^*/2}.
	\end{aligned}
	\]
	Taking the \(2^*\)-th root gives the factor
	\(b^{{\ell}/2^*+(2-{\ell})/2}
	=b^{1-{\ell}/n}\).
	If \({\ell}=2\), the power estimate is replaced by a logarithm.
	In this case,
	\[
	\int_t^b\rho(r)^{-1}\,dr\leq C\log(b/t)
	\]
	and the substitution \(t=bu\) gives
	\[
	\int_0^b t\big(\log(b/t)\big)^{2^*/2}\,dt
	=b^2\int_0^1u\big(\log(1/u)\big)^{2^*/2}\,du
	\leq Cb^2.
	\]
	Its \(2^*\)-th root is \(Cb^{2/2^*}=Cb^{1-2/n}\). Thus both cases
	yield
	\[
	\|f\|_{L^{2^*}(I,\rho\,dt)}
	\leq C e_f(I)^{1/2}b^{1-{\ell}/n}
	\leq C e_f(I)^{1/2}L_\rho(I)^{(n-1)/n}.
	\]
	If \(a=0\), apply the preceding estimate on
	\([\varepsilon,b]\) and let \(\varepsilon\downarrow0\).
	Monotone convergence gives the same estimate on \([0,b]\).
	{Reflection \(t\mapsto \pi/2-t\)} proves the right-half estimate, with
	\(k\) in place of \(\ell\).
\end{proof}

\begin{proposition}
	\label{prop:minimax-growth-bounds}
	There exist constants \(C_1,C_2>0\), independent of \(j\), such that
	\[
	C_1j^{2p/(p-1)}
	\leq c_j
	\leq C_2j^{2p/(p-1)},
	\qquad j\geq2.
	\]
\end{proposition}

\begin{proof}
	Let \(B\in\Gamma_j\), where \(j\geq2\). Choose positive integers \(j_-\) and \(j_+\) such that
	\[
	j_-+j_+=j,
	\qquad
	|j_--j_+|\leq1.
	\]
	For the profile \(f\)
	of \(w\in B\), put
	\[
	a_w(t)=\rho(t)^{-1/(n-1)}+\rho(t)|f'(t)|^2.
	\]
	By \eqref{eq:invariant-profile-norm}, the normalization
	\(\|w\|=1\), and the integrability of
	\(\rho^{-1/(n-1)}\), there exists \(C>0\), independent of \(j\)
	and \(w\), such that
	\[
	0<
		{\int_0^{\pi/2}\rho(t)^{-1/(n-1)}\,dt}
	\leq
		{\int_0^{\pi/2}a_w(t)\,dt}
	\leq C.
	\]
	{Partition \([0,\pi/4]\) into \(j_-\) intervals of equal
	\(a_w\)-mass and \([\pi/4,\pi/2]\) into \(j_+\) intervals in the same}
	way. Since
	\[
	a_w(t)\geq\rho(t)^{-1/(n-1)}>0
	\qquad{(0<t<\pi/2)},
	\]
	the corresponding cumulative mass functions are continuous and
	strictly increasing. Hence all division points are uniquely
	determined. On each cell \(I\),
	\begin{equation}
		L_\rho(I)+e_f(I)=\int_I a_w(t)\,dt\leq\frac{C}{j}.
		\label{eq:adaptive-cell-mass}
	\end{equation}
	Here we used \(\min\{j_-,j_+\}\geq j/3\), also valid when
	\(j=2\).
	
	To prove the continuity of the division points, let
	\(w_m,w\in\mathcal S\) satisfy \(w_m\to w\) strongly in \(H\),
	and denote their profiles by \(f_m\) and \(f\), respectively. Then 
	\[
	\begin{aligned}
		{\int_0^{\pi/2}}\rho\,\big||f_m'|^2-|f'|^2\big|\,dt
		&\leq
		\left({\int_0^{\pi/2}}\rho|f_m'-f'|^2\,dt\right)^{1/2}\\
		&\quad\times
		\left(\left({\int_0^{\pi/2}}\rho|f_m'|^2\,dt\right)^{1/2}
		+\left({\int_0^{\pi/2}}\rho|f'|^2\,dt\right)^{1/2}\right)
		\longrightarrow0.
	\end{aligned}
	\]
	Hence the cumulative mass functions converge uniformly. Let \(t\)
	be a division point for \(w\), and let \(\varepsilon>0\) be small.
	Strict monotonicity puts the cumulative mass at \(t-\varepsilon\)
	below the prescribed fraction and that at \(t+\varepsilon\) above
	it. Uniform convergence preserves these inequalities for \(w_m\).
	Its corresponding division point therefore lies between
	\(t-\varepsilon\) and \(t+\varepsilon\). Thus every division point
	converges.
	
	{List all interior nodes, including \(\pi/4\), as
	\(t_1(w)<\cdots<t_{j-1}(w)\). They lie in \((0,\pi/2)\).
	On each compact subinterval of \((0,\pi/2)\)}, convergence in \(H\)
	implies convergence in the ordinary \(H^1\) norm and hence uniform
	convergence. The map
	\[
	B\longrightarrow\mathbb R^{j-1},\qquad
	w\longmapsto\big(f(t_1(w)),\ldots,f(t_{j-1}(w))\big)
	\]
	is therefore continuous. It is odd, because
	\(a_{-w}=a_w\) and consequently \(t_i(-w)=t_i(w)\).
	If it never vanished, its normalization would be a continuous
	odd map from \(B\) to \(\mathbb S^{j-2}\). The index properties
	would then give \(i(B)\leq j-1\), contrary to \(B\in\Gamma_j\).
	We may thus choose \(w\in B\) whose profile vanishes at all nodes.
	
	For every left-half cell, the profile vanishes at its right
	endpoint. For every right-half cell, it vanishes at its left
	endpoint. Lemma~\ref{lem:orbit-partition-bound} and
	\eqref{eq:adaptive-cell-mass} imply
	\[
	\|f\|_{L^{2^*}(I,\rho\,dt)}
	\leq Cj^{-(n-1)/n-1/2}.
	\]
	The cells have disjoint interiors. Summing the \(2^*\)-th powers
	and including the fixed sphere-area factors gives
	\[
	\|w\|_{2^*}
	\leq Cj^{1/2^*-(n-1)/n-1/2}=Cj^{-1}.
	\]
	It follows from \eqref{eq:interaction-hls} that
	\[
	\mathcal D(w)\leq Cj^{-2p},
	\qquad
	\sup_{z\in B}\Psi(z)\geq\Psi(w)\geq C^{-1}j^{2p/(p-1)}.
	\]
	Taking the infimum over \(B\in\Gamma_j\) proves the lower bound.
	
	For the upper bound, fix {\([a,a+b]\subset(0,\pi/2)\)}, where \(b>0\),
	and choose a nonzero \(\chi\in C_c^\infty(0,1)\). For
	\(z=(z_1,\ldots,z_j)\in\mathbb S^{j-1}\), set
	\[
	h_i=bz_i^2,\qquad a_i=a+\sum_{r<i}h_r,\qquad
	\sum_{i=1}^jh_i=b.
	\]
	If \(h_i>0\), define on \([a_i,a_i+h_i]\)
	\[
	f_z(t)=\operatorname{sgn}(z_i)h_i
	\chi\left(\frac{t-a_i}{h_i}\right).
	\]
	Omit intervals of length zero and extend the profile by zero
	outside \([a,a+b]\). Since \(\chi\) vanishes near both ends,
	the pieces fit together in the weighted energy space.
	
	The weight is bounded above and below on the fixed interval.
	Changing variables \(t=a_i+h_i u\) on each nonempty piece gives
	\[
	\begin{aligned}
		\int_{a_i}^{a_i+h_i}|f_z'|^2\rho\,dt
		&=h_i\int_0^1|\chi'(u)|^2\rho(a_i+h_i u)\,du,\\
		\int_{a_i}^{a_i+h_i}|f_z|^2\rho\,dt
		&=h_i^3\int_0^1|\chi(u)|^2\rho(a_i+h_i u)\,du,\\
		\int_{a_i}^{a_i+h_i}|f_z|^p\rho\,dt
		&=h_i^{p+1}\int_0^1|\chi(u)|^p\rho(a_i+h_i u)\,du.
	\end{aligned}
	\]
	Let \(v_z\) be the corresponding invariant function.
	Since \(\int_0^1|\chi'|^2>0\), these identities give
	\begin{equation}
		C^{-1}\leq\|v_z\|\leq C,
		\qquad
		\|v_z\|_p^p\geq C^{-1}\sum_{i=1}^j h_i^{p+1}
		\geq C^{-1}j^{-p}.
		\label{eq:variable-bump-estimates}
	\end{equation}
	The last step uses convexity and \(\displaystyle \sum_i h_i=b\).
	The constants do not depend on \(j\) or \(z\).
	
	The map \(z\mapsto v_z\) is odd. It is continuous on each set
	where the nonzero coordinates have fixed signs, by continuity
	of translations and dilations of the smooth bump in \(H^1\).
	If a coordinate tends to zero, the squared energy norm of its
	piece is at most \(C(h_i+h_i^3)\), which tends to zero. The
	remaining endpoints vary continuously. This proves continuity
	on the whole sphere, including all its coordinate faces.
	
	By \eqref{eq:variable-bump-estimates}, the normalized image
	\[
	B_j=\{v_z/\|v_z\|:z\in\mathbb S^{j-1}\}
	\]
	is a compact symmetric subset of \(\mathcal S\). The continuous
	odd map \(\mathbb S^{j-1}\to B_j\) gives \(i(B_j)\geq j\).
	For every \(w\in B_j\), \eqref{eq:kernel-lower} and
	\eqref{eq:variable-bump-estimates} yield
	\[
	\mathcal D(w)\geq C^{-1}j^{-2p},
	\qquad \Psi(w)\leq Cj^{2p/(p-1)}.
	\]
	Since \(B_j\in\Gamma_j\), this proves the upper bound.
\end{proof}

\section{Attainment and sign change}
\label{sec:attainment}
\noindent

The following minimax lemma relies on the monotonicity, subadditivity,
and neighborhood properties of the cohomological index, together with
the \(\mathbb Z_2\)-equivariant deformation lemma.

\begin{lemma}
	\label{lem:symmetric-minimax}
	Let $F$ be an even $C^1$ functional on the unit sphere of an
	{infinite-dimensional} Hilbert space.  Assume that $\inf_{\mathcal S}F>0$ and that $F$ satisfies
	the {Palais--Smale} condition.  For
	\[
	d_j={\inf_{B\in\Gamma_j}}\sup_{w\in B}F(w),
	\]
	where {$\Gamma_j$ is defined by
		\eqref{eq:admissible-family}}, every finite $d_j$ is a
	critical value.  If $d_j=\cdots=d_{j+s}=d$, then the critical set $K_d$ has
	$i(K_d)\geq s+1$, where \(
	K_d
	=
	\{w\in\mathcal S:F(w)=d,\ dF(w)=0\}.
	\)
\end{lemma}
\begin{proof}
	Set \(F^a=\{w\in\mathcal S:F(w)\leq a\}\).
	The \(\mathbb Z_2\)-equivariant deformation lemma
	\cite[Lemma~3.1]{Willem1996}, applied on \(\mathcal S\) to a
	tangential pseudogradient, shows that if \(d_j\) were not a critical
	value, then, for some \(\varepsilon>0\), there would exist an odd
	homeomorphism \(\eta:\mathcal S\to\mathcal S\) such that
	\(
	\eta\bigl(F^{d_j+\varepsilon}\bigr)
	\subset F^{d_j-\varepsilon}.
	\)
	Choose \(B\in\Gamma_j\) such that
	\(\sup_BF<d_j+\varepsilon\).
	Since \(\eta\) is an odd homeomorphism, \(\eta(B)\) is compact and
	symmetric, and
	\(
	i(\eta(B))=i(B)\geq j.
	\)
	Thus \(\eta(B)\in\Gamma_j\). On the other hand,
	\(
	\sup_{\eta(B)}F\leq d_j-\varepsilon<d_j,
	\)
	which contradicts the definition of 
	\(
	d_j.
	\)
	Therefore, every finite \(d_j\) is a critical value.
	
	Suppose that
	\(
	d_j=\cdots=d_{j+s}=d
	\)
	and \(i(K_d)\leq s\). By the Palais--Smale condition, \(K_d\) is
	compact. The continuity property of the cohomological index
	\cite[Theorem~5.1(5)]{FadellRabinowitz1978} yields a symmetric open
	neighborhood \(U\) of \(K_d\) such that
	\(
	i(\overline U)\leq s.
	\)
	The equivariant deformation lemma then sends
	\(F^{d+\varepsilon}\setminus U\) into \(F^{d-\varepsilon}\) for
	sufficiently small \(\varepsilon>0\). Choose
	\(B\in\Gamma_{j+s}\) such that
	\(
	\sup_BF<d+\varepsilon.
	\)
	By the monotonicity and subadditivity properties of the index
	\cite[Theorem~5.1(4),(6)]{FadellRabinowitz1978}, we have
	\[
	j+s
	\leq i(B)
	\leq i(B\setminus U)+i(B\cap\overline U)
	\leq i(B\setminus U)+s.
	\]
	Hence \(i(B\setminus U)\geq j\). If \(\eta\) denotes the odd
	homeomorphism supplied by the deformation lemma, then
	\[
	i\bigl(\eta(B\setminus U)\bigr)
	=
	i(B\setminus U)
	\geq j.
	\]
	Thus \(\eta(B\setminus U)\in\Gamma_j\), whereas
	\[
	\sup_{\eta(B\setminus U)}F
	\leq d-\varepsilon
	<d_j,
	\]
	contradicting the definition of \(d_j\). Therefore,
	\(
	i(K_d)\geq s+1.
	\)
\end{proof}

We first classify the one-signed solutions in \(H\). 
\begin{lemma}
	\label{lem:positive-invariant}
	Let \(v\in H\setminus\{0\}\) be a weak solution of
	\eqref{eq:intro-sphere}. If \(v\) does not change sign, then
	\(
	v= a_0
	\) or \(
	v=-a_0,
	\)
where \(a_0\) is defined in \eqref{eq:constant-state}.
	The energy of these constant solutions is
	\begin{equation}
		E(a_0)=E(-a_0)
		=
		\frac{p-1}{2p}
		\frac{n(n-2)}{4}
		a_0^2|\mathbb S^n|
		=
		c_1.
		\label{eq:constant-energy}
	\end{equation}
	Moreover, the only \(G\)-invariant solutions at the first minimax level
	\(c_1\) are \(a_0\) and \(-a_0\).
\end{lemma}

\begin{proof}
	{
		Suppose first that \(v\geq0\), and set \(u=Pv\).
		Proposition~\ref{prop:conformal-identities} shows that \(u\) is a
		nonnegative, nonzero weak solution of \eqref{eq:intro-euclidean}.
		By \cite[Theorem~2.3]{AlmutairiGluck2026},
		\(
		u\in C_{\mathrm{loc}}^{2,\alpha}(\mathbb R^n)
		\)
		for some \(\alpha\in(0,1)\). Since
		\[
		-\Delta u
		=
		\left(
		\int_{\mathbb R^n}
		\frac{|u(y)|^p}{|x-y|^\mu}\,dy
		\right)u^{p-1}
		\geq0,
		\]
		the strong maximum principle {\cite[Theorem~3.5]{GilbargTrudinger1977}}
		and \(u\not\equiv0\) yield $u>0$ on $\mathbb R^n$.
		{Apply \cite[Corollary~1.5]{DaiLiuQin2021} with their
			$s=1$, $\sigma=\mu$, $\alpha=2$, $m=0$ and their two powers
			equal to $p$ and $p-1$. Indeed, the two critical powers in that
			classification become
			\[
			\frac{2n-\mu}{n-2}=p,\qquad
			\frac{n+2-\mu}{n-2}=p-1>0.
			\]
			The resulting Euclidean bubble, pulled back
			by $P^{-1}$, has the form}
		\begin{equation}
			v(\xi)
			=
			C(1-a\cdot\xi)^{-(n-2)/2},
			\qquad
			C>0,\quad
			a\in\mathbb R^{n+1},\quad
			|a|<1.
			\label{eq:classified-positive-state}
		\end{equation}
		The $G$-invariance of $v$ yields $v(g\xi)=v(\xi)$ for all $g\in G$.
		Substituting \eqref{eq:classified-positive-state} yields
		\[
		1-a\cdot g\xi
		=
		1-a\cdot\xi,
		\qquad
		g\in G,
		\quad
		\xi\in\mathbb S^n,
		\]
	Thus \(g^{-1}a=a\) for every \(g\in G\). The action of
	\(O(k)\times O(\ell)\) on \(\mathbb R^k\times\mathbb R^\ell\)
	has no nonzero fixed vector. Hence \(a=0\), and \eqref{eq:classified-positive-state}
	shows that \(v\) is constant. Write \(v\equiv b\) with \(b>0\).
	Then
	\[
	Av=\frac{n(n-2)}4b,
	\qquad
	J_\mu(|v|^p)|v|^{p-2}v
	=\kappa_0b^{2p-1}.
	\]
	Consequently, \eqref{eq:intro-sphere} reduces to
	\[
	\frac{n(n-2)}4b=\kappa_0b^{2p-1}.
	\]
	Dividing by \(b>0\) gives \(b=a_0\).
	If \(v\leq0\), then \(-v\) is a nonnegative nonzero solution. The
		preceding argument gives \(v=-a_0\).
	}
	
	{
	For \(v=a_0\), testing \eqref{eq:intro-sphere} against \(v\)
	gives \(\|v\|^2=\mathcal D(v)\). Consequently,
	\[
	E(a_0)
	=
	\frac12\|a_0\|^2-\frac1{2p}\mathcal D(a_0)
	=
	\frac{p-1}{2p}\|a_0\|^2
	=
	\frac{p-1}{2p}\frac{n(n-2)}4
	a_0^2|\mathbb S^n|.
	\]
	Since \(E\) is even, \(E(-a_0)=E(a_0)\).
	It remains to prove that this common value is \(c_1\) and to
	identify all invariant solutions at this level. Since every \(B\in\Gamma_1\) is contained in \(\mathcal S\), while
	\(\{w,-w\}\in\Gamma_1\) and \(\Psi(-w)=\Psi(w)\) for every
	\(w\in\mathcal S\), it follows from the definition of \(c_1\) that
	\begin{equation}
		c_1=\inf_{w\in\mathcal S}\Psi(w).
		\label{eq:first-level-infimum}
	\end{equation}
		Let $S_n>0$ be the best constant in
		\[
		S_n\|w\|_{L^{2^*}(\mathbb S^n)}^2\leq\|w\|^2,
		\qquad w\in H^1(\mathbb S^n),
		\]
		and let $C_{\mathrm{HLS}}(n,\mu)$ denote the best constant in the
		spherical HLS inequality. The sharp HLS inequality
		\cite{Lieb1983}, in the spherical formulation of
		\cite[Theorem~2.2]{FrankLieb2012}, and the sharp Sobolev inequality
		{\cite[Theorem~3.1]{FrankLieb2012} and \cite{Beckner1993}} give
		\begin{align}
			\mathcal D(w)
			&\leq
			C_{\mathrm{HLS}}(n,\mu)
			\|w\|_{L^{2^*}(\mathbb S^n)}^{2p}
			\leq
			C_{\mathrm{HLS}}(n,\mu)
			S_n^{-p}\|w\|^{2p}.
			\label{eq:sharp-combined-inequality}
		\end{align}
		Set
		\(
		C_*
		=
		C_{\mathrm{HLS}}(n,\mu)S_n^{-p}.
		\)
		For \(w\in\mathcal S\), \eqref{eq:sharp-combined-inequality} becomes
		\(
		\mathcal D(w)\leq C_*.
		\)
		Constant functions attain equality in both sharp inequalities.
		Therefore,
		\[
		\sup_{w\in\mathcal S}\mathcal D(w)=C_*.
		\]
		Equations \eqref{eq:intro-reduced} and
		\eqref{eq:first-level-infimum} now give
		\begin{equation}
			c_1
			=
			\frac{p-1}{2p}
			C_*^{-1/(p-1)}.
			\label{eq:first-level-sharp-value}
		\end{equation}
		Let \(w_0\in\mathcal S\) be the positive normalized constant. Then
		\[
		\mathcal D(w_0)=C_*,
		\qquad
		\Psi(w_0)=c_1.
		\]
		The radial rescaling in Lemma~\ref{lem:radial-reduction} transforms
		\(w_0\) into the positive constant solution \(a_0\). Hence,
		\[
		E(a_0)
		=
		\Psi(w_0)
		=
		c_1.
		\]
	}
	
	{
		Finally, let \(v\in H\setminus\{0\}\) be any critical point of \(E|_H\)
		satisfying
		\(
		E(v)=c_1.
		\)
		{Testing the equation with $v$ gives}
		\(
		\|v\|^2=\mathcal D(v).
		\)
		Define
		\(
		w={v}/{\|v\|}\in\mathcal S.
		\)
		The \(2p\)-homogeneity of \(\mathcal D\) yields
		\[
		\mathcal D(w)
		=
		\frac{\mathcal D(v)}{\|v\|^{2p}}
		=
		\|v\|^{2-2p}.
		\]
		Therefore,
		\[
		\begin{aligned}
			\Psi(w)
			&=
			\frac{p-1}{2p}
			\mathcal D(w)^{-1/(p-1)}=
			\frac{p-1}{2p}\|v\|^2=
			E(v)
			=
			c_1.
		\end{aligned}
		\]
		By the definition of \(\Psi\), the identity
		\(\Psi(w)=E(v)=c_1\), and
		\eqref{eq:first-level-sharp-value}, one has
		\[
		\frac{p-1}{2p}
		\mathcal D(w)^{-1/(p-1)}
		=
		c_1
		=
		\frac{p-1}{2p}
		C_*^{-1/(p-1)}.
		\]
		Therefore,
		\(
		\mathcal D(w)=C_*.
		\)
		Moreover, \(w\in\mathcal S\), and hence \(\|w\|=1\).
		Equation~\eqref{eq:sharp-combined-inequality} now gives
		\[
		\begin{aligned}
			C_*
			=
			\mathcal D(w)
			&\leq
			C_{\mathrm{HLS}}(n,\mu)
			\|w\|_{L^{2^*}(\mathbb S^n)}^{2p}\leq
			C_{\mathrm{HLS}}(n,\mu)
			S_n^{-p}\|w\|^{2p}=
			C_*.
		\end{aligned}
		\]
		Thus equality holds throughout \eqref{eq:sharp-combined-inequality}.
		In particular, equality holds in the sharp Sobolev inequality.
		The equality cases of the Sobolev inequality
		{\cite[Theorem~3.1]{FrankLieb2012} and \cite{Beckner1993}} give
		\[
		w(\xi)
		=
		C(1-a\cdot\xi)^{-(n-2)/2}
		\]
		up to an overall sign, where \(|a|<1\). The \(G\)-invariance of \(w\)
		again implies
		\(
		a=0.
		\)
		Hence, \(w\) and \(v\) are constant. Substitution into
		\eqref{eq:intro-sphere} gives
		\(
		v= a_0
	\) or \(
		v=-a_0.
		\)
	}
\end{proof}
	\begin{proof}[Proof of Theorem~\ref{thm:levels-intro}]
		Proposition~\ref{prop:minimax-growth-bounds} gives
		\[
		C_1j^{2p/(p-1)}
		\leq c_j
		\leq C_2j^{2p/(p-1)},
		\qquad j\geq2.
		\]
		Together with Lemma~\ref{lem:positive-invariant}, which gives \(c_1=E(a_0)<\infty\),
		this shows that every \(c_j\) is finite. Moreover,
		\(
		c_j\longrightarrow\infty\)
		 as \(j\longrightarrow\infty.
		\)
		By Lemma~\ref{lem:psi-ps}, the functional
		\(\Psi\in C^1(\mathcal S)\) is even, bounded away from zero, and
		satisfies the Palais--Smale condition. Hence
		Lemma~\ref{lem:symmetric-minimax} shows that every \(c_j\) is a
		critical value of \(\Psi\). Thus, for each \(j\geq1\), there exists
		\(w_j\in\mathcal S\) such that
		\[
		d\Psi(w_j)=0,
		\qquad
		\Psi(w_j)=c_j.
		\]
		Set
		\[
		v_j
		=
		\mathcal D(w_j)^{-1/(2p-2)}w_j.
		\]
		Lemma~\ref{lem:radial-reduction} gives
		\[
		(E|_H)'(v_j)=0,
		\qquad
		E(v_j)=\Psi(w_j)=c_j.
		\]
		It remains to verify criticality in the full space. Let \(dg\) be
		the normalized Haar measure on \(G\), and let
		\(\phi\in H^1(\mathbb S^n)\). Define
		\[
		\overline\phi(\xi)
		=
		\int_G\phi(g\xi)\,dg.
		\]
		Then \(\overline\phi\in H\). Since \(v_j\) is \(G\)-invariant and
		the action of \(G\) preserves \(E\),
		\[
		E'(v_j)[\phi]
		=
		E'(v_j)[\phi(g\,\cdot)]
		\qquad(g\in G).
		\]
		Averaging this identity over \(G\) yields
		\[
		E'(v_j)[\phi]
		=
		E'(v_j)[\overline\phi]
		=
		0.
		\]
		Therefore, \(v_j\) is a weak solution of
		\eqref{eq:intro-sphere} in the full space.
		By Lemma~\ref{lem:radial-reduction} and Lemma~\ref{lem:positive-invariant}, the critical set of \(\Psi\) at \(c_1\) is
		\(
		K_{c_1}=\{w_0,-w_0\},
		\)
		where \(w_0\) is the positive normalized constant. Hence
		\(
		i(K_{c_1})=1.
		\)
		If \(c_2=c_1\), the plateau conclusion in Lemma~\ref{lem:symmetric-minimax}, applied with
		\(j=1\) and \(s=1\), would give
		\(
		i(K_{c_1})\geq2,
		\)
		a contradiction. Therefore, \(c_2>c_1\).
		Since \(\Gamma_{j+1}\subset\Gamma_j\), the sequence \((c_j)\) is
		nondecreasing. Consequently,
		\[
		c_j\geq c_2>c_1,
		\qquad j\geq2.
		\]
		Lemma~\ref{lem:positive-invariant} also shows that every nonzero
		one-signed solution in \(H\) has energy \(c_1\). Therefore, every
		invariant solution at a level \(c_j\), \(j\geq2\), changes sign.
		Finally, since \((c_j)\) is nondecreasing and \(c_j\to\infty\),
		it contains a strictly increasing subsequence. The corresponding
		critical points yield infinitely many distinct pairs of
		\(G\)-invariant sign-changing solutions.
	\end{proof}

	\section{Regularity and nodal geometry}
	\label{sec:regularity}
	\noindent
	
	{We first prove regularity and describe the zeros of the orbit profile.
		The proof of Theorem~\ref{thm:nodal-intro} is completed in
		Section~\ref{sec:hessian-nodal}.}
	
	\begin{lemma}
		\label{lem:regularity}
		Let {\(v\in H^1(\mathbb S^n)\)} be a weak solution of \eqref{eq:intro-sphere}. Then
		\(v\in L^\infty(\mathbb S^n)\). Moreover, there exists
		\(\gamma\in(0,1)\) such that
		\begin{equation}
			J_\mu(|v|^p)\in C^{0,\gamma}(\mathbb S^n),
			\qquad
			v\in C^{2,\gamma}(\mathbb S^n).
			\label{eq:regularity-result}
		\end{equation}
	\end{lemma}
	
	{
		\begin{proof}
			In the chart determined by \(\pi\), let \(u=Pv\).
			Proposition~\ref{prop:conformal-identities} shows that \(u\) is a
			Euclidean weak solution. By the rotational invariance of the equation, rotating the omitted
			pole shows that the same argument applies in every stereographic
			chart. Set
			\[
			I_u(x)
			:=
			\int_{\mathbb R^n}
			\frac{|u(y)|^p}{|x-y|^\mu}\,dy.
			\]
			The regularity result
			\cite[Theorem~2.3 and Appendix~A]{AlmutairiGluck2026}
			gives, for some \(\alpha\in(0,1)\),
			\[
			u\in C_{\mathrm{loc}}^{2,\alpha}(\mathbb R^n),
			\qquad
			I_u\in C_{\mathrm{loc}}^{0,\alpha}(\mathbb R^n).
			\]
			Moreover, the conformal change of variables gives
			\[
			I_u(x)
			=
			\Omega(x)^{\mu/2}
			J_\mu(|v|^p)(\pi^{-1}(x)).
			\]
			Thus \(v\) and \(J_\mu(|v|^p)\) have the stated local H\"older
			regularity in every chart. Choose finitely many coordinate
			neighborhoods whose closures are contained in stereographic charts
			and which cover \(\mathbb S^n\). On each of them, the conformal factor
			is smooth and bounded away from zero. Taking the minimum of their
			H\"older exponents gives, for some \(\gamma\in(0,1)\),
			\[
			J_\mu(|v|^p)\in C^{0,\gamma}(\mathbb S^n),
			\qquad
			v\in C^{2,\gamma}(\mathbb S^n).
			\]
			The compactness of \(\mathbb S^n\) also gives
			\(v\in L^\infty(\mathbb S^n)\).
		\end{proof}
	}
	
	{For the remainder of this section, assume $p\geq2$ and fix a
		nonzero $G$-invariant weak solution $v\in H$.}
	For $v(\xi)=f(t)$, define
	\begin{equation}
		V_f(t)=J_\mu(|v|^p)(\cos t\,\omega,\sin t\,\theta).
		\label{eq:profile-potential}
	\end{equation}
	The value in \eqref{eq:profile-potential} is independent of $\omega$ and
	$\theta$. By \eqref{eq:regularity-result}, it is continuous and nonnegative.
	If $v\neq0$, strict positivity of the kernel gives $V_f(t)>0$
	for every $t\in[0,\pi/2]$.
	
	\begin{lemma}
		\label{lem:profile-equation}
		The profile \(f\) of the fixed nonzero \(G\)-invariant weak
		solution \(v\) satisfies
		\begin{equation}
			-(\rho f')'+\frac{n(n-2)}4\rho f
			=\rho V_f|f|^{p-2}f
			\quad\hbox{on }(0,\pi/2),
			\label{eq:profile-equation}
		\end{equation}
		and
		\begin{equation}
			f'(0)=f'(\pi/2)=0.
			\label{eq:profile-neumann}
		\end{equation}
	\end{lemma}
	
	\begin{proof}
		By Lemma~\ref{lem:reduced-laplacian},
		\[
		\Delta_{\mathbb S^n}v
		=
		\rho^{-1}(\rho f')'.
		\]
		Substituting this identity and the definition of \(V_f\) into
		\eqref{eq:intro-sphere} gives
		\[
		-\rho^{-1}(\rho f')'
		+
		\frac{n(n-2)}{4}f
		=
		V_f|f|^{p-2}f
		\quad\hbox{on }(0,\pi/2).
		\]
		Multiplication by \(\rho\) yields
		\eqref{eq:profile-equation}.
		
		It remains to verify the endpoint conditions. Fix
		\(\omega\in\mathbb S^{k-1}\) and define
		\[
		F(y)
		=
		v\bigl(\sqrt{1-|y|^2}\,\omega,y\bigr)
		\]
		for \(y\in\mathbb R^\ell\) near \(0\). The \(O(\ell)\)-invariance of
		\(v\) gives \(F(-y)=F(y)\). Differentiating this identity at \(y=0\)
		yields \(DF(0)=0\). For every
		\(\theta\in\mathbb S^{\ell-1}\),
		\[
		f(t)=F(\sin t\,\theta),
		\qquad
		f'(t)
		=
		DF(\sin t\,\theta)[\cos t\,\theta].
		\]
		Evaluating at \(t=0\) gives \(f'(0)=DF(0)[\theta]=0\).
		To treat the other endpoint, fix
		\(\theta\in\mathbb S^{\ell-1}\) and define
		\[
		G(x)
		=
		v\bigl(x,\sqrt{1-|x|^2}\,\theta\bigr)
		\]
		for \(x\in\mathbb R^k\) near \(0\). The \(O(k)\)-invariance of \(v\)
		gives \(G(-x)=G(x)\), and hence \(DG(0)=0\). For every
		\(\omega\in\mathbb S^{k-1}\),
		\[
		f(t)=G(\cos t\,\omega),
		\qquad
		f'(t)
		=
		DG(\cos t\,\omega)[-\sin t\,\omega].
		\]
		Evaluating at \(t=\pi/2\) gives
		\(f'(\pi/2)=DG(0)[-\omega]=0\). Therefore,
		\eqref{eq:profile-neumann} holds.
	\end{proof}
	
	\begin{lemma}
		\label{lem:simple-zeros}
		Assume \(p\geq2\). A nonzero solution profile has no zero at \(0\) or \(\pi/2\). Every zero
		in \((0,\pi/2)\) is simple, and the number of zeros is finite.
	\end{lemma}
	
	\begin{proof}
		{In the open orbit interval $\rho>0$, so \eqref{eq:profile-equation}
			rewrites as}
		\[
		f''+\frac{\rho'}{\rho}f'-c(t)f=0,\qquad c(t)=\frac{n(n-2)}{4}-V_f(t)|f(t)|^{p-2}.
		\]
		The coefficient \(c\) is continuous by the regularity established
		above and the assumption \(p\geq2\), with the convention
		\(|f|^0\equiv1\) when \(p=2\).
			Uniqueness for a linear ODE with continuous coefficients
			\cite[Chapter~1]{CoddingtonLevinson1955} implies that
			$f(t_0)=f'(t_0)=0$ at an interior point forces $f\equiv0$ on
			$(0,\pi/2)$, and then on its closure by continuity. Every interior
			zero of a nontrivial profile is therefore simple.
		
		We exclude endpoint zeros next. Suppose that \(f(0)=0\).
		Equation~\eqref{eq:profile-equation} becomes
		\[
		(\rho f')'=\rho c(t)f.
		\]
		The coefficient \(c\) is bounded near \(0\), and
		\eqref{eq:profile-neumann} gives \(f'(0)=0\). Integration over
		\((0,t)\) yields
		\[
		\rho(t)f'(t)
		=
		\int_0^t\rho(s)c(s)f(s)\,ds.
		\]
		For sufficiently small \(t>0\), there exist constants {\(C_1,C_2>0\)}
		such that
		\[
		{C_1}t^{\ell-1}
		\leq
		\rho(t)
		\leq
		{C_2}t^{\ell-1}.
		\]
		Define \(M(t)=\max_{0\leq s\leq t}|f(s)|\). The preceding estimates give
		\[
		\begin{aligned}
			|f'(t)|
			&\leq
			\frac{C}{\rho(t)}
			\int_0^t\rho(s)|f(s)|\,ds \leq
			Ct^{1-\ell}M(t)
			\int_0^t s^{\ell-1}\,ds \leq
			CtM(t).
		\end{aligned}
		\]
		Since \(f(0)=0\), integration over \((0,s)\), with \(0\leq s\leq t\),
		gives
		\[
		|f(s)|
		\leq
		\int_0^s|f'(r)|\,dr
		\leq
		CM(t)\int_0^s r\,dr
		\leq
		Ct^2M(t).
		\]
		Taking the maximum over \(s\in[0,t]\) yields
		\[
		M(t)\leq Ct^2M(t).
		\]
		Choose \(t>0\) such that \(Ct^2<1\). Then \(M(t)=0\), and hence
		\(f\equiv0\) on \([0,t]\). Interior uniqueness extends this identity to
		\([0,\pi/2]\), which contradicts the assumption that \(f\) is nonzero.
		Therefore, \(f(0)\neq0\).
		
		To treat the other endpoint, set \(r=\pi/2-t\),
		\(\widetilde f(r)=f(\pi/2-r)\), and
		\(\widetilde\rho(r)=\rho(\pi/2-r)\). Near \(r=0\),
		\[
		{\widetilde C_1}r^{k-1}
		\leq
		\widetilde\rho(r)
		\leq
		{\widetilde C_2}r^{k-1}.
		\]
		Applying the preceding argument to \(\widetilde f\), together with
		\(f'(\pi/2)=0\) from \eqref{eq:profile-neumann}, gives
		\(f(\pi/2)\neq0\).
		
		Finally, suppose that the interior zeros have an accumulation point
		\(t_*\in(0,\pi/2)\). For a sequence of distinct zeros \(t_m\to t_*\),
		continuity gives \(f(t_*)=0\), while differentiability gives
		\[
		f'(t_*)
		=
		\lim_{m\to\infty}
		\frac{f(t_m)-f(t_*)}{t_m-t_*}
		=
		0.
		\]
		This contradicts the simplicity of the interior zeros. Accumulation at
		\(0\) or \(\pi/2\) is excluded by \(f(0)\neq0\),
		\(f(\pi/2)\neq0\), and continuity. An infinite subset of the compact
		interval \([0,\pi/2]\) has an accumulation point. Therefore, \(f\) has
		only finitely many zeros.
	\end{proof}

	\section{The Hessian and nodal domain bounds}
	\label{sec:hessian-nodal}
	\noindent
	
Assume \(p\geq2\) throughout this section, and set
\(
H^{-1}(\mathbb S^n)
:=
\bigl(H^1(\mathbb S^n)\bigr)^*.
\)
At \(p=2\), the convention \(|v|^0=1\) is used in the
second-variation formulas. Recall that \(\langle\cdot,\cdot\rangle_A\) denotes the inner product
associated with the energy norm in \eqref{eq:intro-operators},
and that \(\langle\cdot,\cdot\rangle_H\) is its restriction to \(H\).
The weak realization of the conformal Laplacian
\(A\) defined in \eqref{eq:conformal-laplacian} is
\[
A:H^1(\mathbb S^n)\longrightarrow H^{-1}(\mathbb S^n),
\]
characterized by
\begin{equation}
	\langle A\phi,\psi\rangle_{H^{-1},H^1}
	=
	\langle\phi,\psi\rangle_A,
	\qquad
	\phi,\psi\in H^1(\mathbb S^n).
	\label{eq:A-weak-realization}
\end{equation}
The second variation is given by the following lemma.
	
	\begin{lemma}
		\label{lem:twice-differentiability}
		If $p\geq2$, then $\mathcal D$ and $E$ are of class $C^2$ on
		$H^1(\mathbb S^n)$.  For $v,\phi,\psi\in H^1(\mathbb S^n)$,
		\begin{align}
			\mathcal D''(v)[\phi,\psi]={}&2p(p-1)\int_{\mathbb S^n}
			J_\mu(|v|^p)|v|^{p-2}\phi\psi\,dV\nonumber\\
			&+2p^2\iint_{\mathbb S^n\times\mathbb S^n}
			\frac{|v(\xi)|^{p-2}v(\xi)\phi(\xi)
				|v(\eta)|^{p-2}v(\eta)\psi(\eta)}
			{|\xi-\eta|^\mu}\,dV_\xi dV_\eta.
			\label{eq:D-second-bilinear}
		\end{align}
	\end{lemma}
	
\begin{proof}
	Since $pq_\mu=2^*$, the HLS and H\"older inequalities give
	\[
	\left|
	\int_{\mathbb S^n}
	J_\mu(|v|^p)|v|^{p-2}\phi\psi\,dV
	\right|
	\leq
	C\|v\|_{2^*}^{2p-2}
	\|\phi\|_{2^*}\|\psi\|_{2^*}
	\]
	and
	\[
	\left|
	\int_{\mathbb S^n}
	J_\mu(|v|^{p-2}v\phi)
	|v|^{p-2}v\psi\,dV
	\right|
	\leq
	C\|v\|_{2^*}^{2p-2}
	\|\phi\|_{2^*}\|\psi\|_{2^*}.
	\]
	Let ${\mathcal B_v}$ denote the bilinear form on the
	right-hand side of~\eqref{eq:D-second-bilinear}. The estimates above give
	\(
	\|{\mathcal B_v}\|
	\leq C\|v\|^{2p-2}
	\)
	in the bilinear operator norm.
	
	We next prove the continuity of
	$v\mapsto{\mathcal B_v}$.
	Let $v_m\to v$ in $H^1(\mathbb S^n)$. Then
	$v_m\to v$ in $L^{2^*}(\mathbb S^n)$, and the power maps satisfy
	\[
	|v_m|^p\to|v|^p
	\quad\hbox{in }L^{q_\mu}(\mathbb S^n),\quad\quad
	|v_m|^{p-2}v_m\to|v|^{p-2}v
	\quad\hbox{in }L^{2^*/(p-1)}(\mathbb S^n).
	\]
	If $p>2$, one also has
\(
|v_m|^{p-2}\to|v|^{p-2}\)
in \(L^{2^*/(p-2)}(\mathbb S^n).
\)
For $p=2$, this factor is identically one under the convention
\(|v|^0\equiv1\). Applying the HLS and H\"older
inequalities to the differences gives
\(
\|
{\mathcal B_{v_m}}
-{\mathcal B_v}
\|\longrightarrow0.
\)
The same power-map estimates also give the continuity of
\(
\mathcal D':
H^1(\mathbb S^n)\longrightarrow H^{-1}(\mathbb S^n).
\)
For smooth $v$ and $h$, differentiation along $v+th$ gives
\[
\mathcal D'(v+h)-\mathcal D'(v)
-{\mathcal B_v}[h,\cdot]
=
\int_0^1
\left(
{\mathcal B_{v+th}}
-{\mathcal B_v}
\right)[h,\cdot]\,dt.
\]
Approximation in $H^1(\mathbb S^n)$ extends this identity to every
$v,h\in H^1(\mathbb S^n)$. Consequently,
\[
\begin{aligned}
	&\left\|
	\mathcal D'(v+h)-\mathcal D'(v)
	-{\mathcal B_v}[h,\cdot]
	\right\|_{(H^1(\mathbb S^n))^*}
	\leq
	\sup_{0\leq t\leq1}
	\left\|
	{\mathcal B_{v+th}}
	-{\mathcal B_v}
	\right\|
	\|h\|
	=o(\|h\|).
\end{aligned}
\]
The last relation follows from the continuity of
\(v\mapsto\mathcal B_v\) in the bilinear operator norm.
Thus
\(
\mathcal D''(v)={\mathcal B_v}\) and 
\(
\mathcal D\in C^2.
\)
Formula~\eqref{eq:intro-energy} then gives $E\in C^2$.
\end{proof}
	
	At a critical point $v$, the Hessian quadratic form is
	\begin{align}
		Q_v(\phi)={}&\|\phi\|^2
		-(p-1)\int_{\mathbb S^n}J_\mu(|v|^p)|v|^{p-2}\phi^2\,dV\nonumber\\
		&-p\iint_{\mathbb S^n\times\mathbb S^n}
		\frac{|v(\xi)|^{p-2}v(\xi)\phi(\xi)
			|v(\eta)|^{p-2}v(\eta)\phi(\eta)}
		{|\xi-\eta|^\mu}\,dV_\xi dV_\eta.
		\label{eq:hessian-form}
	\end{align}
	We use two operator realizations of the Hessian form. The
	energy-space realization
	\[
	T_v:H^1(\mathbb S^n)\longrightarrow H^1(\mathbb S^n)
	\]
	is defined, by the Riesz representation theorem, through
	\begin{equation}
		\langle T_v\phi,\psi\rangle_A
		=
		E''(v)[\phi,\psi],
		\qquad
		\phi,\psi\in H^1(\mathbb S^n).
		\label{eq:energy-space-linearization}
	\end{equation}
	\begin{lemma}
		\label{lem:second-variation}
		Let $p\geq2$ and let $v\in H^1(\mathbb S^n)$ be a weak solution.
		Formula \eqref{eq:hessian-form} holds for every
		$\phi\in H^1(\mathbb S^n)$. Moreover, there exists a compact self-adjoint operator
		\(K_v:H^1(\mathbb S^n)\to H^1(\mathbb S^n)\) such that
		\[
		T_v=I-K_v.
		\]
		Consequently, \(v\) has finite Morse index and nullity.
	\end{lemma}
	
	\begin{proof}
		{Differentiating the energy \eqref{eq:intro-energy} gives}
			\[
			E''(v)[\phi,\psi]=\langle\phi,\psi\rangle_A
			-\frac{1}{2p}\mathcal D''(v)[\phi,\psi].
			\]
			Lemma~\ref{lem:twice-differentiability} therefore yields
			\eqref{eq:hessian-form}. By Lemma~\ref{lem:regularity},
		for some \(\gamma\in(0,1)\),
		\[
		v\in C^{2,\gamma}(\mathbb S^n),
		\qquad
		J_\mu(|v|^p)\in C^{0,\gamma}(\mathbb S^n).
		\]
		Since \(\mathbb S^n\) is compact, both functions are bounded. Hence,
		\begin{equation}
			a_v:=|v|^{p-2}v,
			\qquad
			V_v:=J_\mu(|v|^p)|v|^{p-2}
			\label{eq:linearized-coefficients}
		\end{equation}
		belong to \(L^\infty(\mathbb S^n)\).
		
		Let \(\phi_m\rightharpoonup0\) in \(H^1(\mathbb S^n)\). The
		{Rellich compact embedding on the compact sphere
			\cite{Aubin1998}} gives
		\[
		\phi_m\longrightarrow0
		\quad\text{in }L^2(\mathbb S^n).
		\]
		For every \(\psi\in H^1(\mathbb S^n)\), H\"older's inequality and the
		continuous embedding
		\(H^1(\mathbb S^n)\hookrightarrow L^2(\mathbb S^n)\) yield
		\[
		\begin{aligned}
			\left|
			\int_{\mathbb S^n}V_v\phi_m\psi\,dV
			\right|
			&\leq
			\|V_v\|_{L^\infty(\mathbb S^n)}
			\|\phi_m\|_{L^2(\mathbb S^n)}
			\|\psi\|_{L^2(\mathbb S^n)}\\
			&\leq
			C\|V_v\|_{L^\infty(\mathbb S^n)}
			\|\phi_m\|_{L^2(\mathbb S^n)}
			\|\psi\|_{H^1(\mathbb S^n)}.
		\end{aligned}
		\]
		Consequently,
		\begin{equation}
			\sup_{\|\psi\|_{H^1(\mathbb S^n)}\leq1}
			\left|
			\int_{\mathbb S^n}V_v\phi_m\psi\,dV
			\right|
			\longrightarrow0.
			\label{eq:local-hessian-compactness}
		\end{equation}
		For \(K(\xi,\eta)=|\xi-\eta|^{-\mu}\), rotational invariance and
		\eqref{eq:constant-state} give
		\[
		\sup_\xi\int_{\mathbb S^n}K(\xi,\eta)\,dV_\eta
		=
		\sup_\eta\int_{\mathbb S^n}K(\xi,\eta)\,dV_\xi
		=
		\kappa_0.
		\]
		The Cauchy--Schwarz inequality with respect to the measure
		\(K(\xi,\eta)\,dV_\eta\), followed by Fubini's theorem, yields
		\begin{equation}
			\|J_\mu(g)\|_{L^2(\mathbb S^n)}
			\leq
			\kappa_0\|g\|_{L^2(\mathbb S^n)},
			\qquad
			g\in L^2(\mathbb S^n).
			\label{eq:riesz-L2-bound}
		\end{equation}
		Moreover,
		\[
		\|a_v\phi_m\|_{L^2(\mathbb S^n)}
		\leq
		\|a_v\|_{L^\infty(\mathbb S^n)}
		\|\phi_m\|_{L^2(\mathbb S^n)}
		\longrightarrow0.
		\]
		Applying \eqref{eq:riesz-L2-bound}, H\"older's inequality, and the
		continuous embedding
		\(H^1(\mathbb S^n)\hookrightarrow L^2(\mathbb S^n)\), one obtains
		\[
		\begin{aligned}
			\left|
			\int_{\mathbb S^n}
			J_\mu(a_v\phi_m)a_v\psi\,dV
			\right|
			&\leq
			\|J_\mu(a_v\phi_m)\|_{L^2(\mathbb S^n)}
			\|a_v\psi\|_{L^2(\mathbb S^n)}\\
			&\leq
			C\|a_v\phi_m\|_{L^2(\mathbb S^n)}
			\|a_v\|_{L^\infty(\mathbb S^n)}
			\|\psi\|_{H^1(\mathbb S^n)}.
		\end{aligned}
		\]
		Hence,
		\begin{equation}
			\sup_{\|\psi\|_{H^1(\mathbb S^n)}\leq1}
			\left|
			\int_{\mathbb S^n}
			J_\mu(a_v\phi_m)a_v\psi\,dV
			\right|
			\longrightarrow0.
			\label{eq:nonlocal-hessian-compactness}
		\end{equation}
		By the Riesz representation theorem, define
		\(K_v:H^1(\mathbb S^n)\to H^1(\mathbb S^n)\) through the energy
		inner product by
		\[
		\begin{aligned}
			\langle K_v\phi,\psi\rangle_A
			={}&
			(p-1)\int_{\mathbb S^n}V_v\phi\psi\,dV+
			p\int_{\mathbb S^n}
			J_\mu(a_v\phi)a_v\psi\,dV.
		\end{aligned}
		\]
		Equations \eqref{eq:local-hessian-compactness} and
		\eqref{eq:nonlocal-hessian-compactness} imply
		\[
		\phi_m\rightharpoonup0
		\ \text{in }H^1(\mathbb S^n)
		\quad\Longrightarrow\quad
		K_v\phi_m\longrightarrow0
		\ \text{in }H^1(\mathbb S^n).
		\]
		Thus, \(K_v\) is compact.
		
		The multiplication term is symmetric in \(\phi\) and \(\psi\). The
		symmetry of the Riesz kernel gives
		\[
		\int_{\mathbb S^n}J_\mu(f)g\,dV
		=
		\int_{\mathbb S^n}fJ_\mu(g)\,dV,
		\]
		so the nonlocal term is also symmetric. Therefore, \(K_v\) is
		{self-adjoint} with respect to the energy inner product.
		
	By the definitions of \(T_v\) and \(K_v\),
	\[
	\langle T_v\phi,\psi\rangle_A
	=
	\langle (I-K_v)\phi,\psi\rangle_A,
	\qquad
	\phi,\psi\in H^1(\mathbb S^n).
	\]
	Hence,
	\(
	T_v=I-K_v.
	\)
	Since \(K_v\) is compact and self-adjoint, the spectral theorem
	shows that \(K_v\) has only finitely many eigenvalues greater than
	\(1\), counted with multiplicity, and that its eigenspace
	corresponding to the eigenvalue \(1\) is finite-dimensional.
	Therefore, \(T_v\) has only finitely many negative eigenvalues and
	\[
	\dim\ker T_v<\infty.
	\]
	Hence the Morse index and nullity of \(v\) are finite.
	\end{proof}
		
We next introduce the $L^2(\mathbb S^n)$-realization of the same
Hessian form. Using the coefficients introduced in~\eqref{eq:linearized-coefficients}, define
		\begin{equation}
			L_v\phi
			=
			A\phi-(p-1)V_v\phi-p\,a_vJ_\mu(a_v\phi),
			\qquad
			\operatorname{dom}(L_v)=H^2(\mathbb S^n).
			\label{eq:L2-linearization}
		\end{equation}
	The boundedness of $a_v$ and $V_v$, together with
	\eqref{eq:riesz-L2-bound}, shows that the two lower-order terms in
	\eqref{eq:L2-linearization} define bounded operators on
	$L^2(\mathbb S^n)$. Thus
	\[
	L_v=A+B_v,
	\qquad
	B_v\phi:=-(p-1)V_v\phi-p\,a_vJ_\mu(a_v\phi).
	\]
	The symmetry identity used above shows that the nonlocal part of
	$B_v$ is self-adjoint. Hence $B_v$ is bounded and self-adjoint on
	$L^2(\mathbb S^n)$. Consequently, $L_v$ is self-adjoint on
	$L^2(\mathbb S^n)$ with
	$\operatorname{dom}(L_v)=H^2(\mathbb S^n)$
	\cite[Chapter~V, \S4, Theorem~4.3]{Kato1995}.
	
	The operator \(L_v\) also has compact resolvent. For a sufficiently
	large \(c>0\), the resolvent identity is
	\[
	(L_v+c)^{-1}=(A+c)^{-1}
	\big(I+B_v(A+c)^{-1}\big)^{-1}.
	\]
	The first factor is compact on \(L^2(\mathbb S^n)\), since
	\(A\) is elliptic on the compact sphere. The second factor is bounded
	when \(c\) is chosen so that \(\|B_v(A+c)^{-1}\|<1\).
	The product is therefore compact.
	
	The closed bilinear form associated with \(L_v\)
	{\cite[Chapter~VI, \S1.6 and \S2.1]{Kato1995}} has domain
	\(H^1(\mathbb S^n)\) and is equal to \(E''(v)\). More precisely,
	\eqref{eq:energy-space-linearization} and
	\eqref{eq:L2-linearization} give
	\begin{equation}
		L_v
		=
		AT_v:
		H^1(\mathbb S^n)
		\longrightarrow
		H^{-1}(\mathbb S^n)
		\label{eq:two-linearizations}
	\end{equation}
	in the weak sense.
	
	If \(T_v\phi=0\), then \eqref{eq:two-linearizations} gives
	\(L_v\phi=0\) in \(H^{-1}(\mathbb S^n)\). Equivalently,
	\[
	A\phi
	=
	(p-1)V_v\phi
	+
	p\,a_vJ_\mu(a_v\phi).
	\]
	The {right-hand} side belongs to \(L^2(\mathbb S^n)\), so elliptic
	{regularity \cite[Section~9.5]{GilbargTrudinger1977}} yields \(\phi\in H^2(\mathbb S^n)\). Therefore,
	\(\phi\in\ker L_v\). Conversely, if \(\phi\in\ker L_v\), then
	\[
	E''(v)[\phi,\psi]
	=
	0,
	\qquad
	\psi\in H^1(\mathbb S^n),
	\]
	and \eqref{eq:energy-space-linearization} gives \(T_v\phi=0\). Hence,
	\[
	\ker T_v
	=
	\ker L_v
	\subset
	H^2(\mathbb S^n).
	\]
	
	The common Hessian form gives the same Morse index for both
	realizations. The kernel identity gives the same nullity. Their
	numerical eigenvalues generally differ because the inner products
	differ. We use \(T_v\) for the Fredholm argument in the energy space
	and \(L_v\) for the spectral analysis in \(L^2(\mathbb S^n)\).
	
	The nodal estimate below relies on the strict positivity of \(J_\mu\).
	
	\begin{lemma}
		\label{lem:riesz-positive}
		For every nonzero {real-valued} function
		\(h\in L^2(\mathbb S^n)\),
		\begin{equation}
			\iint_{\mathbb S^n\times\mathbb S^n}
			\frac{h(\xi)h(\eta)}{|\xi-\eta|^\mu}
			\,dV_\xi\,dV_\eta
			>
			0.
			\label{eq:riesz-positive}
		\end{equation}
	\end{lemma}
	
	\begin{proof}
		{
			Let
			\(
			h=\displaystyle \sum_{m=0}^{\infty}h_m
			\)
			be the orthogonal decomposition of \(h\) into spherical harmonics, where
			\(h_m\) belongs to the space of spherical harmonics of degree \(m\).
			The {{Funk--Hecke} formula} gives
			\begin{equation}
				J_\mu h_m
				=
				\kappa_m h_m,
				\qquad
				\kappa_m
				=
				2^{n-\mu}\pi^{n/2}
				\frac{\Gamma((n-\mu)/2)}{\Gamma(\mu/2)}
				\frac{\Gamma(m+\mu/2)}
				{\Gamma(m+n-\mu/2)}.
				\label{eq:riesz-eigenvalues}
			\end{equation}
			For the normalization in \eqref{eq:riesz-eigenvalues}, use
			{\cite[Corollary~4.3]{FrankLieb2012}
			in dimension $n$ and with $\alpha=\mu/2$.} Since
			$|\xi-\eta|^{-\mu}=2^{-\mu/2}(1-\xi\cdot\eta)^{-\mu/2}$,
			the multiplier there must be multiplied by $2^{-\mu/2}$.
			For nonintegral \(\alpha=\mu/2\), the Gamma recurrence gives
			\[
			(-1)^m\frac{\Gamma(1-\alpha)}{\Gamma(1-\alpha-m)}
			=\frac{\Gamma(m+\alpha)}{\Gamma(\alpha)}.
			\]
			The degree-zero multiplier is fixed by direct integration,
			\[
			\begin{aligned}
				\kappa_0
				&=|\mathbb S^{n-1}|\int_{-1}^{1}
				(2-2t)^{-\mu/2}(1-t^2)^{(n-2)/2}\,dt\\
				&=2^{n-\mu-1}|\mathbb S^{n-1}|
				\int_0^1 u^{(n-\mu)/2-1}(1-u)^{n/2-1}\,du\\
				&=2^{n-\mu}\pi^{n/2}
				\frac{\Gamma((n-\mu)/2)}{\Gamma(n-\mu/2)}.
			\end{aligned}
			\]
			Together with the multiplier ratio from the Funk--Hecke formula,
			this gives \eqref{eq:riesz-eigenvalues}. Integer values of
			\(\mu/2\) follow by continuity.
			The positivity of the Gamma function on
			\((0,\infty)\) gives \(\kappa_m>0\) for every \(m\geq0\). Moreover,
			{the Gamma quotient asymptotic \cite[Section~5.11(iii)]{DLMF}}
			applied to
			\eqref{eq:riesz-eigenvalues} yields
			\[
			\kappa_m
			=
			O(m^{\mu-n})
			\qquad
			\text{as }m\to\infty.
			\]
			The condition \(0<\mu<n\) therefore implies
			\[
			\sup_{m\geq0}\kappa_m<\infty,
			\qquad
			\kappa_m\longrightarrow0.
			\]
			For each integer $M\geq0$, set \(h^{(M)}=\displaystyle \sum_{m=0}^{M}h_m\). The orthogonality of
			distinct spherical harmonic subspaces and
			\eqref{eq:riesz-eigenvalues} give
			\begin{equation}
				\int_{\mathbb S^n}
				J_\mu(h^{(M)})h^{(M)}\,dV
				=
				\sum_{m=0}^{M}
				\kappa_m
				\|h_m\|_{L^2(\mathbb S^n)}^2.
				\label{eq:riesz-finite-expansion}
			\end{equation}
			By~\eqref{eq:riesz-L2-bound}, the operator
			$J_\mu:L^2(\mathbb S^n)\to L^2(\mathbb S^n)$ is bounded.
			Moreover, the boundedness of $\{\kappa_m\}_{m\geq0}$ ensures that
			the series
			\(
			\displaystyle\sum_{m=0}^{\infty}\kappa_m
			\|h_m\|_{L^2(\mathbb S^n)}^2
			\)
			converges. Since $h^{(M)}\to h$ in $L^2(\mathbb S^n)$, passing to
			the limit in~\eqref{eq:riesz-finite-expansion} yields
			\[
			\begin{aligned}
				\iint_{\mathbb S^n\times\mathbb S^n}
				\frac{h(\xi)h(\eta)}{|\xi-\eta|^\mu}
				\,dV_\xi\,dV_\eta
				&=
				\int_{\mathbb S^n}J_\mu(h)h\,dV=
				\sum_{m=0}^{\infty}
				\kappa_m
				\|h_m\|_{L^2(\mathbb S^n)}^2.
			\end{aligned}
			\]
			If \(h\neq0\), the completeness of the spherical harmonic decomposition
			implies that \(h_m\neq0\) for at least one \(m\). Together with
			\(\kappa_m>0\), this gives
			\[
			\sum_{m=0}^{\infty}
			\kappa_m
			\|h_m\|_{L^2(\mathbb S^n)}^2
			>
			0,
			\]
			which proves \eqref{eq:riesz-positive}.
		}
	\end{proof}
	
	{
		For the remainder of this section, fix a nonzero $G$-invariant weak
		solution $v\in H$. The profile has finitely many simple interior
		zeros and no endpoint zeros. The complementary orbit intervals have
		constant sign. Since both sphere factors are connected, their images
		in \(\mathbb S^n\) are connected and separated by the zero layers.
		They are therefore the nodal domains \(\Omega_1,\ldots,\Omega_s\),
		and each is \(G\)-invariant. Put
		\begin{equation}
			v_i=v\,\mathbf1_{\Omega_i},\qquad
			B_{ij}=\iint_{\Omega_i\times\Omega_j}
			\frac{|v(\xi)|^p|v(\eta)|^p}{|\xi-\eta|^\mu}
			\,dV_\xi dV_\eta.
			\label{eq:nodal-pieces}
		\end{equation}
		Let $I_i$ be the corresponding interval on which $f$ has a fixed sign.
		The next lemma proves that the functions in \eqref{eq:nodal-pieces}
		belong to $H$.
	}
	
	\begin{lemma}
		{
			\label{lem:nodal-zero-extension}
			For every $i$, one has $v_i\in H$ and
			\begin{equation}
				\nabla v_i=\mathbf1_{\Omega_i}\nabla v
				\quad\hbox{almost everywhere on }\mathbb S^n.
				\label{eq:nodal-gradient}
			\end{equation}
			If $i\neq j$, then $v_i$ and $v_j$ are orthogonal in the energy scalar
			product.
		}
	\end{lemma}
	
	\begin{proof}
		{
			Let \(I_i=(a_i,b_i)\) be the orbit interval corresponding to
			\(\Omega_i\), and set \(f_i=f\mathbf 1_{I_i}\). By
			Lemma~\ref{lem:simple-zeros},
			\[
			f(a_i)=0\quad\text{if }a_i>0,
			\qquad
			f(b_i)=0\quad\text{if }b_i<\frac{\pi}{2}.
			\]
			For every \(\zeta\in C_c^\infty(0,\pi/2)\), integration by parts on
			\(I_i\) gives
			\[
			\begin{aligned}
				-\int_0^{\pi/2}f_i(t)\zeta'(t)\,dt
				&=
				-\int_{a_i}^{b_i}f(t)\zeta'(t)\,dt\\
				&=
				\int_{a_i}^{b_i}f'(t)\zeta(t)\,dt
				-
				\bigl[f(t)\zeta(t)\bigr]_{a_i}^{b_i}\\
				&=
				\int_{I_i}f'(t)\zeta(t)\,dt.
			\end{aligned}
			\]
			The boundary term vanishes at every interior endpoint by
			Lemma~\ref{lem:simple-zeros}. If \(a_i=0\) or \(b_i=\pi/2\), it vanishes
			from the compact support of \(\zeta\). Therefore,
			\[
			f_i'
			=
			\mathbf 1_{I_i}f'
			\quad\text{in }\mathcal D'(0,\pi/2).
			\]
			It follows that
			\[
			\int_0^{\pi/2}
			\left(
			|f_i'(t)|^2+|f_i(t)|^2
			\right)\rho(t)\,dt
			\leq
			\int_0^{\pi/2}
			\left(
			|f'(t)|^2+|f(t)|^2
			\right)\rho(t)\,dt
			<
			\infty.
			\]
		}
		
		{
			{In any interior product chart, the weak derivative identity
				for $f_i$ above and the absence of angular derivatives give
				$v_i\in H^1(U)$ for every relatively compact interior product
				chart $U$.} Near a singular orbit, $v_i$ either agrees with
			$v$ or vanishes identically, so it belongs to $H^1(U)$ on each
			corresponding singular chart $U$ as well. A finite covering and a subordinate partition of unity therefore yield
			\[
			v_i\in H^1(\mathbb S^n).
			\]
			Both \(v\) and \(\Omega_i\) are \(G\)-invariant, so \(v_i\) is
			\(G\)-invariant. Consequently,
			\[
			v_i\in H.
			\]
			The distributional identity for \(f_i'\), together with the vanishing of
			all angular derivatives, gives \eqref{eq:nodal-gradient}. The singular
			orbits have zero \(n\)-dimensional measure and therefore do not affect
			the almost everywhere identity.
		}
		
		For \(i\neq j\), the nodal domains \(\Omega_i\) and \(\Omega_j\) are
		disjoint. Hence,
		\[
		v_iv_j=0
		\quad\text{and}\quad
		\nabla_{\mathbb S^n}v_i\cdot\nabla_{\mathbb S^n}v_j=0
		\quad\text{almost everywhere on }\mathbb S^n,
		\]
		where the second identity follows from \eqref{eq:nodal-gradient}.
		Therefore,
		\[
		\begin{aligned}
			\langle v_i,v_j\rangle_A
			&=
			\int_{\mathbb S^n}
			\left(
			\nabla_{\mathbb S^n}v_i\cdot\nabla_{\mathbb S^n}v_j
			+
			\frac{n(n-2)}{4}v_iv_j
			\right)dV=
			0.
		\end{aligned}
		\]
		Thus, \(v_i\) and \(v_j\) are orthogonal in the energy inner product.
	\end{proof}
	
	\begin{lemma}
		{
			\label{lem:negative-nodal-space}
			Assume \(p\geq2\). Then \(Q_v\) is negative definite on
			\begin{equation*}
				\mathcal V
				=
				\operatorname{span}\{v_1,\ldots,v_s\}.
			\end{equation*}
		}
	\end{lemma}
	
	\begin{proof}
		{
			Lemma~\ref{lem:nodal-zero-extension} gives
			\(v_i\in H\), so \(v_i\) is an admissible test function
			in the weak formulation of \eqref{eq:intro-sphere}. Moreover, the
			orthogonality established in Lemma~\ref{lem:nodal-zero-extension} gives
			\[
			\langle v,v_i\rangle_A
			=
			\left\langle\sum_{j=1}^{s}v_j,v_i\right\rangle_A
			=
			\|v_i\|^2.
			\]
			Testing \eqref{eq:intro-sphere} with \(v_i\) therefore yields
			\begin{align}
				\|v_i\|^2
				&=
				\int_{\mathbb S^n}
				J_\mu(|v|^p)|v|^{p-2}vv_i\,dV\nonumber\\
				&=
				\int_{\Omega_i}
				J_\mu(|v|^p)|v|^p\,dV\nonumber\\
				&=
				\sum_{j=1}^{s}
				\iint_{\Omega_i\times\Omega_j}
				\frac{|v(\xi)|^p|v(\eta)|^p}
				{|\xi-\eta|^\mu}
				\,dV_\xi\,dV_\eta\nonumber\\
				&=
				\sum_{j=1}^{s}B_{ij}, 
				\label{eq:nodal-testing}
			\end{align}
			where $B_{ij}$ is defined in \eqref{eq:nodal-pieces}.
			Let \(\phi=\displaystyle \sum_{i=1}^{s}c_iv_i\in\mathcal V\), where
			\(c_1,\ldots,c_s\in\mathbb R\). The orthogonality of the functions
			\(v_i\) gives
			\[
			\|\phi\|^2
			=
			\sum_{i=1}^{s}c_i^2\|v_i\|^2.
			\]
			Furthermore,
			\[
			\begin{aligned}
				\int_{\mathbb S^n}
				J_\mu(|v|^p)|v|^{p-2}\phi^2\,dV
				&=
				\sum_{i=1}^{s}c_i^2
				\int_{\Omega_i}J_\mu(|v|^p)|v|^p\,dV=
				\sum_{i=1}^{s}c_i^2\|v_i\|^2,
			\end{aligned}
			\]
			where \eqref{eq:nodal-testing} was used in the last equality. On
			\(\Omega_i\),
			\(
			|v|^{p-2}v\phi
			=
			c_i|v_i|^p.
			\)
			Consequently,
			\[
			\int_{\mathbb S^n}
			J_\mu\bigl(|v|^{p-2}v\phi\bigr)
			|v|^{p-2}v\phi\,dV
			=
			\sum_{i,j=1}^{s}c_ic_jB_{ij}.
			\]
			Substituting these identities into \eqref{eq:hessian-form} gives
			\begin{equation}
				Q_v(\phi)
				=
				(2-p)\sum_{i=1}^{s}c_i^2\|v_i\|^2
				-
				p\sum_{i,j=1}^{s}c_ic_jB_{ij}.
				\label{eq:nodal-quadratic}
			\end{equation}
			Define
			\(
			h
			=
			\displaystyle\sum_{i=1}^{s}c_i|v_i|^p.
			\)
			Lemma~\ref{lem:regularity} gives \(v\in L^\infty(\mathbb S^n)\),
			so $h\in L^2(\mathbb S^n)$. The nonzero sets of
			$|v_1|^p,\ldots,|v_s|^p$ are pairwise disjoint, and each has positive
			measure. Hence, \(h\neq0\) whenever
			\((c_1,\ldots,c_s)\neq0\). By the definition of \(B_{ij}\) and Lemma~\ref{lem:riesz-positive}, we obtain
			\[
			\iint_{\mathbb S^n\times\mathbb S^n}
			\frac{h(\xi)h(\eta)}
			{|\xi-\eta|^\mu}
			\,dV_\xi\,dV_\eta=\sum_{i,j=1}^{s}c_ic_jB_{ij}
			>0
			\qquad
			\text{for }(c_1,\ldots,c_s)\neq0.
			\]
			Since \(p\geq2\), the first term in \eqref{eq:nodal-quadratic} is
			nonpositive. The strict positivity just proved makes its second term
			negative. Thus
			\[
			Q_v(\phi)<0
			\qquad
			\text{for every }\phi\in\mathcal V\setminus\{0\}.
			\]
			Thus, \(Q_v\) is negative definite on \(\mathcal V\). Moreover, the functions $v_1,\ldots,v_s$ are linearly independent
			because they have pairwise disjoint nonzero supports. Hence
			\(
			\dim\mathcal V=s.
			\)
		}
	\end{proof}
	
	{
		\begin{proof}[Proof of Theorem~\ref{thm:nodal-intro}]
			Lemma~\ref{lem:regularity} gives
			\(v\in C^{2,\gamma}(\mathbb S^n)\) for some
			\(\gamma\in(0,1)\). Lemma~\ref{lem:simple-zeros} shows that the
			orbit profile has finitely many simple interior zeros and no endpoint
			zeros.
			By Lemma~\ref{lem:simple-zeros}, either $f$ has no zeros, in which
			case $q=0$, or its zeros satisfy $0<t_1<\cdots<t_q<\pi/2$.
			For $q\geq1$, the corresponding level sets are described below.
			The orbit coordinates in \eqref{eq:orbit-coordinates} give
			\[
			\{\xi\in\mathbb S^n:t(\xi)=t_r\}
			=
			\mathbb S^{k-1}(\cos t_r)
			\times
			\mathbb S^{\ell-1}(\sin t_r).
			\]
			Since \(v(\xi)=f(t(\xi))\), it follows that
			\[
			\{v=0\}
			=
			\bigcup_{r=1}^{q}
			\mathbb S^{k-1}(\cos t_r)
			\times
			\mathbb S^{\ell-1}(\sin t_r),
			\]
			which proves \eqref{eq:intro-nodal-set}.
			
			A simple zero changes the sign of \(f\). Hence, \(f\) has a constant
			sign on each of the \(q+1\) complementary intervals, and the signs on
			adjacent intervals are opposite. For \(q\geq1\), denote these intervals
			by
			\[
			I_0=[0,t_1),
			\qquad
			I_r=(t_r,t_{r+1})\quad(1\leq r\leq q-1),
			\qquad
			I_q=(t_q,\pi/2].
			\]
			If \(q=0\), set \(I_0=[0,\pi/2]\).
			For each \(r\), the subset of \(\mathbb S^n\) corresponding to
			\(t\in I_r\) is the continuous image of
			\[
			I_r\times\mathbb S^{k-1}\times\mathbb S^{\ell-1}
			\]
			under the map
			\[
			(t,\omega,\theta)
			\longmapsto
			(\cos t\,\omega,\sin t\,\theta).
			\]
			The interval \(I_r\) is connected, and
			\(\mathbb S^{k-1}\) and \(\mathbb S^{\ell-1}\) are connected under the
			assumption \(k,\ell\geq2\). Therefore, each corresponding subset of
			\(\{v\neq0\}\) is connected. Distinct subsets are separated by the
			nodal layers \(t=t_r\), and their union equals \(\{v\neq0\}\).
			Consequently, these subsets are precisely the nodal domains of \(v\),
			and
			\[
			N_G(v)=q+1.
			\] 
			For the Morse-index estimate, note that
			{
				\(v_1,\ldots,v_s\) are nonzero and pairwise orthogonal in the energy
				inner product by Lemma~\ref{lem:nodal-zero-extension}. They are
				therefore linearly independent, and
				\(
				\dim\mathcal V=s.
				\)
				By Lemma~\ref{lem:negative-nodal-space}, \(Q_v\) is negative definite
				on the \(s\)-dimensional subspace
				\(\mathcal V\subset H\). The definition of the
				restricted Morse index therefore gives
				\[
				s\leq m_G(v).
				\]
				Since \(s=N_G(v)=q+1\),
				\[
				q+1\leq m_G(v),
				\]
				which proves the index inequality in \eqref{eq:intro-nodal-morse}.
			}
			
			We now prove the minimax energy comparison. Put
			\[
			A_i=\|v_i\|^2,\qquad B=(B_{ij})_{1\leq i,j\leq s},\qquad
			\mathbf1=(1,\ldots,1)^T.
			\]
			Equation~\eqref{eq:nodal-testing} states that \(A_i=(B\mathbf1)_i\).
			The proof of Lemma~\ref{lem:negative-nodal-space} also gives
			\[
			y^TBy>0\qquad\text{for every }y\in\mathbb R^s\setminus\{0\}.
			\]
			Hence \(\lambda_{\min}(B)>0\). For \(t=(t_1,\ldots,t_s)\in\mathbb R^s\),
			set \(x_i=|t_i|^p\). The disjoint supports give
			\[
			\left|\sum_i t_i v_i\right|^p=\sum_i x_i|v_i|^p,
			\qquad
			\left\|\sum_i t_i v_i\right\|^2=\sum_i A_i x_i^{2/p}.
			\]
			Therefore,
			\[
			E\left(\sum_i t_i v_i\right)
			=\frac12\sum_i A_i x_i^{2/p}-\frac1{2p}x^TBx.
			\]
			Since \(0<2/p\leq1\), concavity on \([0,\infty)\) gives
			\[
			x_i^{2/p}-1\leq\frac2p(x_i-1).
			\]
			It follows that
			\[
			\begin{aligned}
				E\left(\sum_i t_i v_i\right)-E(v)
				&\leq\frac1p\sum_i A_i(x_i-1)
				-\frac1{2p}\big(x^TBx-\mathbf1^TB\mathbf1\big)\\
				&=\frac1p\mathbf1^TB(x-\mathbf1)
				-\frac1{2p}\big(x^TBx-\mathbf1^TB\mathbf1\big).
			\end{aligned}
			\]
			Completing the quadratic form yields
			\begin{equation}
				E\left(\sum_i t_i v_i\right)-E(v)
				\leq-\frac1{2p}(x-\mathbf1)^TB(x-\mathbf1).
				\label{eq:nodal-energy-gap}
			\end{equation}
			This inequality includes \(p=2\); there concavity is an equality.
			
			The energy unit sphere \(\mathcal V\cap\mathcal S\) is compact and
			has cohomological index \(s\). For any \(w\) on this sphere,
			\(rw\in\mathcal V\) for all \(r>0\). By the radial maximization
			formula and \eqref{eq:nodal-energy-gap},
			\[
			\Psi(w)=\max_{r>0}E(rw)\leq E(v).
			\]
			Taking the supremum on this sphere and then the minimax infimum gives
			\[
			c_s\leq\sup_{\mathcal V\cap\mathcal S}\Psi\leq E(v).
			\]
			If \(1\leq j<s\), choose \(j\) of the nodal pieces and let
			\(\mathcal V_j\) be their span. At least one coordinate of \(x\)
			vanishes for every vector in \(\mathcal V_j\). Thus
			\[
			|x-\mathbf1|^2\geq1,
			\qquad
			E\left(\sum_i t_i v_i\right)
			\leq E(v)-\frac{\lambda_{\min}(B)}{2p}
			\quad\text{on }\mathcal V_j.
			\]
			The sphere \(\mathcal V_j\cap\mathcal S\) has index \(j\), so
			\[
			c_j\leq\sup_{\mathcal V_j\cap\mathcal S}\Psi
			\leq E(v)-\frac{\lambda_{\min}(B)}{2p}<E(v).
			\]
			In particular, if \(E(v)=c_j\), then \(s>j\) would contradict
			this strict inequality. Hence \(N_G(v)=s\leq j\). No isolation
			of the critical point or critical level is used here.
			
			Every sign-changing invariant solution has at least two nodal domains.
			For such a solution, monotonicity of the minimax levels gives
			\[
			c_2\leq c_s\leq E(v).
			\]
			Theorem~\ref{thm:levels-intro} gives a sign-changing solution at \(c_2\).
			Every solution at that level satisfies \(2\leq N_G(v)\leq2\).
			Thus \(c_2\) is the least energy of sign-changing invariant solutions,
			and all invariant solutions at \(c_2\) have exactly two nodal domains.
		\end{proof}
	}
	
	\begin{remark}
		\label{rem:strong-kernel}
		{For $1<p<2$, the coefficient $|v|^{p-2}$ can be singular at
			nodal zeros. The second-variation argument here applies for $p\geq2$.
			The minimax bounds and sign change at every level $c_j$, $j\geq2$, apply
			throughout the parameter range of Theorem~\ref{thm:levels-intro}.}
	\end{remark}
	
	\section{Cohomological minimax labels and the Morse index}
	\label{sec:morse-label}
	\noindent 
	
	{
		This section proves Theorem~\ref{thm:morse-intro}. The index $i(B)$ is
		defined using Alexander--Spanier cohomology with coefficients in
		$\mathbb Z_2$ and the characteristic class of $B\to B/\mathbb Z_2$
		\cite{FadellRabinowitz1978}. For local critical groups we use singular
		cohomology with the same coefficients. The proof below passes between these two theories on regular
		sublevel sets and finite complexes.}
	\begin{lemma}
		\label{lem:index-localization}
		Let \(F\in C^2(\mathcal S)\) be even and satisfy the {Palais--Smale}
		condition, and let \(d_j\) be defined as in
		Lemma~\ref{lem:symmetric-minimax}. Assume \(d_j\in\mathbb R\).
		Suppose that every critical point
		at \(d_j\) is isolated in the entire critical set of \(F\).
		Assume that each constrained Hessian at a critical point of level $d_j$
		is represented by a {self-adjoint} Fredholm operator with finite Morse index. Then there
		exists a critical point \(w\) of \(F\) at level \(d_j\) such that
		\begin{equation}
			\widetilde m(w)
			\leq
			j-1
			\leq
			\widetilde m(w)+\widetilde\nu(w),
			\label{eq:constrained-index-sandwich}
		\end{equation}
		where \(\widetilde m(w)\) and \(\widetilde\nu(w)\) are the Morse index
		and nullity of the constrained Hessian of \(F\) on
		\(T_w\mathcal S\).
	\end{lemma}
	
	{
		\begin{proof}
			Put \(c=d_j\), and write \(K_c=\{F=c,\ dF=0\}\).
			Every sequence in \(K_c\) is a Palais--Smale sequence, so \(K_c\)
			is compact. If it contained infinitely many distinct points, a
			convergent subsequence would contradict isolation at its limit.
			Thus \(K_c\) is a finite union of antipodal pairs.
			We next show that \(c\) is an isolated critical value. Otherwise
			there would be critical points \(w_m\) with
			\[
			F(w_m)\ne c,\qquad F(w_m)\longrightarrow c,\qquad dF(w_m)=0.
			\]
			Palais--Smale gives a subsequence converging to \(w\in K_c\).
			The points \(w_m\) are different from \(w\), contradicting isolation
			in the full critical set. This proves isolation of the level.

			Set $M=\mathcal S/\mathbb Z_2$, let
			$\overline F:M\to\mathbb R$ be induced by $F$, and define
			\(
			M^a=\{\overline F\leq a\},
			~
			F^a=\{F\leq a\}.
			\)
			Let
			\(
			\omega\in H^1(M;\mathbb Z_2)
			\)
			be the Alexander--Spanier characteristic class of the double cover
			$\mathcal S\to M$. For every regular value $a$, we have
			\begin{equation}
				\left.\omega^{j-1}\right|_{M^a}\neq0
				\quad\Longleftrightarrow\quad
				\text{there exists }B\in\Gamma_j
				\text{ such that }B\subset F^a.
				\label{eq:compact-detection}
			\end{equation}
			Indeed, if $B\in\Gamma_j$ and $B\subset F^a$, naturality of the
			characteristic class gives
			\(
			\omega_B^{j-1}
			=
			\iota_B^*
			\left(
			\left.\omega^{j-1}\right|_{M^a}
			\right)\neq0,
			\)
			where $\iota_B:B/\mathbb Z_2\to M^a$ is induced by inclusion. This
			proves the implication from right to left.
			Conversely, suppose that
			\(
			\alpha
			:=
			\left.\omega^{j-1}\right|_{M^a}
			\neq0.
			\)
			The regular sublevel $M^a$ is a metrizable, locally contractible
			Hilbert manifold with boundary, or is empty. Hence Alexander--Spanier
			and singular cohomology agree naturally on $M^a$; see
			\cite[Chapter~6, Section~9, Theorem~1]{Spanier1966}.
			{Over $\mathbb Z_2$, the universal coefficient theorem yields a
				homology class on which the singular class corresponding to $\alpha$
				evaluates nontrivially. By the definition of singular homology, this
				homology class is represented by a finite singular cycle.} Gluing its finitely many simplices along
			cancelling faces gives a finite $\Delta$-complex $X$ and a map
			$f:X\to M^a$ such that
			\(
			f^*\alpha\neq0.
			\)
			Pull back the double cover $F^a\to M^a$ along $f$, and let
			$\widetilde X\to X$ be the resulting compact double cover. Its natural
			map $g:\widetilde X\to F^a$ is equivariant. Therefore,
			\(
			B:=g(\widetilde X)
			\)
			is a compact symmetric subset of $F^a$. {If $\widetilde\omega$ denotes
				the characteristic class of the pullback cover $\widetilde X\to X$, then
				naturality first gives $\widetilde\omega^{j-1}=f^*\alpha$. Applying
				naturality again to the equivariant map $g$ gives}
			\(
			f^*\alpha
			=
			\overline g^{\,*}\left(\omega_B^{j-1}\right)\neq0,
			\)
			where $\overline g:X\to B/\mathbb Z_2$ is induced by $g$. Thus
			$\omega_B^{j-1}\neq0$, so $i(B)\geq j$ and $B\in\Gamma_j$. This proves
			\eqref{eq:compact-detection}.
			
			Since \(c\) is isolated, choose \(\varepsilon>0\) so that
			\(c\) is the only critical value in
			\([c-2\varepsilon,c+2\varepsilon]\) and \(c\pm\varepsilon\) are
			regular values. The definition of \(d_j\) and
			\eqref{eq:compact-detection} give
			\[
			\left.\omega^{j-1}\right|_{M^{c+\varepsilon}}\neq0,
			\qquad
			\left.\omega^{j-1}\right|_{M^{c-\varepsilon}}=0.
			\]
			The exact cohomology sequence of the pair therefore yields
			\begin{equation}
				H^{j-1}
				\left(M^{c+\varepsilon},M^{c-\varepsilon};\mathbb Z_2\right)
				\neq0.
				\label{eq:relative-nonzero}
			\end{equation}
			Let \(z_1,\ldots,z_m\) be the critical points of \(\overline F\) at
			level \(c\). {Choose pairwise disjoint isolating neighborhoods
				$U_r$ of $z_r$ and set}
			{
				\[
				C^q(\overline F,z_r)=H^q\bigl(M^c\cap U_r,(M^c\cap U_r)
				\setminus\{z_r\};\mathbb Z_2\bigr).
				\]}
			The {Palais--Smale} condition, the isolation of these points
			and the second deformation lemma localize the relative group to their
			local critical groups:
			\begin{equation}
				H^q
				\left(M^{c+\varepsilon},M^{c-\varepsilon};\mathbb Z_2\right)
				\cong
				\bigoplus_{r=1}^m C^q(\overline F,z_r).
				\label{eq:local-critical-group-sum}
			\end{equation}
			{This is the critical-group decomposition from the second
				deformation lemma and excision
				\cite[Chapter~I, Section~4]{Chang1993}; see also
				\cite[Theorem~III.3]{ChangGhoussoub1996} for the local
				Gromoll--Meyer description.} Hence
			\eqref{eq:relative-nonzero} and
			\eqref{eq:local-critical-group-sum} imply
			\[
			C^{j-1}(\overline F,z_r)\neq0
			\]
			for at least one \(r\).
			
			Choose either lift \(w\) of this \(z_r\). Locally the covering map is a
			diffeomorphism, so the constrained Hessians at \(z_r\) and \(w\) have
			the same Morse index and nullity. The Fredholm splitting lemma and the
			shifting theorem give
			\[
			C^q(\overline F,z_r)
			\cong
			C^{q-\widetilde m(w)}(h,0),
			\]
			where {$h$ is the reduced functional on the
				$\widetilde\nu(w)$-dimensional null space, with an isolated critical
				point at $0$}. Consequently these local
			groups vanish for
			\[
			q<\widetilde m(w)
			\quad\text{or}\quad
			q>\widetilde m(w)+\widetilde\nu(w);
			\]
			see \cite[Chapter~I, Section~5]{Chang1993}. Taking \(q=j-1\) proves
			\[
			\widetilde m(w)\leq j-1
			\leq\widetilde m(w)+\widetilde\nu(w),
			\]
			which is \eqref{eq:constrained-index-sandwich}.
		\end{proof}
	}
	
	We now compare the constrained Hessian of \(\Psi\) with the Hessian
	of \(E\). Let \(w\in\mathcal S\) be a critical point of \(\Psi\), and {set
		\[
		t_0:=\mathcal D(w)^{-1/(2p-2)},
		\qquad v:=t_0w,
		\]
		as in \eqref{eq:rescaling-critical}.}
	
	\begin{lemma}
		\label{lem:radial-index-shift}
		At the corresponding critical points \(w\in\mathcal S\) and
		\(v=t_0w\in H\), one has
		\begin{equation}
			m_G(v)=\widetilde m(w)+1,
			\qquad
			\nu_G(v)=\widetilde\nu(w).
			\label{eq:radial-index-shift}
		\end{equation}
	\end{lemma}
	
	\begin{proof}
		For \(z=t\zeta\in H\setminus\{0\}\), where \(t>0\) and
		\(\zeta\in\mathcal S\), the homogeneity of \(\mathcal D\) gives
		\begin{equation}
			E(t\zeta)
			=
			\frac{1}{2}t^2
			-
			\frac{1}{2p}t^{2p}\mathcal D(\zeta).
			\label{eq:polar-energy}
		\end{equation}
		The orthogonal decomposition
		\[
		H
		=
		\operatorname{span}\{w\}
		\oplus
		T_w\mathcal S
		\]
		separates the radial and tangential directions.
		
		Let \(\eta\in T_w\mathcal S\). Differentiating
		\eqref{eq:polar-energy} once in the radial direction and once in the
		tangential direction gives
		\[
		\partial_t\partial_\eta E(t_0w)
		=
		-t_0^{2p-1}
		d(\mathcal D|_{\mathcal S})(w)[\eta].
		\]
		The criticality of \(w\) for \(\Psi|_{\mathcal S}\), together with the
		definition of \(\Psi\), implies
		\(d(\mathcal D|_{\mathcal S})(w)=0\). Hence,
		\[
		\partial_t\partial_\eta E(t_0w)=0.
		\]
		Thus the Hessian of \(E\) at \(v\) has no mixed radial-tangential
		term.
		
		Differentiating \eqref{eq:polar-energy} twice with respect to \(t\) and
		using \(t_0^{2p-2}\mathcal D(w)=1\), which follows from
		\eqref{eq:rescaling-critical}, yields
		\[
		\begin{aligned}
			\partial_t^2E(t_0w)
			&=
			1-(2p-1)t_0^{2p-2}\mathcal D(w)\\
			&=
			-2(p-1)<0.
		\end{aligned}
		\]
		Therefore, the radial direction contributes exactly one negative
		direction and no null direction.
		
		The relation
		\(
		t_0=\mathcal D(w)^{-1/(2p-2)}
		\)
		gives
		\(
		t_0^{2p}
		=
		\mathcal D(w)^{-p/(p-1)}.
		\)
		At the constrained critical point \(w\), the second derivative of
		\(\Psi\) along \(\mathcal S\) is
		\begin{equation}
			\operatorname{Hess}_{\mathcal S}\Psi(w)[\eta,\eta]
			=
			-\frac{1}{2p}
			\mathcal D(w)^{-p/(p-1)}
			\operatorname{Hess}_{\mathcal S}\mathcal D(w)[\eta,\eta].
			\label{eq:psi-tangential-hessian}
		\end{equation}
		Applying \eqref{eq:polar-energy} along a geodesic of \(\mathcal S\)
		through \(w\) with initial velocity \(\eta\) gives
		\begin{equation}
			Q_v(t_0\eta)
			=
			-\frac{1}{2p}
			t_0^{2p}
			\operatorname{Hess}_{\mathcal S}\mathcal D(w)[\eta,\eta].
			\label{eq:E-tangential-hessian}
		\end{equation}
		The identity
		\(
		t_0^{2p}=\mathcal D(w)^{-p/(p-1)}
		\),
		together with
		\eqref{eq:psi-tangential-hessian} and
		\eqref{eq:E-tangential-hessian}, yields
		\[
		Q_v(t_0\eta)
		=
		\operatorname{Hess}_{\mathcal S}\Psi(w)[\eta,\eta].
		\]
		Since \(t_0>0\), the map
		\[
		{S_{t_0}:T_w\mathcal S\longrightarrow T_w\mathcal S},
		\qquad
		{S_{t_0}\eta=t_0\eta},
		\]
		is a linear isomorphism. Polarizing the quadratic-form identity gives
		\[
		E''(v)[{S_{t_0}\eta,S_{t_0}\xi}]
		=
		\operatorname{Hess}_{\mathcal S}\Psi(w)[\eta,\xi],
		\qquad
		\eta,\xi\in T_w\mathcal S.
		\]
		Thus {\(S_{t_0}\)} carries the tangential Hessian of \(E\) to the constrained
		Hessian of \(\Psi\). Their Morse indices and nullities agree. Adding
		the negative radial direction gives
		\[
		m_G(v)=\widetilde m(w)+1,
		\qquad
		\nu_G(v)=\widetilde\nu(w),
		\]
		which proves \eqref{eq:radial-index-shift}.
	\end{proof}
	
	\begin{proof}[Proof of Theorem~\ref{thm:morse-intro}]
		Lemma~\ref{lem:psi-ps} shows that $\Psi$ is even and satisfies the
		Palais--Smale condition. Lemma~\ref{lem:twice-differentiability}, the positivity of $\mathcal D$,
		and~\eqref{eq:intro-reduced} imply that $\Psi\in C^2(\mathcal S)$.
		For each critical point $w$ of $\Psi$ at level $c_j$, set
		\[
		t=\mathcal D(w)^{-1/(2p-2)},
		\qquad
		v=tw.
		\]
		By Lemma~\ref{lem:radial-reduction}, $v$ is a critical point of $E|_H$. Since $E$ is
		$G$-invariant, the principle of symmetric criticality~\cite{Palais1979} implies
		that $v$ is a critical point of $E$ on the full space
		$H^1(\mathbb S^n)$, and hence a $G$-invariant weak solution of~\eqref{eq:intro-sphere}.
		
		Lemma~\ref{lem:second-variation} therefore applies to $v$. Moreover, since $v$ is
		$G$-invariant, the coefficients in the Hessian and the compact
		operator $K_v$ preserve the invariant subspace $H$. Consequently,
		the Hessian of $E|_H$ at $v$ is represented on $H$ by the identity
		minus a compact self-adjoint operator. Hence its representing operator is self-adjoint and Fredholm, with
		finite Morse index and nullity. The radial--tangential decomposition
		and the linear isomorphism established in Lemma~\ref{lem:radial-index-shift} show that the
		constrained Hessian of $\Psi$ on $T_w\mathcal S$ is also represented
		by a self-adjoint Fredholm operator, with finite Morse index and
		nullity.

	By the isolation hypotheses of Theorem~\ref{thm:morse-intro}, every critical point of
	$\Psi$ at level $c_j$ is isolated in the full critical set of
	$\Psi$. Consequently, all the assumptions of Lemma~\ref{lem:index-localization} are
	satisfied with $F=\Psi$ and $d_j=c_j$. Thus, there exists a critical point
		$w_j\in\mathcal S$ at level $c_j$ such that
		\[
		\widetilde m(w_j)
		\leq
		j-1
		\leq
		\widetilde m(w_j)+\widetilde\nu(w_j).
		\]
		Set
		\[
		t_j=\mathcal D(w_j)^{-1/(2p-2)},
		\qquad
		v_j=t_jw_j.
		\]
	By the preceding argument, $v_j$ is a $G$-invariant weak solution
	of~\eqref{eq:intro-sphere} with energy $E(v_j)=c_j$.
		By \eqref{eq:radial-index-shift},
		\[
		m_G(v_j)=\widetilde m(w_j)+1,
		\qquad
		\nu_G(v_j)=\widetilde\nu(w_j).
		\]
		Substitution into the preceding inequalities yields
		\[
		m_G(v_j)
		\leq
		j
		\leq
		m_G(v_j)+\nu_G(v_j).
		\]
		Theorem~\ref{thm:nodal-intro}, together with the preceding estimate, gives
		\[
		N_G(v_j)
		\leq
		m_G(v_j)
		\leq
		j.
		\]
		Finally, if $\nu_G(v_j)=0$, then
		\[
		m_G(v_j)\leq j\leq m_G(v_j),
		\]
		and hence $m_G(v_j)=j$.
	\end{proof}
	
	\begin{corollary}
		\label{cor:plateau}
		If
		\(
		c_j=c_{j+1}=\cdots=c_{j+s},
		\)
		then the critical set at this common level has cohomological index at
		least \(s+1\). Consequently, a finite union of antipodal critical pairs
		can correspond to at most one minimax label.
	\end{corollary}
	
	\begin{proof}
		Let \(c=c_j=\cdots=c_{j+s}\), and denote the critical set of \(\Psi\) at
		level \(c\) by \(K_c\). Lemma~\ref{lem:symmetric-minimax} gives
		\(i(K_c)\geq s+1\). 
		Suppose that \(K_c\) is a nonempty finite union of antipodal pairs. Then
		\(K_c/\mathbb Z_2\) is a nonempty finite discrete space, and hence
		\[
		H^0(K_c/\mathbb Z_2;\mathbb Z_2)\neq0,
		\qquad
		H^q(K_c/\mathbb Z_2;\mathbb Z_2)=0
		\quad\text{for }q\geq1.
		\]
		The definition of the cohomological index therefore gives \(i(K_c)=1\).
		Combining this identity with \(i(K_c)\geq s+1\) yields \(s=0\).
	\end{proof}

The isolation hypothesis of Theorem~\ref{thm:morse-intro} is verified at $c_1$ in
Proposition~\ref{prop:constant-spectrum}, but remains open for the higher minimax levels.
Accordingly, the asymptotic conclusions in Remark~\ref{rem:conditional-index-asymptotics} are conditional.

	\section{Full-space Hessian and angular decomposition}
	\label{sec:full-hessian}
	\noindent

		We decompose the full Hessian into the angular sectors of
		$O(k)\times O(\ell)$. Each sector is determined by a radial quadratic
		form. We identify the symmetry-generated kernel and use the Jacobi
		spectrum to locate the sectors that contribute to the index and
		nullity. We then prove the transverse index bound and compute the
		three specified low blocks for minimizers odd under factor interchange
		when $k=\ell$ and $\mu$ is small.
		
		{
		For $a,b\geq0$, let $\mathcal H_a(\mathbb S^{k-1})$ and
		$\mathcal H_b(\mathbb S^{\ell-1})$ denote the spaces of spherical
		harmonics of degrees $a$ and $b$, respectively.
	}
	
	\subsection{Symmetry-generated kernel modes}
	
	For a Euclidean solution \(u\), the associated linearized operator is
	\begin{align*}
		\mathcal L_u\phi
		={}&
		-\Delta\phi
		-(p-1)
		\left(
		\int_{\mathbb R^n}
		\frac{|u(y)|^p}{|x-y|^\mu}\,dy
		\right)
		|u(x)|^{p-2}\phi(x)\\
		&-
		p
		\left(
		\int_{\mathbb R^n}
		\frac{|u(y)|^{p-2}u(y)\phi(y)}
		{|x-y|^\mu}\,dy
		\right)
		|u(x)|^{p-2}u(x).
	\end{align*}
	
	\begin{proposition}
		{
			\label{prop:conformal-blocks}
			Assume \(p\geq2\), and let \(v\) be a nonzero \(G\)-invariant solution
			with orbit profile \(f\). Define
			\begin{equation}
				\begin{aligned}
					g_{1,0}(t)
					&=
					-\sin t\,f'(t)
					-\frac{n-2}{2}\cos t\,f(t),\\
					g_{0,1}(t)
					&=
					\cos t\,f'(t)
					-\frac{n-2}{2}\sin t\,f(t).
				\end{aligned}
				\label{eq:conformal-block-profiles}
			\end{equation}
			Neither radial profile is identically zero. For every
			\(Y_1\in\mathcal H_1(\mathbb S^{k-1})\) and
			\(Z_1\in\mathcal H_1(\mathbb S^{\ell-1})\), the functions
			\begin{equation}
				g_{1,0}(t)Y_1(\omega),
				\qquad
				g_{0,1}(t)Z_1(\theta)
				\label{eq:conformal-block-fields}
			\end{equation}
			belong to the kernel of the full spherical Hessian. These two angular
			sectors therefore contribute at least
			\(k+\ell=n+1\) dimensions to the full nullity.
		}
		
		If \(f\) is nonconstant, then
		\begin{equation}
			f'(t)Y_1(\omega)Z_1(\theta)
			\label{eq:rotation-block-fields}
		\end{equation}
		is a nonzero kernel field for every pair of nonzero {degree-one}
		harmonics \(Y_1\) and \(Z_1\). The \((1,1)\) sector consequently contributes at least
		\(k\ell\) additional dimensions to the full nullity. Hence, the tangent space
		\(\mathcal K(v)\) to the conformal orbit of \(v\) satisfies
		\begin{equation}
			\dim\mathcal K(v)
			=
			\begin{cases}
				n+1, & f\text{ is constant},\\
				n+1+k\ell, & f\text{ is nonconstant}.
			\end{cases}
			\label{eq:conformal-orbit-dimension}
		\end{equation}
	\end{proposition}
	
	\begin{proof}
		For \(q\in\mathbb R^{n+1}\), define
		\(\zeta_q(\xi)=q\cdot\xi\). {The tangential gradient of $\zeta_q$ is
			$q-\zeta_q\xi$. Taking the tangential part of its Euclidean derivative
			along a tangent vector gives}
		\[
		\nabla_{\mathbb S^n}^2\zeta_q
		=
		-\zeta_q g_{\mathbb S^n}.
		\]
		Hence, for \(X_q=\nabla_{\mathbb S^n}\zeta_q\), {the Lie derivative
			of the metric along $X_q$ satisfies}
		\[
		\mathcal L_{X_q}g_{\mathbb S^n}
		=
		2\nabla_{\mathbb S^n}^2\zeta_q
		=
		-2\zeta_qg_{\mathbb S^n}.
		\]
		Thus \(X_q\) is a conformal vector field, with infinitesimal
		length factor \(-\zeta_q\). Let \(F_s\) be its flow and let
		\(F_s^*g_{\mathbb S^n}=\lambda_s^2g_{\mathbb S^n}\).
		The metric identity gives
		\[
		\partial_s\log\lambda_s(\xi)=-\zeta_q(F_s(\xi)),
		\qquad \lambda_0=1.
		\]
		For \(x_s=F_s(\xi)\) and \(y_s=F_s(\eta)\), one has
		\[
		\begin{aligned}
			\frac{d}{ds}|x_s-y_s|^2
			&=2(x_s-y_s)\cdot
			\bigl(q-(q\cdot x_s)x_s-q+(q\cdot y_s)y_s\bigr)\\
			&=-\bigl(\zeta_q(x_s)+\zeta_q(y_s)\bigr)|x_s-y_s|^2.
		\end{aligned}
		\]
	Integration and the initial value at \(s=0\) give
	\[
	|F_s(\xi)-F_s(\eta)|^2
	=
	\lambda_s(\xi)\lambda_s(\eta)|\xi-\eta|^2,
	\qquad
	F_s^*dV
	=
	\lambda_s^n dV.
	\]
		By Lemma~\ref{lem:regularity}, \(v\in C^{2,\gamma}(\mathbb S^n)\) for some
		\(\gamma\in(0,1)\). In particular, \(v\) is bounded on
		\(\mathbb S^n\), and the conformal orbit map below is differentiable
		as an \(H^1(\mathbb S^n)\)-valued map.
		Put \(h=(n-2)/2\) and \(v_s(\xi)=\lambda_s(\xi)^h v(F_s(\xi))\).
		The conformal covariance of \(A\) gives \(\|v_s\|=\|v\|\).
		In the nonlocal integral the exponent at each integration variable is
		\[
		hp+\frac\mu2-n
		=\frac{2n-\mu}{2}+\frac\mu2-n=0.
		\]
		A change of variables therefore gives
		\(\mathcal D(v_s)=\mathcal D(v)\), and hence \(E(v_s)=E(v)\).
		The action is a linear isomorphism of the energy space for each
		\(s\), so it maps critical points to critical points.
		Differentiating the conformal action on \(v\) gives the tangent field
		\begin{equation}
			\Phi_q
			=
			\nabla_{\mathbb S^n}\zeta_q
			\cdot
			\nabla_{\mathbb S^n}v
			-
			\frac{n-2}{2}\zeta_qv.
			\label{eq:spherical-conformal-field}
		\end{equation}
		{For the curve of solutions $v_s$ defined above, we have
		\[
		\left.\frac{d}{ds}\right|_{s=0}v_s=\Phi_q,
		\qquad E'(v_s)=0.
		\]
			Differentiating the latter identity at $s=0$ and
			using Lemma~\ref{lem:twice-differentiability} give
			\[
			E''(v)[\Phi_q,\phi]=0,
			\qquad
			\phi\in H^1(\mathbb S^n).
			\]}
		For \(1\leq i\leq k\), {take}
		\(\zeta_q=\xi_i=\cos t\,\omega_i\). Then
		\(\partial_t\xi_i=-\sin t\,\omega_i\), and
		\eqref{eq:spherical-conformal-field} yields
		\[
		{\Phi_q}
		=
		\left(
		-\sin t\,f'(t)
		-\frac{n-2}{2}\cos t\,f(t)
		\right)\omega_i
		=
		g_{1,0}(t)\omega_i.
		\]
		For \(1\leq j\leq\ell\), {take}
		\(\zeta_q=\xi_{k+j}=\sin t\,\theta_j\). The identity
		\(\partial_t\xi_{k+j}=\cos t\,\theta_j\) gives
		\[
		{\Phi_q}
		=
		\left(
		\cos t\,f'(t)
		-\frac{n-2}{2}\sin t\,f(t)
		\right)\theta_j
		=
		g_{0,1}(t)\theta_j.
		\]
		The coordinate functions
		\(\omega_1,\ldots,\omega_k\) and
		\(\theta_1,\ldots,\theta_\ell\) form bases of
		\(\mathcal H_1(\mathbb S^{k-1})\) and
		\(\mathcal H_1(\mathbb S^{\ell-1})\), respectively. This proves the
		kernel assertion in \eqref{eq:conformal-block-fields}.
		
		Suppose that \(g_{1,0}\equiv0\). Then
		\[
		f'(t)
		=
		-\frac{n-2}{2}\cot t\,f(t),
		\]
		and integration gives
		\[
		f(t)
		=
		C(\sin t)^{-(n-2)/2}.
		\]
		The boundedness of \(v\) at the singular orbit \(t=0\) forces \(C=0\),
		contradicting \(v\neq0\). Hence \(g_{1,0}\not\equiv0\).
		Similarly, \(g_{0,1}\equiv0\) would imply
		\[
		f(t)
		=
		C(\cos t)^{-(n-2)/2},
		\]
		and boundedness at \(t=\pi/2\) would again force \(C=0\). Therefore,
		\(g_{0,1}\not\equiv0\). The corresponding angular multiplicities are
		\(k\) and \(\ell\).
		
		{
			For \(1\leq i\leq k\) and \(1\leq j\leq\ell\), consider the infinitesimal
			rotation
			\[
			R_{i,k+j}
			=
			\xi_i\partial_{\xi_{k+j}}
			-
			\xi_{k+j}\partial_{\xi_i}.
			\]
			This vector field is tangent to \(\mathbb S^n\). From
			\[
			R_{i,k+j}(\sin^2t)
			=
			2\sin t\cos t\,\omega_i\theta_j
			\]
			one obtains \(R_{i,k+j}t=\omega_i\theta_j\), and hence
			\[
			R_{i,k+j}v
			=
			f'(t)\omega_i\theta_j.
			\]
			Rotations preserve both the equation and the energy. Differentiating the
			rotated family of solutions at the identity gives
			\[
			E''(v)[R_{i,k+j}v,\phi]
			=
			0,
			\qquad
			\phi\in H^1(\mathbb S^n).
			\]
		The regularity stated above shows that the fields in
		\eqref{eq:rotation-block-fields} belong to \(H^1(\mathbb S^n)\).
		}
		
		If \(f\) is nonconstant, then \(f'\not\equiv0\). The products
		\(\omega_i\theta_j\), \(1\leq i\leq k\), \(1\leq j\leq\ell\), form a
		basis of
		\[
		\mathcal H_1(\mathbb S^{k-1})
		\otimes
		\mathcal H_1(\mathbb S^{\ell-1}),
		\]
		so the resulting \(k\ell\) kernel fields are linearly independent.
		The \((1,1)\) sector is orthogonal to the \((1,0)\) and \((0,1)\)
		sectors.

		{
			Rotations within $\mathbb R^k$ or $\mathbb R^\ell$ annihilate the
			\(G\)-invariant function \(v\).
			{Every conformal vector field $X$ on the round sphere is
				the sum of a rotational field and a gradient $\nabla\zeta_q$.
				Indeed, if $\mathcal L_X g_{\mathbb S^n}=2\sigma g_{\mathbb S^n}$, differentiating the
				scalar curvature transformation
				\cite[Section~2]{LeeParker1987} gives
				$-2(n-1)(\Delta_{\mathbb S^n}+n)\sigma=0$, as the round scalar
				curvature is constant. Thus $\sigma=-\zeta_q$ for some
				$q\in\mathbb R^{n+1}$ by the first spherical eigenspace
				\cite[Chapter~1]{DaiXu2013}. The field $X-\nabla\zeta_q$ is
				Killing and therefore generates rotations.}
			Every rotational field is a linear combination of
			rotations within the two factors and the mixed rotations
			\(R_{i,k+j}\). The within-factor rotations annihilate \(v\).
			Thus the gradient conformal fields and the mixed rotations span
			\(\mathcal K(v)\). When \(f\) is constant, the mixed rotation
			fields vanish, leaving \(n+1\) independent directions. When
			\(f\) is nonconstant, the mixed rotations give \(k\ell\)
			additional independent directions in the orthogonal \((1,1)\)
			sector. This proves \eqref{eq:conformal-orbit-dimension}.
		}
	\end{proof}
	
	\subsection{Angular block decomposition and endpoint domains}
	\label{sec:angular-domains}
	
	{
		For nonnegative integers \(a,b\), write
		\[
		d_{k,a}=\dim\mathcal H_a(\mathbb S^{k-1}),
		\qquad
		d_{\ell,b}=\dim\mathcal H_b(\mathbb S^{\ell-1})
		\]
		for the corresponding multiplicities.}

	Let \(\Pi_{a,b}\) be the orthogonal projection in the angular variables
	onto
	\[
	\mathcal H_a(\mathbb S^{k-1})
	\otimes
	\mathcal H_b(\mathbb S^{\ell-1}).
	\]
	{The spaces
		\[
		\mathcal H_a(\mathbb S^{k-1})
		\otimes
		\mathcal H_b(\mathbb S^{\ell-1})
		\]
		are pairwise inequivalent irreducible representations of
		$G=O(k)\times O(\ell)$; for a positive harmonic degree
		on an \(O(2)\) factor, reflection interchanges the two Fourier modes
		of opposite frequencies, which together form one irreducible
		\(O(2)\)-representation. The degree-zero space is the
		one-dimensional trivial representation; see
		\cite[Chapter~1]{DaiXu2013}.}
	Since $J_\mu$ is $G$-equivariant, it preserves every angular sector.
	Hence,
	\[
	J_\mu\Pi_{a,b}=\Pi_{a,b}J_\mu.
	\]
	Consequently,
	\[
	\Pi_{a,b}J_\mu\Pi_{a',b'}
	=
	0
	\qquad
	\text{whenever }(a,b)\neq(a',b').
	\]
	
	{
		Set
		\[
		\lambda_{k,a}=a(a+k-2),
		\qquad
		\lambda_{\ell,b}=b(b+\ell-2).
		\]
		The radial form domain associated with the bidegree \((a,b)\) is
	\begin{align}
		\mathcal X_{a,b}
		=\bigg\{g\in H_{\mathrm{loc}}^1(0,\pi/2):&
		\int_0^{\pi/2}\bigg(|g'(t)|^2\rho(t)\nonumber\\
		&\quad+
		\left(
		\frac{\lambda_{k,a}}{\cos^2 t}
		+
		\frac{\lambda_{\ell,b}}{\sin^2 t}
		+
		\frac{n(n-2)}{4}
		\right)
		|g(t)|^2\rho(t)\bigg)\,dt<\infty
		\bigg\}.
		\label{eq:block-domain}
	\end{align}
		For fixed harmonics of unit $L^2$ norm with respect to the standard
		surface measure,
		\(Y_a\in\mathcal H_a(\mathbb S^{k-1})\) and
		\(Z_b\in\mathcal H_b(\mathbb S^{\ell-1})\), consider
		\[
		\phi(t,\omega,\theta)
		=
		g(t)Y_a(\omega)Z_b(\theta),
		\qquad
		g\in\mathcal X_{a,b}.
		\]
	}
	For almost every \((t,s)\in (0,\pi/2) \times (0,\pi/2)\), define
	\begin{align}
		K_{\mu,a,b}(t,s)
		={}&
		\int_{\mathbb S^{k-1}}
		\int_{\mathbb S^{\ell-1}}
		\int_{\mathbb S^{k-1}}
		\int_{\mathbb S^{\ell-1}}
		\frac{
			Y_a(\omega)Z_b(\theta)
			Y_a(\omega')Z_b(\theta')
		}{
			\bigl(
			2-2\cos t\cos s\,\omega\cdot\omega'
			-2\sin t\sin s\,\theta\cdot\theta'
			\bigr)^{\mu/2}
		}\nonumber\\
		&\hspace{25mm}
		dV_{\theta'}\,dV_{\omega'}\,dV_\theta\,dV_\omega.
		\label{eq:projected-kernel}
	\end{align}
	{For fixed $t\ne s$ the angular integral operator commutes with $G$.
		Schur's lemma on the complexified product harmonic space gives a scalar
		operator, so $K_{\mu,a,b}$ is independent of the chosen real unit harmonics.}
	Moreover,
	\[
	K_{\mu,a,b}(t,s)
	=
	K_{\mu,a,b}(s,t),
	\]
	and the integral in \eqref{eq:projected-kernel} is finite whenever
	\(t\neq s\).
	
	At a regular orbit, the angular dimension is \(n-1\). Near coincident
	angular points, the singular part of the integral has the form
	\[
	\int_0^\varepsilon r^{n-2-\mu}\,dr,
	\]
	which diverges when \(\mu\geq n-1\). Thus,
	\(K_{\mu,a,b}(t,t)\) need not be finite. The diagonal
	\[
	\{(t,s)\in(0,\pi/2)^2:t=s\}
	\]
	has two-dimensional Lebesgue measure zero. We may therefore assign
	arbitrary values to \(K_{\mu,a,b}\) there. We use its almost everywhere
	equivalence class, or the bounded projected operator
	\(\Pi_{a,b}J_\mu\Pi_{a,b}\).
	
	For the unit harmonics fixed above, define the radial operator by
	\[
	(J_{\mu,a,b}g)(t)
	=\int_0^{\pi/2}K_{\mu,a,b}(t,s)g(s)\rho(s)\,ds.
	\]
	It is the bounded operator on \(L^2((0,\pi/2),\rho\,dt)\) determined by
	\[
	J_\mu(gY_aZ_b)=(J_{\mu,a,b}g)Y_aZ_b.
	\]
	The integral formula holds almost everywhere. By
	\eqref{eq:riesz-L2-bound} and the unit normalization of the harmonics,
	\[
	\|J_{\mu,a,b}g\|_{L^2(\rho\,dt)}
	\leq\kappa_0\|g\|_{L^2(\rho\,dt)}.
	\]
	
	{
		Write
		\[
		v(\cos t\,\omega,\sin t\,\theta)=f(t),
		\]
		and set \(V_f(t)=J_\mu(|v|^p)(\cos t\,\omega,\sin t\,\theta)\).
		Then \(V_v(t)=V_f(t)|f(t)|^{p-2}\), with \(V_v\) as in
		\eqref{eq:linearized-coefficients}. On
		\(\mathcal X_{a,b}\), define
		\begin{align}
			Q_{a,b}(g)
			={}&
			\int_0^{\pi/2}
			\bigg[
			|g'(t)|^2
			+
			\bigg(
			\frac{\lambda_{k,a}}{\cos^2t}
			+
			\frac{\lambda_{\ell,b}}{\sin^2t}
			+
			\frac{n(n-2)}{4}\nonumber\\
			&\hspace{34mm}
			-(p-1)V_f(t)|f(t)|^{p-2}
			\bigg)
			|g(t)|^2
			\bigg]
			\rho(t)\,dt\nonumber\\
			&-
			p
			\int_0^{\pi/2}
			\int_0^{\pi/2}
			K_{\mu,a,b}(t,s)
			|f(t)|^{p-2}f(t)g(t)\nonumber\\
			&\hspace{36mm}\times
			|f(s)|^{p-2}f(s)g(s)
			\rho(t)\rho(s)\,dt\,ds.
			\label{eq:radial-block-form}
		\end{align}
		The last term is understood as the quadratic form of the bounded
		operator \(\Pi_{a,b}J_\mu\Pi_{a,b}\). Indeed, for
		\(h(t)=|f(t)|^{p-2}f(t)g(t)\), the \(L^2(\mathbb S^n)\)-boundedness of \(J_\mu\),
		{from \eqref{eq:riesz-L2-bound}}, gives
		\[
		\left|
		\left\langle
		J_\mu\bigl(hY_aZ_b\bigr),
		hY_aZ_b
		\right\rangle_{L^2(\mathbb S^n)}
		\right|
		\leq
		\kappa_0
		\|hY_aZ_b\|_{L^2(\mathbb S^n)}^2.
		\]
		The same estimate applied to \(|hY_aZ_b|\) gives
		\[
		\iint_{\mathbb S^n\times\mathbb S^n}
		\frac{|h(t)Y_a(\omega)Z_b(\theta)|
			|h(s)Y_a(\omega')Z_b(\theta')|}{|\xi-\eta|^\mu}
		\,dV_\xi dV_\eta
		\leq\kappa_0\|hY_aZ_b\|_2^2<\infty.
		\]
		Fubini's theorem therefore yields \eqref{eq:radial-block-form}.
		No finite pointwise value of \(K_{\mu,a,b}(t,t)\) is required.
	}
	
	{
		The space \(\mathcal X_{a,b}\) is complete with respect to the weighted
		form norm in \eqref{eq:block-domain}. It is precisely the radial
		\(H^1(\mathbb S^n)\)-form domain of the \((a,b)\) sector. The density at the singular
		orbits requires a logarithmic cutoff when the corresponding normal
		dimension equals two.
	}
	
	First replace \(g\) by its truncation
	\[
	g_T=\max\{-T,\min\{g,T\}\}.
	\]
	Then \(g_T\to g\) in the form norm as \(T\to\infty\). It is therefore
	enough to consider bounded radial profiles. Near either singular orbit,
	let \(r\) denote the normal distance and let \(q\) denote the
	corresponding normal dimension, with \(q=\ell\) near \(t=0\) and
	\(q=k\) near \(t=\pi/2\).
	
	{
		If $q>2$, take a cutoff that is zero for $r\leq\varepsilon$ and one
		for $r\geq2\varepsilon$, with derivative supported in
		$\{\varepsilon<r<2\varepsilon\}$. It satisfies
		\[
		\int
		r^{q-1}|\chi_\varepsilon'(r)|^2\,dr
		=
		O(\varepsilon^{q-2}).
		\]
		If \(q=2\), define
		\[
		\chi_\varepsilon(r)
		=
		\begin{cases}
			0,
			& r\leq\varepsilon^2,\\[1mm]
			\displaystyle
			\frac{\log(r/\varepsilon^2)}
			{|\log\varepsilon|},
			& \varepsilon^2<r<\varepsilon,\\[3mm]
			1,
			& r\geq\varepsilon.
		\end{cases}
		\]
		For this piecewise smooth cutoff,
		\[
		\int_{\varepsilon^2}^{\varepsilon}
		r|\chi_\varepsilon'(r)|^2\,dr
		=\frac{1}{|\log\varepsilon|}\longrightarrow0.
		\]
		Smoothing the corners on comparable annuli preserves the bound
		$C/|\log\varepsilon|$, which is sufficient for the density argument.
		The endpoint weights are comparable to \(r^{q-1}\,dr\). The remaining
		derivative, \(L^2(\mathbb S^n)\), and angular potential terms converge by dominated
		convergence. For bounded \(g\), the derivative term satisfies
		\[
		\begin{aligned}
			\int\rho\,|((1-\chi_\varepsilon)g)'|^2\,dt
			&\leq 2\int_{\{\chi_\varepsilon\ne1\}}\rho\,|g'|^2\,dt
			+2\|g\|_\infty^2\int\rho\,|\chi_\varepsilon'|^2\,dt\longrightarrow0.
		\end{aligned}
		\]
		The angular potential terms are bounded by the corresponding
		integrable terms in \eqref{eq:block-domain}, multiplied by
		\((1-\chi_\varepsilon)^2\). They also tend to zero. Applying the
		cutoff at both endpoints and then smoothing in the interior proves that
		\[
		C_c^\infty(0,\pi/2)
		\quad\text{is dense in }\mathcal X_{a,b}.
		\]
		The same approximation identifies \(\mathcal X_{a,b}\) with the radial
		form domain of the corresponding \(H^1(\mathbb S^n)\) fields across the
		singular orbits.
	}
	
	{
		The local part of the block is the Friedrichs realization in
		\(L^2((0,\pi/2),\rho(t)\,dt)\) of
		\[
		A_{a,b}g
		=
		-\rho^{-1}(\rho g')'
		+
		\left(
		\frac{\lambda_{k,a}}{\cos^2t}
		+
		\frac{\lambda_{\ell,b}}{\sin^2t}
		+
		\frac{n(n-2)}{4}
		\right)g,
		\]
		with form domain \(\mathcal X_{a,b}\). The coefficients \(V_v\) and \(a_v\) are bounded.
		Multiplication by \(V_v\) is therefore bounded and self-adjoint on
		\(L^2(\rho\,dt)\). Moreover, with \(M_{a_v}\) denoting
		multiplication by the radial profile of \(a_v\),
		\[
		\|M_{a_v}J_{\mu,a,b}M_{a_v}g\|_{L^2(\rho\,dt)}
		\leq\kappa_0\|a_v\|_\infty^2\|g\|_{L^2(\rho\,dt)}.
		\]
		The symmetry of \(J_{\mu,a,b}\) makes this nonlocal operator
		self-adjoint as well. Thus the multiplication and nonlocal terms
		in \eqref{eq:radial-block-form} are bounded self-adjoint
		perturbations.
		Therefore, the full \((a,b)\)-block of \(L_v\) has the same form domain
		\(\mathcal X_{a,b}\). The form domain prescribes the natural Friedrichs behavior at the
		singular endpoints. No zero trace is required when the corresponding
		angular label is zero. Functions in the operator domain still satisfy
		the regularity condition selected by this realization.
	}
	
	\begin{proposition}
		\label{prop:block-diagonalization}
		Let \(p\geq2\), and let \(v\neq0\) be a \(G\)-invariant critical point of \(E\) on \(H^1(\mathbb S^n)\).
		The Hessian operator commutes with the action of
		\(G=O(k)\times O(\ell)\) and is block diagonal with respect to the
		bidegrees \((a,b)\). If \(m_{a,b}(v)\) denotes the negative index of
		the radial form \(Q_{a,b}\) defined in
		\eqref{eq:radial-block-form}, then
		\begin{equation}
			m(v)
			=
			\sum_{a,b\geq0}
			d_{k,a}d_{\ell,b}m_{a,b}(v).
			\label{eq:full-index-sum}
		\end{equation}
		Only finitely many summands in \eqref{eq:full-index-sum} are nonzero.
		Moreover,
		\begin{equation}
			m_{0,0}(v)=m_G(v),
			\qquad
			m(v)\geq m_G(v).
			\label{eq:full-restricted-comparison}
		\end{equation}
		If \(\nu_{a,b}(v)\) denotes the nullity of \(Q_{a,b}\), then
		\begin{equation}
			\nu(v)
			=
			\sum_{a,b\geq0}
			d_{k,a}d_{\ell,b}\nu_{a,b}(v),
			\label{eq:full-nullity-sum}
		\end{equation}
		and only finitely many summands in
		\eqref{eq:full-nullity-sum} are nonzero. Furthermore,
		\begin{equation}
			\nu_{1,0}(v)\geq1,
			\qquad
			\nu_{0,1}(v)\geq1.
			\label{eq:block-conformal-nullity}
		\end{equation}
		If \(v\) is nonconstant, then \(\nu_{1,1}(v)\geq1\).
	\end{proposition}
	
	\begin{proof}
		The coefficients of the local terms in
		\eqref{eq:hessian-form} depend only on \(t\). Moreover, the kernel
		\(|\xi-\eta|^{-\mu}\) is invariant under the diagonal action of
		\(O(n+1)\), and hence under \(G\). It follows that the \(L^2(\mathbb S^n)\)-realization \(L_v\) of the
		full Hessian commutes with the unitary action of \(G\).
		
		{
			Choose orthonormal bases
			\[
			\bigl\{Y_{a,r}\bigr\}_{r=1}^{d_{k,a}}
			\subset\mathcal H_a(\mathbb S^{k-1}),
			\qquad
			\bigl\{Z_{b,s}\bigr\}_{s=1}^{d_{\ell,b}}
			\subset\mathcal H_b(\mathbb S^{\ell-1}).
			\]
			{The spectral decompositions on the two sphere factors
				\cite[Chapter~1]{DaiXu2013}} and Fubini's
			theorem yield
			\[
			\phi
			=
			\sum_{a,b\geq0}
			\sum_{r=1}^{d_{k,a}}
			\sum_{s=1}^{d_{\ell,b}}
			g_{a,b,r,s}(t)
			Y_{a,r}(\omega)Z_{b,s}(\theta)
			\]
			for every \(\phi\in H^1(\mathbb S^n)\). {Parseval's identity and the spherical harmonic}
			eigenvalue equations give
			\[
			\begin{aligned}
				\|\phi\|^2
				=
				\sum_{a,b\geq0}
				\sum_{r=1}^{d_{k,a}}
				\sum_{s=1}^{d_{\ell,b}}
				\int_0^{\pi/2}
				\left[
				|g_{a,b,r,s}'|^2
				+
				\left(
				{\frac{n(n-2)}4}
				+\frac{\lambda_{k,a}}{\cos^2t}
				+\frac{\lambda_{\ell,b}}{\sin^2t}
				\right)
				|g_{a,b,r,s}|^2
				\right]\rho(t)\,dt.
			\end{aligned}
			\]
			{The positivity of \(n(n-2)/4\)} shows that each integral is equivalent to the
			weighted norm in \eqref{eq:block-domain}. Consequently,
			\(g_{a,b,r,s}\in\mathcal X_{a,b}\).
		}
		
		Schur's lemma, applied after complexification, makes the Hessian
			matrix elements zero between distinct bidegrees and scalar in the
			angular variables within each bidegree. Restricting to real-valued
			functions gives the same decomposition. In the
			$L^2(\mathbb S^n)$-realization, the restriction of $L_v$ to the
			$(a,b)$-sector is
			\[
			I_{\mathcal H_a(\mathbb S^{k-1})
				\otimes\mathcal H_b(\mathbb S^{\ell-1})}
			\otimes L_{a,b},
			\]
			where $L_{a,b}$ is the self-adjoint radial operator associated with
			the closed form $Q_{a,b}$ on $\mathcal X_{a,b}$.
		
		The local terms in \eqref{eq:hessian-form}, together with the
		spherical harmonic eigenvalue equations, give the first integral in
		\eqref{eq:radial-block-form}. Applying Fubini's theorem to the
		nonlocal term gives the projected kernel
		\(K_{\mu,a,b}\) defined in \eqref{eq:projected-kernel}. Hence the
		quadratic form associated with \(L_{a,b}\) is precisely
		\(Q_{a,b}\).
		
		Every negative direction of \(Q_{a,b}\) therefore occurs in the full
		Hessian with multiplicity \(d_{k,a}d_{\ell,b}\). Summing these
		multiplicities over all bidegrees gives
		\eqref{eq:full-index-sum}. Lemma~\ref{lem:second-variation} shows that
		the negative spectral space of the full Hessian is {finite-dimensional}.
		Thus, only finitely many terms in \eqref{eq:full-index-sum} are
		nonzero.
		
		The same orthogonal decomposition of the kernel gives
		\eqref{eq:full-nullity-sum}. Lemma~\ref{lem:second-variation} also
		shows that the full kernel is {finite-dimensional}, so only finitely
		many terms in \eqref{eq:full-nullity-sum} are nonzero.
		
		{
		The zero-degree harmonic spaces consist of constant functions.
		The map \(g\mapsto gY_0Z_0\), with \(Y_0\) and \(Z_0\)
		the unit constant harmonics, identifies \(\mathcal X_{0,0}\)
		with \(H\) and carries \(Q_{0,0}\) to \(Q_v|_H\).
		Therefore,
		\(
		m_{0,0}(v)=m_G(v).
		\)
		Equation~\eqref{eq:full-index-sum} then gives
		\[
		m(v)\geq m_G(v),
		\]
		which proves \eqref{eq:full-restricted-comparison}.
		}
		
		Proposition~\ref{prop:conformal-blocks} provides nonzero radial kernel
		modes in the \((1,0)\)- and \((0,1)\)-blocks. Therefore,
		\eqref{eq:block-conformal-nullity} holds. If \(v\) is nonconstant, the
		rotation fields in Proposition~\ref{prop:conformal-blocks} provide a
		nonzero radial kernel mode in the \((1,1)\)-block, and hence
		\(\nu_{1,1}(v)\geq1\).
	\end{proof}
	
	\subsection{Jacobi spectra and localization of negative modes}
	
	\begin{proposition}
		{
			\label{prop:angular-jacobi-spectrum}
			The radial Friedrichs realization of \(A\) in the sector \((a,b)\) has
			{the simple eigenvalues}
			{
				\begin{equation}
					\lambda_{m,a,b}
					=
					(a+b+2m)(a+b+2m+n-1)
					+
					\frac{n(n-2)}4,
					\qquad m\geq0.
					\label{eq:angular-block-spectrum}
				\end{equation}
			}
			The corresponding eigenprofiles are
			\begin{equation}
				g_{m,a,b}(t)
				=
				\gamma_{m,a,b}
				(\cos t)^a(\sin t)^b
				P_m^{(\alpha_{a,b},\beta_{a,b})}(\cos 2t),
				\label{eq:angular-jacobi-profiles}
			\end{equation}
			where
			\[
			\alpha_{a,b}
			=
			b+\frac{\ell-2}{2},
			\qquad
			\beta_{a,b}
			=
			a+\frac{k-2}{2},
			\]
			and \(\gamma_{m,a,b}>0\) is chosen so that
			\(\|g_{m,a,b}\|_{L^2((0,\pi/2),\rho\,dt)}=1\). The family
			\(\{g_{m,a,b}\}_{m\geq0}\) is an orthonormal basis of
			\(L^2((0,\pi/2),\rho\,dt)\) and an orthogonal basis of
			\(\mathcal X_{a,b}\) with respect to the local energy form.
		}
		
		The projected Riesz operator has the same eigenprofiles, with
		\begin{equation}
			J_{\mu,a,b}g_{m,a,b}=\kappa_{a+b+2m}g_{m,a,b}.
			\label{eq:angular-riesz-multipliers}
		\end{equation}
	\end{proposition}
	
	\begin{proof}
		Set
		\[
		w_{a,b}(t)
		=
		(\cos t)^a(\sin t)^b,
		\qquad
		g=w_{a,b}h,
		\]
		and define
		\[
		\rho_{a,b}(t)
		=
		w_{a,b}(t)^2\rho(t)
		=
		(\cos t)^{k-1+2a}
		(\sin t)^{\ell-1+2b}.
		\]
		{The identities}
		\[
		\frac{w_{a,b}'}{w_{a,b}}
		=
		-a\tan t+b\cot t,
		\qquad
		\lambda_{k,a}=a(a+k-2),
		\qquad
		\lambda_{\ell,b}=b(b+\ell-2)
		\]
		give
		\[
		-\rho^{-1}(\rho w_{a,b}')'
		+
		\left(
		\frac{\lambda_{k,a}}{\cos^2t}
		+
		\frac{\lambda_{\ell,b}}{\sin^2t}
		\right)w_{a,b}
		=
		{(a+b)(a+b+n-1)w_{a,b}}.
		\]
		Substituting \(g=w_{a,b}h\) and integrating by parts, we obtain
		\begin{align}
			&\int_0^{\pi/2}
			\left[
			|g'|^2
			+
			\left(
			\frac{\lambda_{k,a}}{\cos^2t}
			+
			\frac{\lambda_{\ell,b}}{\sin^2t}
			+
			\frac{n(n-2)}{4}
			\right)g^2
			\right]\rho(t)\,dt
			\nonumber\\
			=&
			\int_0^{\pi/2}
			\left[
			|h'|^2
			+
			\left(
			{(a+b)(a+b+n-1)}
			+
			\frac{n(n-2)}{4}
			\right)h^2
			\right]\rho_{a,b}(t)\,dt
			\label{eq:angular-ground-state-transform}
		\end{align}
		for every \(h\in C_c^\infty(0,\pi/2)\). Density in
		\(\mathcal X_{a,b}\) extends
		\eqref{eq:angular-ground-state-transform} to the corresponding closed
		form domains.
		
		{
			The multiplication operator
			\[
			U_{a,b}h=w_{a,b}h
			\]
			is unitary from \(L^2((0,\pi/2),\rho_{a,b}\,dt)\) onto
			\(L^2((0,\pi/2),\rho\,dt)\), as
			\[
			\|U_{a,b}h\|_{L^2((0,\pi/2),\rho\,dt)}^2
			=
			\int_0^{\pi/2}|h|^2\rho_{a,b}\,dt.
			\]
			Equation \eqref{eq:angular-ground-state-transform} also identifies the
			associated form domains. Under this unitary transformation, the radial
			operator is conjugate to
			\[
			-\rho_{a,b}^{-1}(\rho_{a,b}h')'
			+
			{(a+b)(a+b+n-1)h}
			+
			\frac{n(n-2)}{4}h.
			\]
		}
		
		The weight \(\rho_{a,b}\) is obtained from the invariant radial weight
		by replacing \(k\) and \(\ell\) with \(k+2a\) and \(\ell+2b\),
		respectively. Proposition~\ref{prop:jacobi-spectrum} therefore applies 
		with
		\[
		\alpha_{a,b}
		=
		\frac{\ell+2b-2}{2},
		\qquad
		\beta_{a,b}
		=
		\frac{k+2a-2}{2}.
		\]
		The derivative part has eigenfunctions
		\[
		P_m^{(\alpha_{a,b},\beta_{a,b})}(\cos 2t)
		\]
		and eigenvalues
		\[
		2m\bigl(2m+n-1+2a+2b\bigr).
		\]
		Adding the constant term in the conjugated operator gives
		\[
		\begin{aligned}
			\lambda_{m,a,b}
			&=
			2m\bigl(2m+n-1+2a+2b\bigr)
			+
			(a+b)(a+b+n-1)
			+
			\frac{n(n-2)}{4}\\
			&=
			(2m+a+b)(2m+a+b+n-1)
			+
			\frac{n(n-2)}{4},
		\end{aligned}
		\]
		{which agrees with \eqref{eq:angular-block-spectrum}.} Applying
		\(U_{a,b}\) gives the profiles in
		\eqref{eq:angular-jacobi-profiles}.
		
		The completeness and orthogonality of the Jacobi polynomials in the
		weighted Jacobi space, together with the unitary map \(U_{a,b}\), show
		that \(\{g_{m,a,b}\}_{m\geq0}\) is an orthonormal basis of
		\(L^2((0,\pi/2),\rho\,dt)\). The spectral theorem for the radial
		Friedrichs operator gives orthogonal completeness in the form domain
		\(\mathcal X_{a,b}\). Finally,
		{\(\lambda_{m+1,a,b}>\lambda_{m,a,b}\) for every
			\(m\geq0\). Together with completeness of the displayed
			eigenprofiles, this shows that each radial eigenspace is
			one-dimensional.}
		
		To identify the Riesz multipliers, extend \(Y_a\) and \(Z_b\)
		to homogeneous harmonic polynomials \(H_a\) and \(H_b\) on
		\(\mathbb R^k\) and \(\mathbb R^\ell\). On the sphere,
		\[
		(\cos t)^aY_a(\omega)=H_a(\xi'),
		\qquad (\sin t)^bZ_b(\theta)=H_b(\xi''),
		\qquad \cos2t=|\xi'|^2-|\xi''|^2.
		\]
		Thus \(g_{m,a,b}Y_aZ_b\) is smooth across both singular orbits.
		Its spherical Laplace eigenvalue, computed above, is
		\((a+b+2m)(a+b+2m+n-1)\). It therefore lies in the spherical
		harmonic space of degree \(a+b+2m\). The multiplier formula in
		Lemma~\ref{lem:riesz-positive} gives
		\[
		J_\mu(g_{m,a,b}Y_aZ_b)
		=\kappa_{a+b+2m}g_{m,a,b}Y_aZ_b.
		\]
		This proves \eqref{eq:angular-riesz-multipliers}. The radial Riesz
		operator is therefore the same for every unit angular vector in
		the sector.
	\end{proof}

	For a fixed sector, put \(g_m=g_{m,a,b}\), and use
	\(\langle\cdot,\cdot\rangle_\rho\) for the scalar product in
	\(L^2((0,\pi/2),\rho\,dt)\). The operator associated with the
	radial form is
	\[
	L_{a,b}=A_{a,b}-(p-1)V_v-p\,a_vJ_{\mu,a,b}(a_v\cdot).
	\]
	Its polarized form has Jacobi matrix
	\begin{equation}
		\begin{aligned}
			Q_{a,b}[g_r,g_s]
			={}&\lambda_{r,a,b}\delta_{rs}
			-(p-1)\langle V_vg_r,g_s\rangle_\rho\\
			&-p\sum_{m\geq0}\kappa_{a+b+2m}
			\langle a_vg_r,g_m\rangle_\rho
			\langle a_vg_s,g_m\rangle_\rho .
		\end{aligned}
		\label{eq:jacobi-hessian-matrix}
	\end{equation}
	Indeed, expand each of \(a_vg_r\) and \(a_vg_s\) in the orthonormal
	Jacobi basis and use \eqref{eq:angular-riesz-multipliers}. Since
	\(a_v\in L^\infty\) and \(0<\kappa_m\leq\kappa_0\),
	\[
	\begin{aligned}
		&\sum_{m\geq0}\kappa_{a+b+2m}
		|\langle a_vg_r,g_m\rangle_\rho|
		|\langle a_vg_s,g_m\rangle_\rho|\\
		&\qquad\leq\kappa_0\|a_vg_r\|_{L^2(\rho\,dt)}
		\|a_vg_s\|_{L^2(\rho\,dt)}<\infty.
	\end{aligned}
	\]
	This justifies the series and passage from finite expansions to the
	closed form. The Riesz operator is diagonal in this basis, whereas
	the multiplications by \(a_v\) and \(V_v\) generally are not.
	
	\begin{proposition}
		{
			\label{prop:transverse-localization}
			Let \(p\geq2\), and let \(v\neq0\) be a \(G\)-invariant critical point of \(E\) on \(H^1(\mathbb S^n)\).
			{
				Set
				\begin{equation}
					M_v
					=
					(2p-1)\|V_v\|_{L^\infty(\mathbb S^n)},
					\label{eq:M-v}
				\end{equation}
				where \(V_v\) is defined in \eqref{eq:linearized-coefficients}.}
			Then
			\begin{equation}
				Q_v(\phi)
				\geq
				\|\phi\|^2
				-
				M_v\|\phi\|_{L^2(\mathbb S^n)}^2
				\qquad
				\text{for every }\phi\in H^1(\mathbb S^n).
				\label{eq:hessian-spectral-lower-bound}
			\end{equation}
			For every bidegree \((a,b)\),
			\begin{equation}
				m_{a,b}(v)
				\leq
				\#\bigl\{
				m\geq0:\lambda_{m,a,b}<M_v
				\bigr\}.
				\label{eq:negative-block-count}
			\end{equation}
			The nonpositive spectral space satisfies
			\[
			m_{a,b}(v)+\nu_{a,b}(v)
			\leq\#\{m\geq0:\lambda_{m,a,b}\leq M_v\}.
			\]
			In particular,
			\begin{equation}
				m_{a,b}(v)=0
				\quad\text{if}\quad
				(a+b)(a+b+n-1)
				+
				\frac{n(n-2)}{4}
				\geq M_v.
				\label{eq:negative-block-cutoff}
			\end{equation}
			If the inequality in \eqref{eq:negative-block-cutoff} is strict, then
			\(\nu_{a,b}(v)=0\). If \(\mathcal N_A(s)\) denotes the number of
			eigenvalues of \(A\) below \(s\), counted with multiplicity, then
			\begin{equation}
				m(v)
				\leq
				\mathcal N_A(M_v)
				\leq
				C_n\bigl(1+M_v^{n/2}\bigr).
				\label{eq:static-morse-count}
			\end{equation}
		}
	\end{proposition}
	
	\begin{proof}
		The pointwise kernel is positive, and Lemma~\ref{lem:riesz-positive}
		gives
		\(\int J_\mu(a_v\phi)a_v\phi\,dV\geq0\). For the upper bound, apply
		\(2|xy|\leq x^2+y^2\) under the double integral with
		\[
		x=|v(\xi)|^{(p-2)/2}|v(\eta)|^{p/2}\phi(\xi),
		\qquad
		y=|v(\eta)|^{(p-2)/2}|v(\xi)|^{p/2}\phi(\eta).
		\]
		The two square terms have the same integral after interchanging
		\(\xi\) and \(\eta\). Hence
		\[
		0\leq\int_{\mathbb S^n}J_\mu(a_v\phi)a_v\phi\,dV
		\leq\int_{\mathbb S^n}J_\mu(|v|^p)|v|^{p-2}\phi^2\,dV
		=\int_{\mathbb S^n}V_v\phi^2\,dV.
		\]
		Substituting into \eqref{eq:hessian-form} yields
		\[
		\begin{aligned}
			Q_v(\phi)
			&\geq\|\phi\|^2-\bigl((p-1)+p\bigr)
			\int_{\mathbb S^n}V_v\phi^2\,dV\\
			&\geq\|\phi\|^2-M_v\|\phi\|_2^2.
		\end{aligned}
		\]
		Moreover,
		\[
		M_v\leq(2p-1)\kappa_0\|v\|_\infty^{2p-2}.
		\]
		This proves \eqref{eq:hessian-spectral-lower-bound}, with the sharper
		threshold \(M_v\) of \eqref{eq:M-v}.
		
		{
			For $g\in\mathcal X_{a,b}$, write
			\[
			\|g\|_{A_{a,b}}^2=
			\int_0^{\pi/2}\left[|g'|^2+
			\left(\frac{\lambda_{k,a}}{\cos^2t}
			+\frac{\lambda_{\ell,b}}{\sin^2t}+\frac{n(n-2)}4\right)|g|^2\right]
			\rho(t)\,dt.
			\]
			Fix a bidegree $(a,b)$, and define
			\[
			\mathcal Y_{a,b}
			=
			\overline{
				\operatorname{span}
				\bigl\{
				g_{m,a,b}:
				\lambda_{m,a,b}\geq M_v
				\bigr\}
			}^{\,\mathcal X_{a,b}}.
			\]
			Proposition~\ref{prop:angular-jacobi-spectrum} gives
			\[
			\operatorname{codim}_{\mathcal X_{a,b}}
			\mathcal Y_{a,b}
			=
			\#\bigl\{
			m\geq0:
			\lambda_{m,a,b}<M_v
			\bigr\}.
			\]
			For
			\[
			g
			=
			\sum_{\lambda_{m,a,b}\geq M_v}
			c_mg_{m,a,b}
			\in\mathcal Y_{a,b},
			\]
			the spectral decomposition of the local radial operator gives
			\[
			\begin{aligned}
				\|g\|_{A_{a,b}}^2
				&=
				\sum_{\lambda_{m,a,b}\geq M_v}
				\lambda_{m,a,b}|c_m|^2\\
				&\geq
				M_v
				\sum_{\lambda_{m,a,b}\geq M_v}
				|c_m|^2
				=
				M_v\|g\|_{L^2((0,\pi/2),\rho\,dt)}^2.
			\end{aligned}
			\]
			The continuity of \(Q_{a,b}\) and
			\eqref{eq:hessian-spectral-lower-bound} therefore imply
			\[
			Q_{a,b}(g)\geq0
			\qquad
			\text{for every }g\in\mathcal Y_{a,b}.
			\]
		}
		Let \(\mathcal V\subset\mathcal X_{a,b}\) be a subspace on which
		\(Q_{a,b}\) is negative definite. If
		\(\dim\mathcal V>\operatorname{codim}\mathcal Y_{a,b}\), then the
		{finite-codimension} intersection principle gives
		\[
		\mathcal V\cap\mathcal Y_{a,b}\neq\{0\}.
		\]
		This contradicts the nonnegativity of \(Q_{a,b}\) on
		\(\mathcal Y_{a,b}\). Hence,
		\[
		m_{a,b}(v)
		\leq
		\operatorname{codim}\mathcal Y_{a,b},
		\]
		which proves \eqref{eq:negative-block-count}.
		For the bound on the nonpositive spectral space, replace
		\(\lambda_{m,a,b}\geq M_v\) by \(\lambda_{m,a,b}>M_v\).
		On the resulting closed subspace every nonzero vector satisfies
		\[
		Q_{a,b}(g)\geq
		\sum_{\lambda_{m,a,b}>M_v}(\lambda_{m,a,b}-M_v)|c_m|^2>0.
		\]
		It has codimension \(\#\{m:\lambda_{m,a,b}\leq M_v\}\).
		A larger nonpositive spectral subspace would intersect this positive
		subspace nontrivially, a contradiction.

		{
			Proposition~\ref{prop:angular-jacobi-spectrum} gives
			\[
			\lambda_{0,a,b}
			=
			(a+b)(a+b+n-1)
			+
			\frac{n(n-2)}{4}.
			\]
			If \(\lambda_{0,a,b}\geq M_v\), then
			\(\mathcal Y_{a,b}=\mathcal X_{a,b}\), and hence
			\(Q_{a,b}\geq0\) on \(\mathcal X_{a,b}\). This proves
			\eqref{eq:negative-block-cutoff}. If
			\(\lambda_{0,a,b}>M_v\), then
			\[
			Q_{a,b}(g)
			\geq
			\bigl(\lambda_{0,a,b}-M_v\bigr)
			\|g\|_{L^2((0,\pi/2),\rho\,dt)}^2
			>0
			\]
			for every \(g\neq0\). Therefore, \(\nu_{a,b}(v)=0\).
		}
		
		Let \(\mathcal Y\) be the closure in \(H^1(\mathbb S^n)\) of the span
		of all eigenfunctions of \(A\) whose eigenvalues are at least \(M_v\).
		The spectral theorem gives
		\[
		\operatorname{codim}_{H^1(\mathbb S^n)}\mathcal Y
		=
		\mathcal N_A(M_v).
		\]
		The spectral decomposition of \(A\), together with
		\eqref{eq:hessian-spectral-lower-bound}, yields
		\[
		Q_v(\phi)\geq0
		\qquad
		\text{for every }\phi\in\mathcal Y.
		\]
		The same {finite-codimension} intersection argument then gives
		\[
		m(v)\leq\mathcal N_A(M_v).
		\]
		
		{
			To distinguish the full harmonic degree from the invariant Jacobi
			label, write the eigenvalues of $A$ as
			\[
			\chi_L
			=
			L(L+n-1)
			+
			\frac{n(n-2)}{4},
			\qquad L\geq0,
			\]
			with multiplicity
			\[
			\dim\mathcal H_L(\mathbb S^n)
			=
			\binom{n+L}{n}
			-
			\binom{n+L-2}{n}
			=
			O\bigl((L+1)^{n-1}\bigr).
			\]
			Here \(\binom{q}{n}=0\) for nonnegative integers
			\(q<n\). The inequality \(\chi_L<M_v\) implies
			\(L\leq C_n(1+M_v^{1/2})\). Therefore,
			\[
			\begin{aligned}
				\mathcal N_A(M_v)
				&\leq
				C_n
				\sum_{0\leq L\leq C_n(1+M_v^{1/2})}
				(L+1)^{n-1}\\
				&\leq
				C_n\bigl(1+M_v^{n/2}\bigr),
			\end{aligned}
			\]
			which completes the proof of
			\eqref{eq:static-morse-count}.
		}
	\end{proof}
	
	\subsection{Transverse negative directions}

	\begin{proposition}
		\label{prop:constant-spectrum}
		Let \(p\geq2\). Each constant solution \(v=a_0\) or \(v=-a_0\)
		satisfies
		\[
		m(v)=m_G(v)=1,
		\qquad
		\nu(v)=n+1,
		\qquad
		\nu_G(v)=0.
		\]
		The kernel of the full Hessian is precisely the space of {degree-one}
		spherical harmonics. Moreover, the first minimax level \(c_1\) is
		isolated in the set of critical values of \(\Psi\), and the critical
		set of \(\Psi\) at this level consists of one antipodal pair.
	\end{proposition}
	
	\begin{proof}
		{
			{The constant state identity}
			\(
			{\kappa_0a_0^{2p-2}={n(n-2)}/4}
			\)
			and the {Funk--Hecke} multiplier formula
			\eqref{eq:riesz-eigenvalues} show that the \(L^2(\mathbb S^n)\)-realization \(L_v\)
			is diagonal with respect to the spherical harmonic decomposition.
			On the space \(\mathcal H_L(\mathbb S^n)\) of degree-\(L\) spherical
			harmonics, the corresponding eigenvalue is
			\[
			q_L
			=
			L(L+n-1)
			+
			{(2-p)\frac{n(n-2)}4}
			-
			{p\frac{n(n-2)}4\frac{\kappa_L}{\kappa_0}}.
			\]
			To determine the Morse index and nullity, we transfer this spectral
			information to the energy space operator $T_v$.  The operator $T_v$
			has eigenvalue
			\[
			{\frac{q_L}{L(L+n-1)+n(n-2)/4}}
			\]
			on the same harmonic space. The denominator is strictly positive, so
			the two realizations have the same signs and the same zero
			eigenspaces.
		}
		
		Formula \eqref{eq:riesz-eigenvalues} gives
		\[
		\frac{\kappa_1}{\kappa_0}
		=
		\frac{\mu}{2n-\mu}.
		\]
		Using \(p=(2n-\mu)/(n-2)\), direct substitution yields
		\[
		{q_0=-2(p-1)\frac{n(n-2)}4<0},
		\qquad
		q_1=0.
		\]
		For every \(L\geq2\),
		\[
		\begin{aligned}
			q_L-q_1
			&=
			L(L+n-1)-n
			+
			{p\frac{n(n-2)}4\frac{\kappa_1-\kappa_L}{\kappa_0}}.
		\end{aligned}
		\]
		The Laplace eigenvalues \(L(L+n-1)\) increase strictly with \(L\), while
		the multipliers \(\kappa_L\) decrease strictly by
		\eqref{eq:riesz-eigenvalues}. Hence,
		\[
		q_L-q_1>0
		\qquad
		\text{for every }L\geq2.
		\]
		Thus, the {degree-zero} harmonic space is the only negative eigenspace,
		and the {degree-one} harmonic space is the entire kernel. Since
		\(\dim\mathcal H_0(\mathbb S^n)=1\) and
		\(\dim\mathcal H_1(\mathbb S^n)=n+1\), it follows that
		\[
		m(v)=1,
		\qquad
		\nu(v)=n+1.
		\]
		
		Proposition~\ref{prop:jacobi-spectrum} shows that the
		\(G\)-invariant spherical harmonics have even degrees \(L=2m\). The
		negative {degree-zero} mode is \(G\)-invariant, whereas the {degree-one}
		kernel contains no nonzero \(G\)-invariant function. All invariant
		modes of positive even degree have positive Hessian eigenvalues.
		Therefore,
		\[
		m_G(v)=1,
		\qquad
		\nu_G(v)=0.
		\]
		
		{
			Let
			\[
			w_0=\frac{a_0}{\|a_0\|}\in\mathcal S.
			\]
			Lemma~\ref{lem:positive-invariant}, together with the radial
			correspondence between critical points of \(\Psi\) and nonzero
			critical points of \(E|_H\), gives
			\[
			\{w\in\mathcal S:
			d(\Psi|_{\mathcal S})(w)=0,\ \Psi(w)=c_1\}
			=
			\{w_0,-w_0\}.
			\]
			By Lemma~\ref{lem:radial-index-shift},
			\[
			\widetilde\nu(w_0)
			=
			\nu_G(a_0)
			=
			0,
			\qquad
			\widetilde\nu(-w_0)
			=
			\nu_G(-a_0)
			=
			0.
			\]
			Thus, \(w_0\) and \(-w_0\) are nondegenerate critical points of
			{$\Psi$. The derivative of its gradient in a local Hilbert
				chart is an invertible Fredholm Hessian at each of the two points.
				The inverse function theorem therefore makes each critical point
				isolated.}
		}
		
		Suppose that \(c_1\) were not isolated in the set of critical values.
		Then there would exist distinct critical values \(d_j\to c_1\) and
		critical points \(w_j\in\mathcal S\) satisfying
		\[
		\Psi(w_j)=d_j,
		\qquad
		d(\Psi|_{\mathcal S})(w_j)=0.
		\]
		The sequence \(\{w_j\}\) is a {Palais--Smale} sequence at level \(c_1\).
		Lemma~\ref{lem:psi-ps} provides a subsequence converging strongly in
		\(H\) to a critical point \(w\) of \(\Psi\) with \(\Psi(w)=c_1\).
		Hence, \(w=w_0\) or \(w=-w_0\). This contradicts the local isolation
		of the two critical points. Therefore, \(c_1\) is an isolated critical
		value.
	\end{proof}
	
	{
		\begin{proof}[Proof of Theorem~\ref{thm:conformal-kernel-intro}]
			Let $u=Pv$. By Proposition~\ref{prop:conformal-identities}, $u$ is a
			Euclidean weak solution. {Lemma~\ref{lem:regularity} and the
				formula defining $P$ yield}
			\(
			u\in C_{\mathrm{loc}}^{2,\gamma}(\mathbb R^n)
			\)
			for some $\gamma\in(0,1)$, which justifies differentiating the
			translation and dilation families introduced below.
			For \(\lambda>0\) and \(a\in\mathbb R^n\), define
			\begin{equation}
				u_{\lambda,a}(x)
				=
				\lambda^{(n-2)/2}u\bigl(\lambda(x-a)\bigr).
				\label{eq:translation-dilation-family}
			\end{equation}
			Set \(X=\lambda(x-a)\) and \(z=\lambda(y-a)\). Then
			\[
			-\Delta u_{\lambda,a}(x)
			=
			\lambda^{(n+2)/2}(-\Delta u)(X).
			\]
			Moreover,
			\[
			\begin{aligned}
				\int_{\mathbb R^n}
				\frac{|u_{\lambda,a}(y)|^p}{|x-y|^\mu}\,dy
				&=
				\lambda^{(n-2)p/2+\mu-n}
				\int_{\mathbb R^n}
				\frac{|u(z)|^p}{|X-z|^\mu}\,dz\\
				&=
				\lambda^{\mu/2}
				\int_{\mathbb R^n}
				\frac{|u(z)|^p}{|X-z|^\mu}\,dz,
			\end{aligned}
			\]
			where \(p=(2n-\mu)/(n-2)\) was used in the last equality. The same
			identity gives
			\[
			|u_{\lambda,a}(x)|^{p-2}u_{\lambda,a}(x)
			=
			\lambda^{(n+2-\mu)/2}
			|u(X)|^{p-2}u(X).
			\]
			Hence, both sides of the equation are multiplied by
			\(\lambda^{(n+2)/2}\), and \(u_{\lambda,a}\) is a weak solution for every
			\(\lambda>0\) and \(a\in\mathbb R^n\).
			
			The family \eqref{eq:translation-dilation-family} is induced by smooth
			conformal transformations of \(\mathbb S^n\). The
			\(C^{2,\gamma}\)-regularity of \(v\) gives differentiability of the
			corresponding parameter action in \(H^1(\mathbb S^n)\), and
			Proposition~\ref{prop:conformal-identities} transfers this
			differentiability to \(D^{1,2}(\mathbb R^n)\).
			Lemma~\ref{lem:twice-differentiability} justifies differentiation of the
			weak equation with respect to the parameters. At
			\((\lambda,a)=(1,0)\),
			\[
			\left.
			\partial_\lambda u_{\lambda,0}
			\right|_{\lambda=1}
			=
			\frac{n-2}{2}u+x\cdot\nabla u
			=
			Z_0,
			\qquad
			\left.
			\partial_{a_i}u_{1,a}
			\right|_{a=0}
			=
			-\partial_{x_i}u
			=
			-Z_i.
			\]
			Consequently, \(Z_0,\ldots,Z_n\in D^{1,2}(\mathbb R^n)\), and
			differentiation of the weak equation gives
			\[
			\mathcal L_uZ_i=0,
			\qquad
			0\leq i\leq n.
			\]

			To prove linear independence, suppose that
			\begin{equation}
				b_0Z_0+b\cdot\nabla u=0,
				\qquad
				b=(b_1,\ldots,b_n)\in\mathbb R^n.
				\label{eq:linear-relation-conformal}
			\end{equation}
			Assume first that \(b_0=0\). If \(b\neq0\), an orthogonal change of
			coordinates transforms \(b\cdot\nabla\) into
			\(|b|\partial_{x_1}\). Relation
			\eqref{eq:linear-relation-conformal} then implies that \(u\) is
			independent of \(x_1\). The Sobolev inequality gives
			\(u\in L^{2^*}(\mathbb R^n)\), whereas Fubini's theorem shows that a
			nonzero \(x_1\)-independent function cannot belong to
			\(L^{2^*}(\mathbb R^n)\). Thus \(u=0\), contradicting the assumption
			that \(u\) is nonzero. Therefore, \(b=0\) whenever \(b_0=0\).
			
			Assume next that \(b_0\neq0\), and set \(x_0=-b/b_0\). Substitution of
			the definition of \(Z_0\) into
			\eqref{eq:linear-relation-conformal} gives
			\begin{equation}
				(x-x_0)\cdot\nabla u
				+
				\frac{n-2}{2}u
				=
				0.
				\label{eq:homogeneity-relation}
			\end{equation}
			Along almost every ray \(x=x_0+r\theta\),
			\eqref{eq:homogeneity-relation} becomes
			\[
			r\frac{d}{dr}u(x_0+r\theta)
			+
			\frac{n-2}{2}u(x_0+r\theta)
			=
			0.
			\]
			Integration yields
			\[
			u(x_0+r\theta)
			=
			r^{-(n-2)/2}A(\theta).
			\]
			If \(A\) is nonzero on a subset of \(\mathbb S^{n-1}\) of positive
			measure, then the identity \((n-2)2^*/2=n\) gives
			\[
			\int_{\mathbb R^n}|u|^{2^*}\,dx
			\geq
			C\int_0^1
			r^{n-1-(n-2)2^*/2}\,dr
			=
			C\int_0^1\frac{dr}{r}
			=
			\infty.
			\]
			This contradicts \(u\in L^{2^*}(\mathbb R^n)\). Hence,
			\(A=0\) almost everywhere and \(u=0\), which again contradicts the
			nontriviality of \(u\). Therefore, \(b_0=0\), and the preceding case
			then gives \(b=0\). Thus \(Z_0,\ldots,Z_n\) are linearly independent.
			
			Proposition~\ref{prop:conformal-identities} shows that
			\(P\) is a linear isometry preserving the full energy, including
			the nonlocal term. Differentiating this full energy identity twice
			identifies the Euclidean and spherical Hessians.
			Thus \(P^{-1}Z_0,\ldots,P^{-1}Z_n\) are independent elements of the
			spherical kernel, and \(\nu(v)\geq n+1\).
			
			For distinct coordinate labels \(i,j\), write
			\(R_{ij}=\xi_i\partial_{\xi_j}-\xi_j\partial_{\xi_i}\).
			The Euclidean and spherical generators are related as follows.
			In the stereographic chart of
			\eqref{eq:conformal-map}, write
			\(\xi_i=\Omega x_i\) for \(i\leq n\) and
			\(\xi_{n+1}=(1-|x|^2)/(1+|x|^2)\). Direct differentiation gives
			\[
			D\pi^{-1}(x)\partial_{x_i}=\nabla\xi_i-R_{i,n+1},
			\qquad
			D\pi^{-1}(x)x=-\nabla\xi_{n+1},
			\]
			and
			\[
			\partial_{x_i}\log\Omega=-\xi_i,
			\qquad 1+x\cdot\nabla\log\Omega=\xi_{n+1}.
			\]
			Applying the product rule to \(u(x)=\Omega(x)^{(n-2)/2}v(\pi^{-1}(x))\)
			therefore gives
			\begin{equation}
				P^{-1}Z_i=\Phi_{e_i}-R_{i,n+1}v\quad(1\leq i\leq n),
				\qquad P^{-1}Z_0=-\Phi_{e_{n+1}}.
				\label{eq:chart-kernel-fields}
			\end{equation}
			Here \(e_i\) denotes the \(i\)-th coordinate vector.
			For an invariant nonconstant \(v\), the terms \(R_{i,n+1}v\)
			with \(i\leq k\) belong to the sector \((1,1)\). Thus individual
			pulled-back Euclidean translations need not have a pure degree-one
			angular profile.
			
		    If \(v\in H\) is nonconstant, Proposition~\ref{prop:conformal-blocks}
			gives the dimension of the conformal orbit. The two gradient families
			lie in the \((1,0)\) and \((0,1)\) sectors, with total dimension
			\(
			k+\ell=n+1.
			\)
			The mixed rotations occupy the \((1,1)\) sector and have dimension
			\(k\ell\). These three angular sectors are mutually orthogonal.
			Consequently,
			\[
			\nu(v)\geq n+1+k\ell,
			\qquad
			\dim\mathcal K(v)=n+1+k\ell.
			\]

			For the same nonconstant \(G\)-invariant solution \(v\), it remains to prove
			\eqref{eq:intro-transverse-index}. Set \(F=|v|^p\) and use the test
			functions \(\phi_i=\xi_i v\)
		for \(1\leq i\leq n+1\). Testing the equation against \(\xi_i^2 v\)
		gives
		\[
		\|\xi_i v\|^2=\int_{\mathbb S^n}V_v \xi_i^2 v^2\,dV
		+\int_{\mathbb S^n}v^2|\nabla\xi_i|^2\,dV .
		\]
		Since \(\displaystyle \sum_i\xi_i^2=1\), \(\displaystyle \sum_i|\nabla\xi_i|^2=n\), and
		\(\xi\cdot\eta=1-\frac12|\xi-\eta|^2\), summing
		\eqref{eq:hessian-form} yields
		\[
		\sum_i Q_v(\xi_i v)=
		n\|v\|_2^2-(2p-2)\mathcal D(v)
		+\frac p2\iint_{\mathbb S^n\times\mathbb S^n}
		|\xi-\eta|^{2-\mu}F(\xi)F(\eta)\,dV_\xi dV_\eta .
		\]
		{Applying \cite[Proposition~4.1]{FrankLieb2012}
		in dimension \(n\) and with \(\alpha=\mu/2\) gives}
		\[
		\iint_{\mathbb S^n\times\mathbb S^n}
		\frac{F(\xi)(\xi\cdot\eta)F(\eta)}{|\xi-\eta|^\mu}
		\,dV_\xi dV_\eta
		\geq\frac{\mu}{2n-\mu}\mathcal D(v).
		\]
		Since \(|\xi-\eta|^2=2(1-\xi\cdot\eta)\),
		\[
		\begin{aligned}
			\iint_{\mathbb S^n\times\mathbb S^n}
			|\xi-\eta|^{2-\mu}F(\xi)F(\eta)\,dV_\xi dV_\eta
			&\leq2\left(1-\frac{\mu}{2n-\mu}\right)\mathcal D(v)=\frac{4(n-\mu)}{2n-\mu}\mathcal D(v).
		\end{aligned}
		\]
		Using \(p=(2n-\mu)/(n-2)\) and \eqref{eq:critical-testing} gives
		\begin{equation}
			\sum_i Q_v(\xi_i v)\leq
			n\|v\|_2^2-\frac4{n-2}\|v\|^2
			=-\frac4{n-2}\int_{\mathbb S^n}|\nabla v|^2\,dV.
			\label{eq:transverse-coordinate-test}
		\end{equation}
			The right-hand side of \eqref{eq:transverse-coordinate-test} is
		negative for a nonconstant solution.
		The test functions belong to \((1,0)\) or \((0,1)\), so at least
		one of these radial blocks has a negative direction.
		{Using \(m_{0,0}(v)=m_G(v)\) and the angular
			multiplicities in \eqref{eq:full-index-sum} gives
			\[
			m(v)\geq m_G(v)+k\,m_{1,0}(v)+\ell\,m_{0,1}(v)
			\geq m_G(v)+\min\{k,\ell\}.
			\]
			This proves \eqref{eq:intro-transverse-index}.}
			When \(k=\ell\) and \(v(\mathcal T\xi)=\pm v(\xi)\), the factor
			interchange identifies the two blocks. Both then contain a negative
			direction.

			For \(v=\pm a_0\), Proposition~\ref{prop:constant-spectrum} shows that the full kernel is exactly the \((n+1)\)-dimensional conformal tangent space.
			This completes the proof of Theorem~\ref{thm:conformal-kernel-intro}.
		\end{proof}
	}
	
	\begin{remark}
		\label{rem:conditional-index-asymptotics}
		{
			Assume $p\geq2$, and let $j\to\infty$ through labels for which the
			isolation hypotheses of Theorem~\ref{thm:morse-intro} hold. For each
			such label, choose the solution $v_j$ supplied by that theorem. Then
			\[
			0
			\leq
			j-m_G(v_j)
			\leq
			\nu_G(v_j),
			\qquad
			N_G(v_j)
			\leq
			m_G(v_j)
			\leq
			m(v_j).
			\]
			If $\nu_G(v_j)=o(j)$ along these labels, then
			$m_G(v_j)/j\to1$. Consequently,
			\[
			\frac{N_G(v_j)}{j}\longrightarrow1
			\quad\Longleftrightarrow\quad
			m_G(v_j)-N_G(v_j)=o(j),
			\]
			and
			\[
			\frac{m(v_j)}{j}\longrightarrow1
			\quad\Longleftrightarrow\quad
			m(v_j)-m_G(v_j)=o(j).
			\]
			By \eqref{eq:full-index-sum},
			\[
			m(v_j)-m_G(v_j)
			=
			\sum_{\substack{a,b\geq0\\(a,b)\neq(0,0)}}
			d_{k,a}d_{\ell,b}m_{a,b}(v_j),
			\]
			which is the total contribution of the noninvariant angular sectors.
			Moreover, Proposition~\ref{prop:transverse-localization} gives
			\[
			m(v_j)
			\leq
			C\left(
			1+
			\|v_j\|_{L^\infty(\mathbb S^n)}^{n(p-1)}
			\right).
			\]}
	\end{remark}
	
	\subsection{The low blocks for equal factors}
	
	\begin{proof}[Proof of Theorem~\ref{thm:swap-blocks}]
		Let \(k=\ell\), so \(n=2k-1\). Write \(p_\mu\), \(E_\mu\), and
		\(\mathcal D_\mu\) for the parameter-dependent quantities. At
		\(\mu=0\), set \(p_0=2^*\) and
		\(J_0 h=\int_{\mathbb S^n}h\,dV\).
		
		\medskip
		\noindent\textbf{Minimizers and their limit.}
		On \(\mathcal S\cap H^-\), minimizing \(\Psi\) is equivalent to
		maximizing \(\mathcal D_\mu\). A normalized multiple of
		\(|\xi'|^2-|\xi''|^2\) belongs to this sphere. A maximizing
		sequence has a subsequence converging weakly in \(H^-\) and
		strongly in \(L^{2^*}\).
		The interaction converges by HLS for \(\mu>0\) and directly
		for \(\mu=0\). The limiting mass is positive, so the limit is
		nonzero. If its norm were less than one, normalization would
		strictly increase the mass. Hence it is a maximizer on
		\(\mathcal S\cap H^-\).
		
		Rescale a maximizer by \eqref{eq:rescaling-critical}. The proof of
		Lemma~\ref{lem:radial-reduction} applies on \(H^-\).
		The energy is invariant under \(v(\xi)\mapsto -v(\mathcal T\xi)\), which fixes
		the rescaled function. Averaging a variation under this involution
		and \(G\) reduces it to a variation in \(H^-\). Thus the rescaled
		function is a full solution. In the orbit variable, \(\mathcal T\) acts
		as \(t\mapsto\pi/2-t\). Put \(t_*=\pi/4\).
		Every profile in \(H^-\) has zero trace at \(t_*\).
		Replacing its profile by the positive absolute value on
		\((0,t_*)\) and the negative absolute value on
		\((t_*,\pi/2)\) preserves its norm and its mass.
		The matching zero traces give an \(H^1\) function.
		Consequently, the resulting maximizer also yields a solution.
		
		For \(\mu>0\) small enough that \(p_\mu>2\), the folded
		profile is nonnegative on the first half. An interior zero there
		would be simple by Lemma~\ref{lem:simple-zeros} and would force a
		sign change, a contradiction. Endpoint zeros are excluded by the
		same lemma. Hence the folded profile is strictly positive on
		\([0,t_*)\) and has exactly one simple zero at \(t_*\).
		Any original maximizer has the same absolute value as its folded
		profile and is continuous, so it has this sign pattern up to an
		overall sign. In the limit problem,
		any rescaled maximizer \(v^{(0)}\) solves
		\[
		A v^{(0)}
		=B_0 |v^{(0)}|^{p_0-2}v^{(0)},\qquad
		B_0=\int_{\mathbb S^n}|v^{(0)}|^{p_0}\,dV>0.
		\]
		By Lemma~\ref{lem:compact-embedding}, \(v^{(0)}\in L^R\)
		for some \(R>2^*\). Its potential
		\(B_0 |v^{(0)}|^{p_0-2}\) is then in \(L^{r_0}\) with \(r_0>n/2\).
		Elliptic \(L^s\) estimates and Sobolev embedding give boundedness.
		In each subcritical step the reciprocal integrability exponent of
		\(v^{(0)}\) decreases by \(2/n-1/r_0>0\). At a critical step one
		may take any finite exponent before applying the estimate again.
		Schauder regularity gives \(C^{2,\gamma}\).
		The profile equation and Lemma~\ref{lem:simple-zeros} apply with the
		constant potential \(B_0\), so the same sign pattern holds at \(\mu=0\).
		
		Fix \(R>p_0\) for which the embedding \(H\hookrightarrow L^R\)
		is compact. Choose \(r>1\) close enough to one that
		\[
		\frac1r>\max\left\{\frac{p_0}{R},
		\frac{p_0-1}{R}+\frac1{2^*}\right\}.
		\]
		Choose
		\(0<\overline\mu<\min\{4,n\}\) so small that
		\(\overline\mu\,r'<n\), where
		\(r'=r/(r-1)\). For \(0\leq\mu\leq\overline\mu\), the operators \(J_\mu\) are uniformly bounded
		from \(L^r\) to \(L^\infty\). In fact,
		\[
		\|J_\mu-J_0\|_{L^r\to L^\infty}
		\leq\left(\int_{\mathbb S^n}
		\big||\xi-\eta|^{-\mu}-1\big|^{r'}\,dV_\eta\right)^{1/r'}
		\longrightarrow0
		\]
		by rotation invariance and dominated convergence.
		
		Let \(\mu_i\to0\) and let \(w_i\) be corresponding
		maximizers with the positive sign on the first half.
		Pass to a subsequence with \(w_i\rightharpoonup w\) in \(H\)
		and \(w_i\to w\) in \(L^R\). We will use the following fact. If \(z_i\to z\) in \(L^R\), \(s_i\to s>0\) with
		\(s_i\leq s\), and \(sq<R\), then \(|z_i|^{s_i}\to |z|^s\)
		in \(L^q\), also with a sign factor. Indeed, there is convergence
		in measure and the powers are bounded in \(L^{q+\varepsilon}\)
		for small \(\varepsilon>0\). In particular,
		\(|w_i|^{p_{\mu_i}}\to |w|^{p_0}\) in \(L^r\) and
		\(\mathcal D_{\mu_i}(w_i)\to\mathcal D_0(w)\).
		If \(M_0=\max_{\mathcal S\cap H^-}\mathcal D_0\), comparison
		with a fixed maximizing profile at \(\mu=0\) gives
		\(\mathcal D_0(w)\geq M_0\). On the other hand,
		\(\mathcal D_0(w)\leq M_0\|w\|^{2p_0}\leq M_0\).
		Thus \(w_i\to w\) in \(H\), and \(w\) is a maximizer at zero.
		The strong \(L^R\) convergence preserves nonnegativity
		on the first half, and \(w\in H^-\). The regularity and zero
		argument for the limit equation therefore give the same strict
		sign pattern for its rescaling \(v^{(0)}\).
		The rescaled solutions converge strongly in \(H\) and in \(L^R\)
		to \(v^{(0)}=\mathcal D_0(w)^{-1/(2p_0-2)}w\).
		
		\medskip
		\noindent\textbf{The three limiting radial operators.}
		The nonlocal term of the Hessian at \(\mu=0\) vanishes on
		every sector \((a,b)\ne(0,0)\). The angular average of
		\(a_{v^{(0)}}gY_aZ_b\) is zero, and \(J_0\) is integration.
		Thus these radial forms are local, with the common potential
		\((p_0-1)B_0|f|^{p_0-2}\), where \(f\) is the profile of
		\(v^{(0)}\). The symmetry profiles from
		Proposition~\ref{prop:conformal-blocks} still solve the kernel
		equations at \(\mu=0\), since its conformal computation also holds
		for \(p=p_0\) and \(J_0\).
		
		{Consider}
		\[
		\mathcal E_f(t)=\frac12 f'(t)^2+\frac{B_0}{p_0}f(t)^{p_0}
		-{\frac{n(n-2)}8} f(t)^2,\qquad 0\leq t\leq t_* .
		\]
		On the open first half,
		\(\mathcal E_f'=-(\rho'/\rho)f'^2\leq0\), since \(k=\ell\).
		At \(t_*\), \(\mathcal E_f(t_*)=f'(t_*)^2/2 >0\).
		At any stationary point on \([0,t_*)\) we therefore have
		{\[B_0 f^{p_0-2}>\frac{p_0n(n-2)}{8}>\frac{n(n-2)}{4}.\]}
		The profile equation implies \(f''<0\) at every interior stationary
		point. At zero, \(f'(0)=0\) and taking the limit in
		\eqref{eq:profile-equation} gives
		\[
		{k f''(0)=\frac{n(n-2)}4 f(0)-B_0 f(0)^{p_0-1}<0} .
		\]
		Thus \(f'<0\) near zero. A first interior zero of \(f'\) would
		have \(f''\geq0\), so \(f'<0\) on the first half. Oddness about
		\(t_*\) and the simple zero at \(t_*\) give \(f'<0\) on the
		whole open interval.
		
		In the mixed sector, \(g=-f'>0\) is a kernel function of the
		local radial operator. For a compactly supported test profile
		\(\phi\), integration by parts in the equation for \(g\) gives
		\[
		Q_{1,1}^{(0)}(\phi)=\int_0^{\pi/2}\rho\, g^2
		\big|(\phi/g)'\big|^2\,dt .
		\]
		Density in the form norm gives nonnegativity.
		If \(\phi\) is in the kernel, approximate it by compactly supported
		profiles in that norm. On every compact subinterval the quotient by
		\(g\) converges in \(H^1\), so the last identity implies
		\((\phi/g)'=0\) there. Hence the mixed block has index zero and
		nullity one.
		
		In the \((1,0)\) sector, the profile \(g_{1,0}\) in
		\eqref{eq:conformal-block-profiles} has exactly one interior zero.
		On \(0<t<t_*\), set \(h=(n-2)/2\) and
		\(y=(\sin t)^h f\). Since
		\(h=k-3/2\), {\(n(n-2)/4=h(h+1)\)}, and
		\(\rho'/\rho=(k-1)(\cot t-\tan t)\),
		\[
		\begin{aligned}
			f'&=(\sin t)^{-h}(y'-h\cot t\,y),\\
			f''&=(\sin t)^{-h}\left[y''-2h\cot t\,y'
			+\left(\frac{h}{\sin^2t}+h^2\cot^2t\right)y\right].
		\end{aligned}
		\]
		The coefficient of \(y\) in the transformed linear part is
		\[
		\begin{aligned}
			\frac{h}{\sin^2t}+h^2\cot^2t
			-h\frac{\rho'}{\rho}\cot t-{\frac{n(n-2)}4}
			&=\frac{h}{\sin^2t}+(h^2-h(k-1))\cot^2t\\
			&\quad+{h(k-1)-\frac{n(n-2)}4}\\
			&=\frac{h}{2\sin^2t}.
		\end{aligned}
		\]
		Here we used \(h^2-h(k-1)=-h/2\) and
		{\(h(k-1)-n(n-2)/4=-h/2\)}. Also \(h(p_0-2)=2\).
		The equation for \(f\) therefore gives
		\[
		y''+\left(\frac{\rho'}{\rho}-2h\cot t\right)y'
		+\frac{(h/2)y+B_0y^{p_0-1}}{\sin^2 t}=0 .
		\]
		The function \(y\) is positive on this half and vanishes at both
		ends. Since \(f(0)>0\), \(f'(0)=0\), and \(h>0\),
		\[
		y'(t)=h f(0)t^{h-1}+o(t^{h-1})>0
		\quad\text{for sufficiently small }t>0.
		\]
		At the other endpoint,
		\(y'(t_*)=(\sin t_*)^h f'(t_*)<0\).
		Thus \(y'\) has a zero. At each zero the transformed equation gives
		\[
		y''=-\frac{(h/2)y+B_0y^{p_0-1}}{\sin^2t}<0.
		\]
		Every zero of \(y'\) is consequently a strict crossing from positive
		to negative. A second zero would require a crossing in the opposite
		direction in between. Hence there is exactly one stationary point.
		Since \(g_{1,0}=-(\sin t)^{1-h}y'\), it has exactly one zero there,
		and its derivative at that point is
		\[
		g_{1,0}'=-(\sin t)^{1-h}y''>0.
		\]
		The zero is therefore simple.
		For \(t_*<t<\pi/2\) both \(f<0\) and \(f'<0\), so \(g_{1,0}>0\)
		there. It is also positive at \(t_*\).
		
		Let \(t_1\) be the zero of \(g_{1,0}\).
		On a closed interval containing \(t_1\) in its interior and contained
		in \((0,\pi/2)\), the weight \(\rho\) is bounded above and below by
		positive constants. The one-dimensional trace inequality therefore gives
		\[
		|\phi(t_1)|\leq C\|\phi\|_{\mathcal X_{1,0}}.
		\]
		This continuous linear functional is nonzero on smooth profiles
		supported near \(t_1\). Its kernel is a closed subspace of codimension one.
		On this subspace,
		the form \(Q_{1,0}^{(0)}\) is nonnegative. Apply the previous
		ground state identity separately on each side of \(t_1\), using
		\(|g_{1,0}|\). To pass to the form domain, cut off near \(t_1\).
		For \(\phi(t_1)=0\), the derivative cost of this cutoff is bounded by
		\[
		C\varepsilon^{-2}\int_{|t-t_1|<2\varepsilon}|\phi(t)|^2\,dt
		\leq C\int_{|t-t_1|<2\varepsilon}|\phi'(t)|^2\,dt\longrightarrow0 .
		\]
		The endpoint cutoff argument in Section~\ref{sec:angular-domains}
		then gives the required density on both subintervals.
		Thus \(m_{1,0}^{(0)}\leq1\).
		The mixed kernel profile belongs to
		\(\mathcal X_{1,1}\subset\mathcal X_{1,0}\), because the
		\((1,1)\) local form contains all the nonnegative terms defining
		\(\mathcal X_{1,0}\). Thus \(f'\) is an admissible test in \((1,0)\),
		and the two local forms differ by the angular potential
		\((k-1)/\sin^2 t\). Consequently,
		\[
		Q_{1,0}^{(0)}(f')=Q_{1,1}^{(0)}(f')
		-(k-1)\int_0^{\pi/2}\frac{|f'|^2}{\sin^2t}\rho\,dt<0 ,
		\]
		so \(m_{1,0}^{(0)}=1\).
		
		Let \(\phi\) be another kernel profile in \((1,0)\). Near zero \(g_{1,0}\) is bounded away from zero,
		and the local equations give
		\(\rho g_{1,0}^2(\phi/g_{1,0})'=C\) for a constant \(C\).
		The derivative of this quotient lies in \(L^2(\rho\,dt)\) near
		zero, since \(g_{1,0}\) is \(C^1\) and \(\phi\in\mathcal X_{1,0}\).
		If \(C\ne0\), its squared integral is infinite because
		\(\int_0^\varepsilon\rho^{-1}\,dt=\infty\) for \(k\geq2\).
		Hence \(\phi\) is proportional to \(g_{1,0}\) near zero,
		and ODE uniqueness extends this identity across the whole interval.
		Reflection about \(t_*\) identifies the \((0,1)\) and \((1,0)\)
		forms, so their indices and nullities agree.
		
		\medskip
		\noindent\textbf{Persistence for small positive \(\mu\).}
		Along the convergent subsequence above, set \(v_i=v_{\mu_i}\)
		and \(U_i=J_{\mu_i}(|v_i|^{p_{\mu_i}})\).
		The choice of \(r\) and the power convergence show
		\(U_i\to B_0\) uniformly. Also,
		\[
		U_i|v_i|^{p_{\mu_i}-2}\longrightarrow
		B_0 |v^{(0)}|^{p_0-2}\quad\text{in }L^{n/2},
		\]
		because \((p_0-2)n/2=p_0<R\). Thus the local Hessian term converges
		in bilinear form norm on \(H^1\).
		For the nonlocal term, write \(b_i=|v_i|^{p_{\mu_i}-2}v_i\).
		The powers converge strongly in \(L^q\), where
		\(1/q=1/r-1/2^*>(p_0-1)/R\).
		Put \(b_0=|v^{(0)}|^{p_0-2}v^{(0)}\) and
		\(T_i\phi=b_i\phi\), \(T_0\phi=b_0\phi\). H\"older's inequality
		and Sobolev embedding give
		\[
		\|T_i-T_0\|_{H^1\to L^r}
		\leq C\|b_i-b_0\|_{L^q}\longrightarrow0,
		\qquad \sup_i\|T_i\|_{H^1\to L^r}<\infty.
		\]
		For \(\phi,\psi\in H^1\), decompose the difference of the nonlocal
		bilinear terms as
		\[
		\begin{aligned}
			&\int J_{\mu_i}(T_i\phi)T_i\psi\,dV
			-\int J_0(T_0\phi)T_0\psi\,dV\\
			&=\int (J_{\mu_i}-J_0)(T_i\phi)T_i\psi\,dV
			+\int J_0((T_i-T_0)\phi)T_i\psi\,dV\\
			&\quad+\int J_0(T_0\phi)(T_i-T_0)\psi\,dV.
		\end{aligned}
		\]
		The finite volume of the sphere gives \(\|z\|_1\leq C\|z\|_r\).
		Each term is therefore bounded using \(L^r\to L^\infty\) for
		\(J_{\mu_i}\) and \(J_0\). In particular, the absolute value of
		the displayed difference is at most
		\[
		C\left(\|J_{\mu_i}-J_0\|_{L^r\to L^\infty}
		+\|b_i-b_0\|_{L^q}\right)\|\phi\|\|\psi\|.
		\]
		For the local term, H\"older's inequality gives
		\[
		\begin{aligned}
			&\left|\int\bigl((p_{\mu_i}-1)U_i|v_i|^{p_{\mu_i}-2}
			-(p_0-1)B_0|v^{(0)}|^{p_0-2}\bigr)\phi\psi\,dV\right|\\
			&\quad\leq C\left\|(p_{\mu_i}-1)U_i|v_i|^{p_{\mu_i}-2}
			-(p_0-1)B_0|v^{(0)}|^{p_0-2}\right\|_{L^{n/2}}
			\|\phi\|\|\psi\|.
		\end{aligned}
		\]
		The right-hand sides tend to zero. Since \(p_{\mu_i}\to p_0\),
		the full Hessians converge in energy operator norm.
		
		Restrict this convergence to each of the three radial energy spaces.
		The limiting energy operators are compact perturbations of the identity.
		Their negative spectra are separated from zero, and zero is a simple
		isolated eigenvalue in each block. Spectral projections for the
		negative spectrum away from zero and for a small interval about zero
		have constant rank under sufficiently small norm perturbations
		{\cite[Chapter~IV]{Kato1995}}. The symmetry fields are
		nonzero kernel elements in the three blocks at each \(\mu_i\).
		Thus the one-dimensional spectral subspace near zero is exactly the
		kernel. The negative counts are therefore the limiting ones.
		To obtain a threshold independent of the chosen minimizer, suppose
		that no \(\mu_*(k)>0\) has the asserted property. We could then
		choose \(\mu_i\to0\) and a minimizer at each \(\mu_i\)
		for which at least one of the three counts fails. Change its sign
		so that the profile is positive on the first half. The compactness
		argument above provides a subsequence converging to a minimizer of
		the limit problem. At that limit, each of the three zero eigenspaces
		is one-dimensional and the negative indices have the stated values.
		Norm convergence and the nonzero symmetry fields force all three
		counts along this subsequence to agree with those limiting values
		for large \(i\), a contradiction. Shrinking the threshold to be
		less than \(\min\{4,n\}\) proves \eqref{eq:intro-swap-inertia}
		for every minimizer in the stated parameter range.
		The one-dimensional kernel in each block is generated by its
		nonzero symmetry profile.
	\end{proof}
	\medskip
	
	\noindent\textbf{Acknowledgments.}
	This work was supported by the National Natural Science Foundation of China (No. 12671143).
	
	\medskip
	\noindent\textbf{Data availability statement.}
	No datasets were generated or analyzed.
	
	\medskip
	\noindent\textbf{Competing interests.}
	The authors declare no competing interests.


\begin{thebibliography}{99}
		\raggedright
		
		\bibitem{AlmutairiGluck2026}
		M. Almutairi and M. Gluck,
		{Infinitely many solutions to a conformally invariant elliptic
			equation with Choquard-type nonlinearity},
		Ann. Mat. Pura Appl. 205 (2026), 1541--1575.
		
		\bibitem{Aubin1998}
		T. Aubin,
		Some Nonlinear Problems in Riemannian Geometry,
		Springer Monographs in Mathematics,
		Springer, Berlin, 1998.
		
		\bibitem{BahriLions1988}
		A. Bahri and P.-L. Lions,
		{Morse index of some min--max critical points. I. Application to
			multiplicity results},
		Comm. Pure Appl. Math. 41 (1988), 1027--1037.
		
		\bibitem{Beckner1993}
		W. Beckner,
		{Sharp Sobolev inequalities on the sphere and the
			Moser--Trudinger inequality},
		Ann. of Math. 138 (1993), 213--242.
		
		\bibitem{BianchiEgnell1991}
		G. Bianchi and H. Egnell,
		{A note on the Sobolev inequality},
		J. Funct. Anal. 100 (1991), 18--24.
		
		\bibitem{CaffarelliGidasSpruck1989}
		L. A. Caffarelli, B. Gidas and J. Spruck,
		{Asymptotic symmetry and local behavior of semilinear elliptic
			equations with critical Sobolev growth},
		Comm. Pure Appl. Math. 42 (1989), 271--297.
		
		\bibitem{Chang1993}
		K.-C. Chang,
		Infinite Dimensional Morse Theory and Multiple Solution Problems,
		vol.~6, Birkh\"auser Boston Inc, Boston, MA, 1993.
		
		\bibitem{ChangGhoussoub1996}
		K.-C. Chang and N. Ghoussoub,
		{The Conley index and the critical groups via an extension of
			Gromoll--Meyer theory},
		Topol. Methods Nonlinear Anal. 7 (1996), 77--93.
		
		\bibitem{CoddingtonLevinson1955}
		E. A. Coddington and N. Levinson,
		Theory of Ordinary Differential Equations,
		McGraw--Hill, New York, 1955.
		
		\bibitem{DaiLiuQin2021}
		W. Dai, Z. Liu and G. Qin,
		{Classification of nonnegative solutions to static
			Schr\"odinger--Hartree--Maxwell type equations},
		SIAM J. Math. Anal. 53 (2021), 1379--1410.
		
		\bibitem{DaiXu2013}
		F. Dai and Y. Xu,
		Approximation Theory and Harmonic Analysis on Spheres and Balls,
		Springer Monographs in Mathematics,
		Springer, New York, 2013.
		
		\bibitem{Ding1986}
		W. Y. Ding,
		{On a conformally invariant elliptic equation on
			$\mathbb R^n$},
		Comm. Math. Phys. 107 (1986), 331--335.
		
		\bibitem{DLMF}
		F. W. J. Olver et al. (eds.),
		NIST Digital Library of Mathematical Functions,
		Section~5.11, National Institute of Standards and Technology,
		\href{https://dlmf.nist.gov/5.11}{dlmf.nist.gov/5.11}.
		
		\bibitem{DuYang2019}
		L. Du and M. Yang,
		{Uniqueness and nondegeneracy of solutions for a critical
			nonlocal equation},
		Discrete Contin. Dyn. Syst. 39 (2019), 5847--5866.
		
		\bibitem{FadellRabinowitz1978}
		E. R. Fadell and P. H. Rabinowitz,
		{Generalized cohomological index theories for Lie group actions
			with an application to bifurcation questions for Hamiltonian systems},
		Invent. Math. 45 (1978), 139--174.
		
		\bibitem{Fang1995}
		G. Fang,
		{Morse indices of critical manifolds generated by min--max methods
			with compact Lie group actions and applications},
		Comm. Pure Appl. Math. 48 (1995), 1343--1368.
		
		\bibitem{FrankLieb2012}
		R. L. Frank and E. H. Lieb,
		{A new, rearrangement-free proof of the sharp
			Hardy--Littlewood--Sobolev inequality},
		in B. M. Brown, J. Lang and I. G. Wood (eds.),
		Spectral Theory, Function Spaces and Inequalities,
		vol.~219, Springer, Basel, 2012, pp.~55--67.
		
		\bibitem{GhimentiMorozVanSchaftingen2017}
		M. Ghimenti, V. Moroz and J. Van Schaftingen,
		{Least action nodal solutions for the quadratic Choquard equation},
		Proc. Amer. Math. Soc. 145 (2017), 737--747.
		
		\bibitem{GhimentiVanSchaftingen2016}
		M. Ghimenti and J. Van Schaftingen,
		{Nodal solutions for the Choquard equation},
		J. Funct. Anal. 271 (2016), 107--135.
		
		\bibitem{GiacomoniWeiYang2020}
		J. Giacomoni, Y. Wei and M. Yang,
		{Nondegeneracy of solutions for a critical Hartree equation},
		Nonlinear Anal. 199 (2020), Article 111969.
		
		\bibitem{GilbargTrudinger1977}
		D. Gilbarg and N. S. Trudinger,
		Elliptic Partial Differential Equations of Second Order,
		1st ed, vol.~224, Springer, Berlin, Heidelberg and New York, 1977.
		
		\bibitem{GuoHuPengShuai2019}
		L. Guo, T. Hu, S. Peng and W. Shuai,
		{Existence and uniqueness of solutions for Choquard equation
			involving Hardy--Littlewood--Sobolev critical exponent},
		Calc. Var. Partial Differential Equations 58 (2019), Article 128.
		
		\bibitem{HebeyVaugon1997}
		E. Hebey and M. Vaugon,
		{Sobolev spaces in the presence of symmetries},
		J. Math. Pures Appl. 76 (1997), 859--881.
		
		\bibitem{Kato1995}
		T. Kato,
		Perturbation Theory for Linear Operators,
		2nd ed, Springer, Berlin, 1995.
		
		\bibitem{LeeParker1987}
		J. M. Lee and T. H. Parker,
		{The Yamabe problem},
		Bull. Amer. Math. Soc. (N.S.) 17 (1987), 37--91.
		
		\bibitem{Lieb1977}
		E. H. Lieb,
		{Existence and uniqueness of the minimizing solution of
			Choquard's nonlinear equation},
		Stud. Appl. Math. 57 (1977), 93--105.
		
		\bibitem{Lieb1983}
		E. H. Lieb,
		{Sharp constants in the Hardy--Littlewood--Sobolev and related
			inequalities},
		Ann. of Math. 118 (1983), 349--374.
		
		\bibitem{LiLiuTangXu2026}
		X. Li, C. Liu, X. Tang and G. Xu,
		{The nondegeneracy of positive bubble solutions for generalized
			energy-critical Hartree equations},
		Sci. China Math. (2026).
		
		\bibitem{MorozVanSchaftingen2017}
		V. Moroz and J. Van Schaftingen,
		{A guide to the Choquard equation},
		J. Fixed Point Theory Appl. 19 (2017), 773--813.
		
		\bibitem{Palais1979}
		R. S. Palais,
		{The principle of symmetric criticality},
		Comm. Math. Phys. 69 (1979), 19--30.
		
		\bibitem{Spanier1966}
		E. H. Spanier,
		Algebraic Topology,
		McGraw--Hill, New York, 1966.
		
		\bibitem{Szego1975}
		G. Szeg\H{o},
		Orthogonal Polynomials,
		4th ed, vol.~23, American Mathematical Society,
		Providence, RI, 1975.
		
		\bibitem{Willem1996}
		M. Willem,
		Minimax Theorems,
		vol.~24, Birkh\"auser Boston Inc, Boston, MA, 1996.
		
	\end{thebibliography}
\end{document}